\documentclass[12pt]{article}

\usepackage{import}
\usepackage{amsmath, amsfonts, amsthm, amssymb}

\usepackage{mathtools}
\usepackage{hyperref}
\usepackage{graphicx}
\usepackage{easy-todo}
\usepackage{enumitem}
\usepackage[hypcap=false]{caption}
\usepackage{parskip}
\usepackage[noabbrev,nameinlink,capitalise]{cleveref}
\usepackage[dvipsnames]{xcolor}
\definecolor{color0}{RGB}{216,27,96} % Color blind friendly
\definecolor{color1}{RGB}{255,193,7}
\definecolor{color2}{RGB}{30,136,229}    

\usepackage{tikz-cd}
\usepackage{tikz}
\usetikzlibrary{shapes, shapes.geometric, arrows, arrows.meta, decorations.markings, decorations.pathreplacing, automata, positioning}
\pgfdeclarelayer{background} 
\pgfdeclarelayer{foreground}
\pgfsetlayers{background,main,foreground}

\theoremstyle{plain}
\newtheorem{thm}{Theorem}[section]
\newtheorem*{thm*}{Theorem}
\newtheorem{lemma}[thm]{Lemma}
\newtheorem*{lemma*}{Lemma}
\newtheorem{cor}[thm]{Corollary}
\newtheorem*{cor*}{Corollary}
\newtheorem{prop}[thm]{Proposition}
\newtheorem*{prop*}{Proposition}

\theoremstyle{definition}
\newtheorem{defi}[thm]{Definition}
\newtheorem*{defi*}{Definition}
\newtheorem{rem}[thm]{Remark}
\newtheorem{exa}[thm]{Example}
\newtheorem*{conj*}{Conjecture}

\newtheorem{innercustomthm}{Theorem}
\newenvironment{customthm}[1]
{\renewcommand\theinnercustomthm{#1}\innercustomthm}
{\endinnercustomthm}

\newtheorem{innercustomprop}{Proposition}
\newenvironment{customprop}[1]
{\renewcommand\theinnercustomprop{#1}\innercustomprop}
{\endinnercustomprop}

\newtheorem{innercustomcor}{Corollary}
\newenvironment{customcor}[1]
{\renewcommand\theinnercustomcor{#1}\innercustomcor}
{\endinnercustomcor}

\newcommand{\say}[1]{``#1"}

\DeclareMathOperator{\id}{id}
\newcommand{\vardot}{\,\cdot\,}

\newcommand{\Rc}{\mathcal R} % Right comb (left comb is already defined)
\newcommand{\Sc}{\mathcal S} % Slalom
\newcommand{\Tc}{\mathcal T} % Trees
\newcommand{\Uc}{\mathcal U}
\newcommand{\Vc}{\mathcal V}
\newcommand{\Sym}{\mathrm{Sym}}
\newcommand{\FSym}{\mathrm{FSym}}
\newcommand{\Homeo}{\mathrm{Homeo}}
\newcommand{\Zc}{\mathcal Z} % Center

\newcommand{\Lc}{\mathcal L}

\newcommand{\coWP}{\mathrm{coWP}}

\DeclareMathOperator{\dist}{Dist}
\newcommand{\supp}{\mathrm{supp}}
\DeclareMathOperator{\rot}{rot}

\DeclareMathOperator{\lcm}{lcm}
\DeclareMathOperator{\ord}{ord}

\newcommand{\E}{\mathbb E}
\DeclareMathOperator{\LDS}{LDS}

\newcommand{\norm}[1]{\left\|#1\right\|}
\newcommand{\abs}[1]{\left|#1\right|}
\newcommand{\la}{\left\langle}
\newcommand{\ra}{\right\rangle}

\renewcommand{\ge}{\geqslant}
\renewcommand{\le}{\leqslant}

\newcommand{\longto}{\longrightarrow}
\newcommand{\onto}{\twoheadrightarrow}
\newcommand{\into}{\hookrightarrow}
\newcommand{\acts}{\curvearrowright}

\newcommand{\N}{\mathbb N}
\newcommand{\Z}{\mathbb Z}
\newcommand{\Q}{\mathbb Q}
\newcommand{\R}{\mathbb R}
\newcommand{\F}{\mathbb F} % finite field
\renewcommand{\P}{\mathbb P} % primes
\newcommand{\bS}{\mathbb S} % circle

\newcommand{\Vmerge}{\begin{tikzpicture}[scale=.3, baseline=1pt]
	\draw[line width=.5pt] (.5,0) -- (.5,.5) -- (0,1);
    \draw[line width=.5pt] (.5,.5) -- (1,1);
	\end{tikzpicture}}

\newcommand{\Vsplit}{\begin{tikzpicture}[scale=.3, baseline=1pt]
	\draw[line width=.5pt] (.5,1) -- (.5,.5) -- (0,0);
    \draw[line width=.5pt] (.5,.5) -- (1,0);
	\end{tikzpicture}}

\newcommand{\Address}{{% additional braces for segregating \footnotesize
		\bigskip
		\small
		
		\textsc{
        University of Newcastle, Australia}
		\par\nopagebreak
        \textit{E-mail address}: \texttt{alexbishop1234@gmail.com} 

        \bigskip
		\par\nopagebreak
		\textsc{%Mathematical Institute,
        University of Oxford, United Kingdom}
        \par\nopagebreak
		\textit{E-mail address}: \texttt{corentin.bodart@maths.ox.ac.uk} 
        
		\bigskip
		\par\nopagebreak
		\textsc{%Section de Mathématiques,
        Université de Genève, Switzerland}
        \par\nopagebreak
        \textit{E-mail address}: \texttt{letizia.issini@unige.ch}, \texttt{davide.perego@unige.ch}
}}

\title{Period growth and co-context-free groups}
\author{Alex Bishop, Corentin Bodart, Letizia Issini, Davide Perego}
\date{\today}

\usepackage[backend=biber, style=numeric, sorting=nyt, %sorting=none,
    doi=false,isbn=false,url=false,eprint=false, maxnames=50]{biblatex}
\begin{document}

\maketitle

\begin{abstract}
    We study period growth in co-context-free groups, giving general results and looking at specific examples such as Thompson groups $T$ and $V$ and the Houghton groups $H_m$. Along the way, we give a refined upper bound on the word metric in Thompson $V$, as well as efficient algorithms to determine if elements of $V$ are torsion, and compute their order. We also adapt our algorithm to compute the rotation number of elements of $T$ and answer a question of D.~Calegari.
    % We also study distortion of cyclic subgroups in co-context-free groups.
\end{abstract}

\section*{Introduction}

Geometric group theory is full of invariants obstructing embeddings between groups. The most striking examples are volume growth, subgroup distortion, F\o{}lner function, Poincaré profiles, etc. In the present work, we study another invariant of a slightly different nature: period growth.

%\begin{defi}
%    Let $Q\subseteq\P$ a set of primes.
%    \begin{itemize}[leftmargin=8mm]
%        \item A positive integer $n$ is $Q$-smooth if all its prime divisors belong to $Q$.
%        \item A group element $g$ is $Q$-torsion if its order $\ord(g)$ is finite and $Q$-smooth.
%    \end{itemize}
%\end{defi}

\medskip

\begin{defi*}[Period growth]
    Let $G$ be a group and $S$ a finite generating set. Consider a subset $D\subseteq \N$ which is closed under divisor, i.e., if $a\mid b$ and $b\in D$, then $a\in D$. We define the period growth as
    \[ p^D_{(G,S)}(n) = \max\left(\left\{\ord(g) \;\big|\; g\in G\text{ and }\norm g_S\le n\right\}\cap D\right). \]
    If $D=\N$, we simply write $p_{(G,S)}(n)$.
\end{defi*}
A function close to $p_{(G,S)}(n)$ was introduced in \cite{grigorchuk1983milnor}, and further studied in \cite{Bartholdi_Sunic, Petschick_2023, bradford2024quantifying}. However, they define $p_{(G,S)}(n)=\infty$ as soon as the ball $B_{(G,S)}(n)$ contains an infinite order element, hence restrict to torsion groups.
%\textbf{Standing assumption.} We will always assume that $D$ is closed under divisor, i.e., if $a\mid b$ and $b\in D$, then $a\in D$.
\newpage

Apart from $D=\N$, typical examples are 
    \begin{itemize}[leftmargin=7mm]
        \item $D=\P\coloneqq\{1\}\cup\{\text{primes}\}$;
        \item $D_Q \coloneqq \left\{m\in\N \;\big|\; \{\text{prime divisors of }m\}\subseteq Q \right\}$ for $Q$ a set of primes;
        \item $D_k \coloneqq \left\{m\in\N \;\big|\; \#\{\text{prime divisors of }m\}\le k \right\}$ for $k$ a positive integer.
        % $\{n \mid \omega(n)\leq k\}$ where $\omega(n)$ is the number of prime divisors of $n$ (without multiplicities)
    \end{itemize}
\medskip

An easy observation (\cref{lem:sub-over-group}) is that the growth type of $p_{(G,S)}^D(n)$ is a group invariant and is monotone under embedding. We relate this with a conjecture that has sparked significant work over the last fifteen years:
\begin{conj*}[Lehnert \cite{lehnert, almost_automorphism}]
    A finitely generated group $G$ is co-context-free if and only if it embeds in Thompson group $V$.
\end{conj*}
In particular, studying period growth for co-context-free groups and Thompson group $V$ could lead to a disproof of this conjecture.
\bigskip

Another important observation links the period growth with the Torsion Problem, i.e., deciding if a word represents an element of finite order or not:
\begin{lemma*} \label{remark:torsion problem}
If $G$ has decidable Word Problem, then $G$ has decidable Torsion Problem if and only if $p_{(G,S)}(n)$ is bounded by a recursive function.
\end{lemma*}
For instance, the period growth of Brin-Thompson group $2V$ is not bounded by any recursive function \cite{torsion_problem_2V}. This contrasts with Thompson group $V$ \cite[Theorem 9.3]{higman1974finitely} and co-context-free groups \cite[Theorem 9]{co-CF} where the Torsion Problem \emph{is} decidable, hence an explicit bound might be achievable.
% Here is a link for Higman's lecture notes: 
% https://www.normalesup.org/~cornulier/Higman_Camberra_1974.pdf
\medskip

These observations motivate the study of period growth within co-context-free groups, and on specific examples. We concentrate on the Houghton groups $H_m$ and Thompson groups $T$ and $V$. (Thompson $F$ is torsionfree.) %, so $p_F(n)\equiv 1$.)
We summarise our main results in the following theorem:
\begin{thm*} Let $G$ be a co-context-free group and $k\in\N$.
    \begin{itemize}[leftmargin=7mm]
        \item $p_G(n)\preceq \exp(n^4)$ and $p_G^{D_k}(n)\preceq \exp(n^2)$. \hfill {\normalfont(\cref{thm:coCF_period})}
        \item $p_{H_m}(n)\asymp\exp(\sqrt{n\log n}\bigr)$ and $p_{H_m}^{D_k}(n)\asymp n^k$.
        \hfill {\normalfont(\cref{thm:period_H})}
        \item $p_T(n)\asymp\exp(n)$.
        \hfill {\normalfont(\cref{thm:period_T})}
        \item $p_V(n)\asymp\exp(n^2)$ and $p_V^{D_k}(n)\asymp\exp(n)$.
        \hfill {\normalfont(\cref{thm:period_V})}
    \end{itemize}
\end{thm*}
Slightly more general results are proven, see the different statements in the main body of text.
%\begin{customthm}{\ref{thm:coCF_period}}
%    Let $G$ be a co-context-free group.
%    \begin{enumerate}[leftmargin=8mm, label={\normalfont(\alph*)}]
%        \item The period growth satisfies $p_G(n)\preceq \exp(n^4)$. 
%        \item If $D\subseteq D_k$ for some integer $k$, then $p_G^D(n)\preceq \exp(n^2)$. 
%    \end{enumerate}
%\end{customthm}
%
%\begin{customthms}{\ref{thm:period_T} and \ref{thm:period_V}}
%    \begin{itemize}[leftmargin=7mm]
%        \item If $D\subseteq D_k$ satisfies a\say{logarithmic Bertrand postulate}, i.e., $\exists C\ge 1$, $\forall M\in \N$, $\exists m\in D$ such that $M\le m\le M^C$, then $p_V^D(n)\asymp \exp(n)$.
%        \item $p_T(n)\asymp \exp(n)$.
%        \item $p_V(n)\asymp \exp(n^2)$.
%    \end{itemize}
%\end{customthms}
In order to prove these results on Thompson $T$ and $V$, we need two main ingredients of independent interest: %good control on the word length of elements, and an efficient way to compute the order of elements.

$\triangleright$ In \cref{sec:metric}, we give an upper bound on the length of elements of $V$ in terms of the size of their diagrams:
\begin{customthm}{\ref{thm:upper_bound_LDS}}
    Let $g=(\hspace{1pt}\Uc,\sigma,\Vc)\in V$ be an element where $\Uc,\Vc$ are two trees with $n$ leaves and $\sigma\in\Sym(n)$. Then
    \[ \norm{g}_V \preceq n\cdot \log(\LDS(\sigma))\]
    where $\LDS(\sigma)=\max\bigl\{\abs A: \sigma|_A\text{ is decreasing}\bigr\}$.
\end{customthm}
Since $1\le \LDS(\sigma)\le n$, this gives an alternative proof of \cite{ThompsonV_metric} with a finer gradation, and improves on a recent result \cite{burillo_note}. Note that this gives a linear upper bound for an exponentially large (but necessarily exponentially negligible) set of permutations $\sigma\in\Sym(n)$.

$\triangleright$ In \cref{sec5_diagrams}, we describe how to compute efficiently the order of elements in $T$ and $V$. Our algorithm uses a variation on abstract strand diagrams, as introduced by Belk and Matucci \cite{strand_diagrams}. Adding some additional labelling, we can also compute the rotation number of elements of $T$.
\medskip

It should be noted that a method to compute and bound the denominator of the rotation number of elements of $T$ was described by D.~Calegari in \cite{calegari2007denominator}. However, he measures the complexity of elements in terms of the depth of the trees in a reduced tree diagram, and not their word length. In particular, this only implies a bound of $p_T(n)\preceq\exp(\exp(n))$. He further proposes better bounds in terms of the depth which we disprove in \cref{subsection:Calegari}.
\bigskip

Finally, we look at \textbf{CF-TR} groups, a class of groups recently introduced in \cite{CFTR}. An important question in regard of Lehnert's conjecture is whether Thompson $V$ is \textbf{CF-TR}. Indeed, a positive answer would reduce Lehnert's conjecture to two seemingly simpler statements: \say{Are all co-context-free groups \textbf{CF-TR}?} and \say{Do all \textbf{CF-TR} groups embed in $V$?}. Combining our result on period growth of $V$ with previous results on period growth of \textbf{CF-TR} groups \cite[Section 5]{CFTR}, we get the following restriction:
\begin{customcor}{\ref{cor:CF_graph_for_V_questionmark}}
    If $\Gamma$ is a context-free graph whose transition group is $V$, then the volume growth $\bar\beta_\Gamma(n)\coloneqq\sup_{v\in \Gamma}\abs{B_\Gamma(v,n)}$ is exponential.
\end{customcor}
Any improvement on \cite[Section 5]{CFTR} would imply that $V$ is not \textbf{CF-TR}. However, we prove that at least one of their estimates is sharp:
\begin{customprop}{\ref{prop:BEN_has_large_period_growth}}
    There exists a \textbf{CF-TR} group $G$ with $p_G^{D_{\{2\}}}(n)\asymp \exp(n)$.
\end{customprop}
We obtain this result by considering a family of groups generalising the Brin-Navas group. These groups were introduced by A.~Erschler \cite{Erschler}.
\section{Generalities on period growth}
\begin{lemma} \label{lem:sub-over-group}
Let $(G,S)$ and $(H,T)$ be two finitely generated groups.
\begin{enumerate}[leftmargin=8mm, label={\normalfont(\alph*)},          ref=\thelemma(\alph*)]
    \item \label{lem:subgroup}  If $G\ge H$, then there exists $C>0$ such that
    \[ \forall D, \forall n\in\N, \quad p_{(G,S)}^D(Cn) \ge p_{(H,T)}^D(n). \]
    \item \label{lem:fi_overgroup} If $G\ge H$ and $H$ has finite-index, then there exists $C>0$ such that
    \[ \forall D,\forall n\in\N, \quad p_{(G,S)}^D(n) \le C\cdot p_{(H,T)}^D(Cn). \]
\end{enumerate}
In particular, the function $p_G^D(n)$ is well-defined up to $\asymp$.
\end{lemma}
\begin{proof}
Part (a) follows easily by taking $C:=\displaystyle\max_{t \in T} \norm t_S$. Indeed, one just has to notice that for any $g \in H$ it holds $\norm g_S \leq C\norm g_T$ and hence $$\{g \in H \mid \ord(g) \in D, \ \norm g_T \leq n \} \subseteq \{g \in G \mid \ord(g) \in D, \ \norm g_S \leq Cn \}.$$ 

Part (b) has the further assumption that $H$ is finite index. It follows that there exist $r,R>0$ such that for all $g \in G$ one has $g^r \in H$ and for all $h \in H$ it holds $\norm h_T \leq R \norm h_S$. Combining the two, $\norm {g^r}_T \leq R \norm {g^r}_S \leq Rr \norm g_S$ and set $C=Rr$. Moreover, $D$ is closed under divisor and hence $\ord(g^r) \in D$. To conclude,
\begin{equation*}
\begin{split}
\max & \{ \ord(g) \;\big|\; \ord(g)\in D \text{ and }\norm g_S\le n \}  \\
& \le \max\left\{r\ord(g^r) \;\big|\; \ord(g^r)\in D, \ g^r \in H \text{ and }\norm {g^r}_T\le Cn\right\} \\
& \le C\max\left\{\ord(h) \;\big|\; \ord(h)\in D, \ h \in H \text{ and }\norm {h}_T\le Cn\right\}. \qedhere
\end{split}
\end{equation*}
\end{proof}

\smallskip
\begin{rem}
    Keeping in mind that $C$ doesn't depend on $D$ might give a more powerful invariant. However, we will not explore this direction.
\end{rem}
\bigskip
Let us briefly investigate the behavior of the period growth with respect to extensions.
\begin{lemma} \label{lem:extension}
Consider the group extension $1 \rightarrow K \rightarrow G \xrightarrow{\pi} H \rightarrow 1$. Unless stated otherwise, we assume $K$ is finitely generated. Consider the (upper) distortion function $\dist_{K}^{G}(n) \coloneqq \max\bigl\{\norm{k}_K \mid k\in K,\; \norm{k}_G\le n\bigr\}$.
    \begin{enumerate}[leftmargin=8mm, label={\normalfont(\alph*)}, ref=\thelemma(\alph*)]
        \item \label{item-semi:1} We have $p_{G}^D(n) \le p_H^D(n)\cdot p_K^D\!\left( \dist_K^{G}(n\cdot p_H^D(n))\right)$. If $\{\ord(k) \mid k \in K\} \cap D$ is finite, then $p_{G}^D(n) \preceq p_H^D(n)$ holds even if $K$ is not finitely generated.
        \item \label{item-semi:2} If it splits, then $p^{D}_{G}(n)\succeq \max\{p^{D}_{K}(n), p^{D}_{H}(n)\}$.
        \item \label{item-semi:3} If it is a wreath-product with $K= \bigoplus_{H} L$, then
        \[ \max\{p^{D}_{H}(n),p^{D}_{L}(n)\} \preceq p_G^D(n) \preceq p_H^D(n)\cdot p_L^D(n)^n. \]
        If moreover $D$ is closed under multiplication, then $p^{D}_{G}(n)\succeq~p^{D}_{H}(n) p^{D}_{L}(n)$.
    \end{enumerate}
    Moreover, if the extension is trivial (i.e.\ $G=K\times H$), then
    \begin{enumerate}[leftmargin=8mm, label={\normalfont(\alph*)}]
    \setcounter{enumi}{3}
        \item \label{item-direct:1} $p^{D}_{G}(n)\preceq p^{D}_{K}(n) \cdot p^{D}_{H}(n)$.
        \item \label{item-direct:2} If $D,D'\subseteq \N$ and $D\cap D'=\{1\}$, then $p^{DD'}_{G}(n)\succeq p^{D}_{K}(n) \cdot p^{D'}_{H}(n)$.
        \item \label{item-direct:3} If $D\subseteq D_1$, then $p^{D}_{G}(n) \asymp \max\{p^{D}_{K}(n), p^{D}_{H}(n)\}$.
    \end{enumerate}
\end{lemma}
\begin{proof}
In order to show \ref{item-semi:1}, let $x\in G$ be such that $\ord(g) \in D$ and $\norm{x}_{G}\leq n$. We have $\ell=\ord(\pi(x))\leq p_H^{D}(n)$ because $D$ is closed under divisor and $\norm{\pi(x)}_H\le\norm{x}_G$. Moreover $x^\ell\in K$ hence
\begin{equation*}
\ord(x)= \ell\ord(x^\ell)\le \ell\cdot p_K^{D}\left(\dist_K^{G}(\norm{x^\ell}_{G})\right)\leq p_H^{D}(n) \cdot p_K^{D}(\dist_K^{G}(n\cdot p_H^{D}(n))).
\end{equation*}
If $M=\sup \bigl(\{\ord(k) \mid k \in K\} \cap D\bigr)$ is finite, then we can bound $\ord(x^\ell)\le M$ so that $\ord(x) = \ell\ord(x^\ell) \le \ell\cdot M \preceq p_H^D(n)$. \\

Part \ref{item-semi:2} follows from Lemma \ref{lem:subgroup} since $K$ and $H$ are subgroups of $G$. \\

For part \ref{item-semi:3}, we have that $H$ and $L$ are subgroups of $G=L\wr H$, so the first lower bound follows directly from Lemma \ref{lem:subgroup}.

\newpage

For the upper bound, consider $g=(\varphi,h)$ with $\varphi\colon H\to L$ (finitely supported) and $h\in H$. We suppose that $\norm g_G\le n$ and $\ord(g)\in D$, hence the order $\ell:=\ord(h)\in D$. More precisely we have $\ord(g)=\ell\cdot \ord(g^\ell)$ where
\[ g^\ell = ( \psi, 1_H) \quad \text{with}\quad \psi(x) = \varphi(x) \varphi(h^{-1}x) \ldots \varphi(h^{-\ell+1}x).\]
We first observe that, if $x=h^ny$ (with $0\le n\le \ell-1$), then
\[ \psi(x)\cdot \bigl(\varphi(x)\ldots \varphi(h^{-n+1}x)\bigr) = \bigl(\varphi(x)\ldots \varphi(h^{-n+1}x)\bigr) \cdot \psi(y). \]
We deduce that non-trivial values of $\psi(x)$ come in at most $\abs{\supp(\varphi)}\le\norm{g}_G$ conjugacy classes, hence $\ord(\psi(x))$ takes at most $\norm{g}_G$ values. Moreover, we have $\norm{\psi(x)}_L\le \sum_{y\in H}\norm{\varphi(y)}_L\le\norm{g}_G$, therefore $\ord(\psi(x))\preceq p_L^D(n)$. We conclude that $\ord(g^\ell)=\lcm\{\ord(\psi(x)):x\in H\}\le p_L^D(n)^n$, hence
\[ \ord(g) \le p_H^D(n)\cdot p_L^D(n)^n. \]
In the case where $D$ is closed under multiplication, let $h\in H$ and $l\in L$ be of length a most $n$ with orders $p_H^D(n)$ and $p_L^D(n)$ respectively. Then $hl$ has order $p_H^D(n)p_L^D(n)$. Since $\norm{hl}_{G}\leq 2n$, we have $p_G^D(2n)\geq p_H^D(n)\cdot p_L^D(n)$.\\

% \begin{equation*}
%  \max\left\{\ord((\delta_{1}^{l},1)) \;\big|\; (\delta_{1}^{l},1) \text{ is $Q$-torsion and }\norm l_{L}\le n \right\}
% \end{equation*}
% where $\delta_{1}^{l}: H\to L$ is the function $1\mapsto l, h\mapsto 1$ for all $h\neq 1$. This is because $\delta_{1}^{l}$ is torsion if and only if $l$ is torsion and because $\norm l_{L}= \norm{(\delta_{1}^{l},1)}_{G \rtimes H}$.

For \ref{item-direct:1}, 
\begin{equation*}
\begin{split}
& \max \bigl\{\ord((k,h))  \;\big|\;  \ord((k,h)) \in D \text{ and } \norm{(k,h)}_{G}\le n \bigr\}  \\
& \;\le \max\left\{\lcm\{\ord(k),\ord(h)\} \;\big|\; \ord(k),\ord(h) \in D \text{ and }\norm{(k,h)}_{G}\le n \right\} \\
& \;\preceq \max\left\{\ord(k)\cdot \ord(h) \;\big|\; \ord(k),\ord(h) \in D \text{ and }\norm{k}_{K}, \norm{h}_{H}\le n \right\} \\
& \;= p^{D}_{K}(n) \cdot p^{D}_{H}(n).
\end{split}
\end{equation*}

To achieve \ref{item-direct:2}, we consider $k \in K$ and $h \in H$ such that $\norm k_K \leq n$ and $\norm h_H \leq n$ that realize $p_{K}^{D}(n)$ and $p_{H}^{D'}(n)$ respectively. Then
$$\ord((k,h)) =\lcm\bigl(\ord(k),\ord(h)\bigr) \overset!=\ord(k)\ord(h)= p_{K}^{D}(n) \cdot p_{H}^{D'}(n)$$
since $D\cap D'=\{1\}$, and $\norm{(k,h)}_{G} \leq 2n$.
\medskip

For \ref{item-direct:3}, we specialize once again what we wrote in part \ref{item-direct:1}. If $\ord((k,h))=\lcm(\ord(k),\ord(h))\in D\subseteq D_1$, then $\ord(k)$ and $\ord(h)$ are powers of the same prime. It follows that $\ord((k,h))=\max\{\ord(k),\ord(h)\}$ and
\begin{equation*}
\begin{split}
& \max \bigl\{\ord((k,h))  \;\big|\; \ord((k,h))\in D \text{ and }\norm{(k,h)}_{G}\le n \bigr\} \\
& \; \preceq \max\left\{\max\{\ord(k),\ord(h)\} \;\big|\; \ord(k),\ord(h)\in D \text{ and }\norm{k}_{K}, \norm{h}_{H}\le n \right\},
\end{split}
\end{equation*}
i.e., $p_G^D(n)\preceq \max\{p_K^D(n),p_H^D(n)\}$.
\end{proof}

Contrary to volume growth, there are no a priori upper or lower bounds on period growth, for good reasons:
\medskip
\begin{prop} \label{prop:arbitrary_torsion_growth}
Let $f\colon\Z_{>0}\to\Z_{>0}$ be an increasing function. There exists a finitely generated group $G$ with $p_G(n)\asymp f(n)$. Moreover, if $f$ is recursive, then $G$ can be taken with decidable Word and Torsion problems.
\end{prop}
\begin{proof} We provide two different constructions:
\begin{enumerate}[leftmargin=7mm]
    \item Take a sequence of words $w_i=a^ib^i$ and define
    \[ G = \la a,b \;\middle|\; w_i^{9f(i)}=1 , \ \forall i\ra.\]
    First, let us note that $p_G(2n) \geq 9f(n)$ just by considering the elements represented by $w_n$. The presentation satisfies the classical $C'(1/8)$ small cancellation condition, this allows a few observations:
    
    \textbf{Claim (Greendlinger \cite[Theorem VIII]{Greendlinger}).} Let $G=\la S\mid R\ra$ be a $C'(1/6)$ presentation. If an element $g$ has finite order $m>1$, then there exists a relator of the form $r=s^m\in R$ and an integer $p$ such that $g\sim s^p$.
    
    \textbf{Claim.} Let $G=\la S\mid R\ra$ be a $C'(1/8)$ presentation, and $h\in F(S)$ be a cyclically $R$-reduced element, i.e., $h$ and its cyclic permutations do not contain more than $\frac12$ of a relation. Then
    \[ \norm{[h]_\sim}_S := \min\bigl\{\norm g_S \;\big|\; g\sim h\bigr\}\ge \frac12\norm{h}_S. \]
    $\triangleright$ Consider an element $g\sim h$ with minimal length, and an associated reduced annular diagram. Combining Theorems 5.3 and 5.5 of \cite{Lyndon_Schupp} as observed in \cite{Logan_MO}, we deduce that each cell must intersect both boundary component and at most two other cells (see \cref*{fig:conjugacy_diagram}).
    \begin{center}
        \tikzset{every picture/.style={line width=0.75pt}} %set default line width to 0.75pt        
\begin{tikzpicture}[x=0.75pt,y=0.75pt,yscale=-.8,xscale=.8]

\draw[draw=blue, double=white, double distance=8pt, rounded corners]
    (208,73)
    .. controls (229,40) and (293,42) .. (308,73)
    .. controls (351,75) and (366,103) .. (360,120)
    .. controls (378,138) and (352,154) .. (316,166)
    .. controls (321,190) and (263,197) .. (226,185)
    .. controls (201,186) and (185,184) ..  (169,160)
    .. controls (150,150) and (142,126) .. (162,108)
    .. controls (170,87) and (189,75) .. (208,73);
\fill[white]
    (208,73)
    .. controls (229,40) and (293,42) .. (308,73)
    .. controls (351,75) and (366,103) .. (360,120)
    .. controls (378,138) and (352,154) .. (316,166)
    .. controls (321,190) and (263,197) .. (226,185)
    .. controls (201,186) and (185,184) ..  (169,160)
    .. controls (150,150) and (142,126) .. (162,108)
    .. controls (170,87) and (189,75) .. (208,73);

\begin{scope}
\clip (208,73)
    .. controls (240,91) and (282,89) .. (308,73)
    .. controls (351,75) and (366,103) .. (360,120)
    .. controls (342,109) and (329,130) ..(316,166)
    .. controls (303,143) and (246,170) .. (226,185)
    .. controls (201,186) and (185,184) ..  (169,160)
    .. controls (190,135) and (183,118) .. (162,108)
    .. controls (170,87) and (189,75) .. (208,73);
    
\draw[draw=purple, double=white, double distance=8pt, rounded corners]
    (208,73)
    .. controls (240,91) and (282,89) .. (308,73)
    .. controls (351,75) and (366,103) .. (360,120)
    .. controls (342,109) and (329,130) ..(316,166)
    .. controls (303,143) and (246,170) .. (226,185)
    .. controls (201,186) and (185,184) ..  (169,160)
    .. controls (190,135) and (183,118) .. (162,108)
    .. controls (170,87) and (189,75) .. (208,73);
\end{scope}
    
\draw (208,73) .. controls (229,40) and (293,42) .. (308,73); %out
\draw (208,73) .. controls (240,91) and (282,89) .. (308,73); %in

\draw (308,73) .. controls (351,75) and (366,103) .. (360,120);

\draw (360,120) .. controls (378,138) and (352,154) .. (316,166); %out
\draw (360,120) .. controls (342,109) and (329,130) ..(316,166); %in

\draw (316,166) .. controls (321,190) and (263,197) .. (226,185); %out
\draw (316,166) .. controls (303,143) and (246,170) .. (226,185); %in

\draw (226,185) .. controls (201,186) and (185,184) ..  (169,160);

\draw (169,160) .. controls (150,150) and (142,126) .. (162,108); %out
\draw (169,160) .. controls (190,135) and (183,118) .. (162,108); %in

\draw (162,108) .. controls (170,87) and (189,75) .. (208,73);

\draw (230,55) -- (234,83);
\draw (257,49) -- (257,85.6);
\draw (282,52.5) -- (280,83.6);

\draw (332,130.6) -- (350,151.6);

\draw (286,157.6) -- (290,188);
\draw (255,168) -- (258,190);

\draw (150,134.6) -- (181.6,131) ;

\node[circle, fill=black, inner sep=1pt] at (208,73) {};
\node[circle, fill=black, inner sep=1pt] at (308,73) {};
\node[circle, fill=black, inner sep=1pt] at (360,120) {};
\node[circle, fill=black, inner sep=1pt] at (316,166) {};
\node[circle, fill=black, inner sep=1pt] at (226,185) {};
\node[circle, fill=black, inner sep=1pt] at (169,160) {};
\node[circle, fill=black, inner sep=1pt] at (162,108) {};

\node[circle, fill=black, inner sep=1pt] at (230,55)  {};
\node[circle, fill=black, inner sep=1pt] at (234,83) {};
\node[circle, fill=black, inner sep=1pt] at (257,49) {};
\node[circle, fill=black, inner sep=1pt] at (257,85.6) {};
\node[circle, fill=black, inner sep=1pt] at (282,52.6) {};
\node[circle, fill=black, inner sep=1pt] at (280,83.6) {};

\node[circle, fill=black, inner sep=1pt] at (332,130.6) {};
\node[circle, fill=black, inner sep=1pt] at (350,151.6) {};

\node[circle, fill=black, inner sep=1pt] at (286,158) {};
\node[circle, fill=black, inner sep=1pt] at (290,188.6) {};
\node[circle, fill=black, inner sep=1pt] at (255,169) {};
\node[circle, fill=black, inner sep=1pt] at (258,190.6) {};

\node[circle, fill=black, inner sep=1pt] at (150.6,134.6) {};
\node[circle, fill=black, inner sep=1pt] at (181.6,131) {};

\node[purple] at (192,95) {$g$};
\node[blue] at (175,65) {$h$};
\end{tikzpicture}
        \captionsetup{font=small}
        
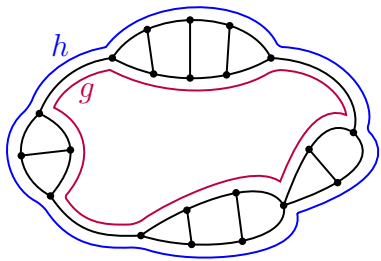
\captionof{figure}{The shape of annular diagrams for $C'(1/8)$ presentations.} \label{fig:conjugacy_diagram}
    \end{center}
    For each cell $C$, we have $\abs{C\cap h}<\frac 12\abs{C}$ while the intersection with each adjacent cell is $<\frac18\abs{C}$, and therefore
    $$ \abs{C\cap g}>\frac14\abs{C} >\frac12\abs{C\cap h}.$$ Summing over all cells, we get $\norm{g}_S\ge\frac12\norm{h}_S$. \hfill$\triangleleft$ 

    Consider a torsion element $g$ of length $n$ and order $\ord(g)>2$. Therefore, there exists $h=w_i^p$ with $\abs{p}<\frac92f(i)$ such $g\sim h$. The second claim gives 
    \[ n = \norm{g}_S\ge \frac12\norm{h}_S \ge i. \]
    Since $\ord(g) \mid 9f(i)$, we conclude that $\ord(g)\le 9f(n)$.
    
    \item Let $g(n)$ be the largest power of $2$ smaller than $f(n)$. We consider the following central extension of $\Z\wr\Z$ introduced 
    %by P.\ Hall
    in \cite{Hall}
    \[ H = \la a,t,(z_i)_{i\ge 1} \;\big|\; z_i=[a,t^iat^{-i}],\; [z_i,a]=[z_i,t]=1 \ra. \]
    This group is center-by-metabelian, hence solvable of class $3$. The quotient $G=H/K$ has the desired period growth, where
    \[ K = \big\langle z_i^{g(i)} : n\ge 1 \big\rangle \le \Zc(H). \]
    The lower bound comes from the elements $z_n$ satisfying $\ord(z_n)=g(n)$ and $\norm{z_n}=4n+4$. For the upper bound, observe that $\Zc(G)\hookrightarrow G\onto \Z\wr\Z$ hence any torsion elements of $G$ must belong to $\Zc(G)$. Moreover,
    \[ B_{G,\{a,t\}}(4n+4)\cap\Zc(G)\le\bigoplus_{i=1}^n\la z_i\ra \simeq \bigoplus_{i=1}^n C_{g(i)},\] and this group has exponent $g(n)$. (We use that $g(i)\mid g(j)$ for $i\le j$.)
\end{enumerate}
Note that, if $f(i)$ is recursive, then the Word Problem is decidable:
\begin{enumerate}[leftmargin=7mm]
    \item $G$ is a small cancellation group, we can use Dehn's algorithm.
    \item Given a word $w\in F(a,t)$, one can decide if $w\in_H\Zc(H)$ (equivalently, $w=1$ in $\Z\wr\Z$). If yes write it as a finite product $w=_H\prod_{i=1}^\ell z_i^{e_i}$ and check if $g(i)\mid e_i$ for all $i$. In this case (and only this case), $w=_G1$. 
\end{enumerate}
Decidability of the Torsion Problem follows from the \hyperref[remark:torsion problem]{Lemma} in the introduction (or more ad hoc methods, eg., just checking if $w=1$ in $\Z\wr\Z$).
\end{proof}
\begin{rem}
    The second construction also allows to construct center-by-metabelian groups with decidable Word Problem and undecidable Torsion Problem. This is done in \cite[Proposition 3.2]{torsion_problem_solv}. %In particular, this group has undecidable Identity Problem, answering a question of Dong \cite{Dong_survey} (end of \S4).
    % For instance, we can quotient by
    %\[ K = \la x_1,\, x_n^2x_{s(n)}^{-1} :  n\ge 2 \ra\]
    %where $s\colon \Z_{\ge 2}\to\Z_{\ge 1}$ is a computable function such that one cannot decide, given $n\ge 2$, if there exists $k$ such that $f^k(n)=1$. The Word Problem in $G/K$ is decidable, for instance one can decide
\end{rem}

\section{Co-context-free groups}

For groups with decidable Word Problem, the period growth is bounded by a recursive function if and only if the Torsion Problem is decidable. In \cite{co-CF}, it is shown that the Torsion Problem is decidable in co-context-free groups. In this section, we prove a uniform and quantitative version of this fact.
%(Note that $2V$ has undecidable Word Problem \cite{torsion_problem_2V}, hence $p_{2V}(n)$ is not upper bounded by any recursive function.) % Acknowledge Giles Gardam

\subsection{Background}

Before proving the result, we need to recall some general definitions and results around regular and context-free languages. We refer to \cite{language_textbook}.

\begin{defi} A \emph{context-free grammar} is a tuple $(V,\Sigma,R,V_0)$ with
    \begin{itemize}[leftmargin=8mm]
        \item $V$ is a finite set of \emph{variables}.
        \item $\Sigma$ is a finite set called the \emph{alphabet}.
        \item $R\subset V\times(V\sqcup \Sigma)^*$ is a finite set of \emph{rewriting rules} (or \emph{production rules}). Elements $(Y,w)\in R$ are written $Y\to w$.
        \item $V_0\in V$ is the start symbol.
\end{itemize}
A grammar is in \emph{Chomsky normal form} if $\abs w\le 2$ for all rules $Y\to w$.
\smallskip

We define a relation (preorder) $\to^*$ on $(V\sqcup \Sigma)^*$, declaring $aYb\to^* awb$ for all $a,b\in(V\sqcup \Sigma)^*$ and $(Y,w)\in R$, and then taking the reflexive-transitive closure. A formal language $\Lc$ is said to be \emph{context-free} if there exists a context-free grammar $(V,\Sigma,R,V_0)$ such that
\[
\Lc= \left\{w\in\Sigma^* \;\big|\; V_0 \to^{*} w\right\}.
\]
\end{defi}

\begin{lemma}[{\cite[\S 7.1.5]{language_textbook}}] \label{lemma:normal-form}
     Given a context-free grammar, there exists a context-free grammar in Chomsky normal form producing the same language.
\end{lemma}
\bigskip
\begin{defi}
    A  \emph{finite-state automaton} is a tuple $(\Sigma,\mathsf{Q},q_0,\delta,\mathsf{F})$ with
    \begin{itemize}[leftmargin=8mm]
        \item $\Sigma$ is a finite set called the \emph{alphabet}.
        \item $\mathsf{Q}$ is a finite set of \emph{states}.
        \item $q_0\in \mathsf{Q}$ is the \emph{initial state}.
        \item $\delta\subseteq \mathsf{Q}\times \Sigma\times \mathsf{Q}$ is a finite set of \emph{transitions}. Elements $(a,s,b)\in \delta$ should be thought of as an oriented edge from $a$ to $b$, labelled by $s$.
        \item $\mathsf{F}\subseteq \mathsf{Q}$ is the set of \emph{final states}.
    \end{itemize}
    An automaton is \emph{deterministic} (DFA) if, for every $a\in \mathsf{Q}$ and $s\in\Sigma$, there exists at most one state $b\in \mathsf{Q}$ such that $(a,s,b)\in \delta$.
    \smallskip
    
    A formal language $\Rc$ is said to be \emph{regular} if there exists a finite-state automaton $(\Sigma,\mathsf{Q},q_0,\delta,\mathsf{F})$ such that
    \[
    \Rc= \left\{ w\in\Sigma^* \;\big|\; w\text{ labels an oriented path }q_0\to \mathsf{F}\right\}.
    \]
\end{defi}
\begin{center}
    \begin{tikzpicture}
        [thick, scale=.9, every loop/.append style={-latex}]
	\node[state, accepting, initial left, minimum size=10pt] (v0) at (0,0) {};
    \node[state, minimum size=10pt] (v1) at (1.6,0) {};
    \node[state, accepting, minimum size=10pt] (v2) at (3.4,-.7) {};
    
    \begin{scope}[shift={(4.26,1.4)}, rotate=-30, scale=1.2]
	\node[state, minimum size=10pt] (a0) at (1,0) {};
    \node[state, accepting, minimum size=10pt] (a1) at (-.5,.87) {};
    \node[state, minimum size=10pt] (a2) at (-.5,-.87) {};
    \end{scope}

    \begin{scope}[shift={(6,-1.63)}, rotate=18, scale=1.2]
    \node[state, accepting, minimum size=10pt] (b0) at (1,0) {};
    \node[state, minimum size=10pt] (b1) at (.31,.95) {};
    \node[state, minimum size=10pt] (b2) at (-.81,.59) {};
    \node[state, minimum size=10pt] (b3) at (-.81,-.59) {};
    \node[state, minimum size=10pt] (b4) at (.31,-.95) {};
    \end{scope}
		
	\path (v0)  edge [-latex] node [below] {$1$} (v1)
        (v1)  edge [-latex] node [below left] {$1$} (v2)
        (v2)  edge [-latex] node [below left] {$1$} (b2)
        
        (b0)  edge [-latex, bend right] (b1)
        (b1)  edge [-latex, bend right] (b2)
        (b2)  edge [-latex, bend right] (b3)
        (b3)  edge [-latex, bend right] (b4)
        (b4)  edge [-latex, bend right] (b0)

        (v1)  edge [-latex] node [above left] {$1$} (a2)
		(a0)  edge [-latex, bend right=50] (a1)
        (a1)  edge [-latex, bend right=50] (a2)
        (a2)  edge [-latex, bend right=50] node [below right] {$1$} (a0);
    \end{tikzpicture}
    \captionsetup{font=small, width=12cm}
    
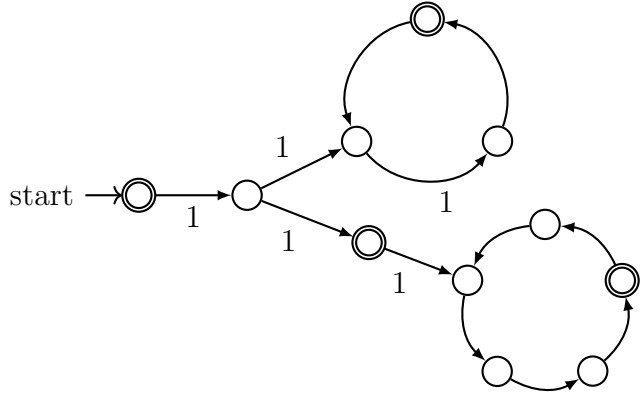
\captionof{figure}{A non-deterministic automaton over the alphabet $\Sigma=\{1\}$ recognising $\Rc=\{1^0,1^2\}\cup\left\{1^{15q+r}\;\big|\; q\ge 0,\; r=4,6,7,10,11,13,16\right\}$.}
    \label{fig:NDA}
\end{center}

Another interesting characterisation of regular language is due to Myhill and Nerode. We first define an equivalence relation on $\Sigma^*$:
\begin{defi}
    Given a language $\Lc\subseteq\Sigma^*$ and two words $u,v\in\Sigma^*$, we define $u\sim_\Lc v$ if we have  $uw\in\Lc\iff vw\in\Lc$ for all $w\in\Sigma^*$.
\end{defi}
\begin{thm}[{\cite[\S 4.4]{language_textbook}}] \label{thm:Myhill_Nerode}
    $\Lc$ is regular if and only if $\Sigma^*/\!\sim_\Lc$ is finite. In this case, $\abs{\Sigma^*/\sim_\Lc}$ is the minimal number of states in a DFA recognising $\Lc$.
\end{thm}
\bigskip

Finally, we recall that intersecting a context-free language with a regular language gives another context-free language. This can be quantified:
\begin{lemma} \label{lemma:regular-cap-context-free}
    Consider two languages $\Lc, \Rc\subseteq \Sigma^*$ such that
    \begin{itemize}[leftmargin=8mm]
        \item $\Lc$ is context-free defined by a grammar with $m$ variables and
        \item $\Rc$ is regular defined by an automata with $n$ states. 
    \end{itemize}
    Then $\Lc\cap \Rc$ is context-free defined by a grammar with $mn^2+1$ variables. Moreover, this construction preserves the length of production rules.
\end{lemma}
\begin{proof}
    Suppose that $\Lc$ is produced by a grammar $(V,\Sigma,R,V_0)$ and $\Rc$ is recognised by an automaton $(\Sigma, \mathsf{Q},q_0, \delta, \mathsf{F})$. Let us define a new context-free grammar $(V',\Sigma,R',V_0')$. We define
    \begin{itemize}[leftmargin=8mm]
        \item $V'=\{V_0'\}\sqcup\{X_{a,b}\mid X\in V, a,b\in \mathsf{Q}\},$
        \item If $R$ contains a rule $X\to w_0Y_1w_1Y_2w_2\ldots Y_nw_n$ with $w_i\in\Sigma^*$ and $Y_i\in V$ and there exists an oriented path $a_i\to b_i$ labelled by $w_i$ in the automaton recognising $\Rc$, then we add to $R'$ the production rule
        \[ X_{a_0,b_n} \to w_0 (Y_1)_{b_0,a_1} w_1 (Y_2)_{b_1,a_2}w_2\ldots (Y_n)_{b_{n-1},a_n}w_n\]
        Finally, we add all the production rules $V_0'\to (V_0)_{q_0,b}$ with $b\in \mathsf{F}$.
    \end{itemize}
    A key observation is that
    \[ \{ w\in \Sigma^* \mid Y_{a,b}\to^* w\}\subseteq \{ w\in\Sigma^*\mid w\text{ labels an oriented path }a\to b\},\]
    by induction. In particular, the language $\Lc'$ produced by $(V',\Sigma,R',V_0')$ is included in $\Rc$. One can check that $\Lc'=\Lc\cap\Rc$ by translating derivations from one grammar to the other.
\end{proof}

\subsection{Period growth}
We are ready to state and prove the main result of this section:
\begin{thm} \label{thm:coCF_period} Let $G$ be a co-context-free group.
    \begin{enumerate}[leftmargin=8mm, label={\normalfont(\alph*)}]
        \item The period growth satisfies $p_G(n)\preceq \exp(n^4)$. 
        \item If $D\subseteq D_k$ for some integer $k$, then $p_G^D(n)\preceq \exp(n^2)$. 
    \end{enumerate}
\end{thm}
\begin{proof}
    We fix a generating set $S$. By \cref{lemma:normal-form}, we may assume that the context-free language $\coWP(G,S)$ is given by a grammar in Chomsky normal form, and let $m$  be the number of variables in this grammar.

    Consider $g\in G$ and fix $w\in S^*$ a word representing it of length $n=\norm g_S$. Due to \cref{lemma:regular-cap-context-free} the language $\coWP(G,S)\cap\{w\}^*$ is context-free, produced by a grammar in Chomsky normal form with $mn^2+1$ variables.
    We apply the map $u\colon s\mapsto 1$ to this language, we get a context-free language
    \begin{align*}
    \Lc_g &:= u(\coWP(G,S)\cap\{w\}^*) \\
    & \;\overset!= \{ 1^{n\ell} : \ell\in\Z_{\ge 0} \text{ such that } \ord(g) \nmid \ell \} \subset \{1\}^*
    \end{align*}
    which is also defined by a grammar with $mn^2+1$ variables. Moreover $\Lc_g$ is unary (meaning the alphabet is a singleton), so it is actually regular (this is a corollary of Parikh's theorem) -- this was quantified in \cite{Unary_CF}.
    
    (a) By \cite[Theorem 6]{Unary_CF}, the language $\Lc_g$ is recognised by a deterministic finite-state automaton with $\abs{\mathsf{Q}}\le 2^{(mn^2+1)^2}$ states. If $g$ has finite order, then the language $\Lc_g$ defines $n\cdot \ord(g)$ Nerode congruences classes. We conclude that $n\cdot \ord(g)\le \abs{\mathsf{Q}}\le 2^{(mn^2+1)^2}$ by \cref{thm:Myhill_Nerode} hence
    \[ p_G(n)\preceq \exp(n^4).\]
    
    (b) By \cite[Theorem 4]{Unary_CF}, the language $\Lc_g$ is recognised by a non-deterministic finite-state automaton with $\abs{\mathsf{Q}}\le 2^{2mn^2+1}+1$ states. Suppose that $\ord(g)$ is finite and belongs to $D_k$. We can write $\ord(g)=p_1^{e_1}p_2^{e_2}\ldots p_k^{e_k}$.
    
    The pumping lemma implies that, for each word $1^x\in \Lc_g$ with $x>\abs{\mathsf{Q}}$, there exists $1\le y_x\le\abs{\mathsf{Q}}$ such that $1^{x+y_xz}\in \Lc_g$ for all $z\ge -1$. In particular,
    \[ \forall x>\abs{\mathsf{Q}},\quad 1^x\in\Lc_g \implies 1^{x+\lcm(1,2,\ldots,\abs{ \mathsf{Q}})}\in\Lc_g, \]
    i.e., the language is eventually $\lcm(1,2,\ldots,\abs{\mathsf{Q}})$-periodic. On the other hand, the minimal period of $\Lc_g$ is clearly $n\cdot\ord(g)$, hence
    \[ p_i^{e_i}\mid n\cdot\ord(g)\mid\lcm(1,2,\ldots,\abs{\mathsf{Q}}).\]
    We deduce that $p_i^{e_i}\le \abs{\mathsf{Q}}$ and finally
    \[ \ord(g)=p_1^{e_1}p_2^{e_2}\ldots p_k^{e_k}\le (2^{2mn^2+1}+1)^k\preceq \exp(n^2). \qedhere\]
\end{proof}
\begin{rem}
    The minimal period of a unary regular language can be larger than the size of a non-deterministic finite-state automaton recognising it. An example is given in \cref{fig:NDA}. In general, an upper bound on the minimal period is given by $g(\abs{ \mathsf{Q}})$, where $g$ is Landau's function defined in \cref*{sec:Houghton}.
\end{rem}

%\begin{rem}
%  \cite[Corollary 1]{Unary_CF} constructs of a unary context-free language defined by a grammar with $\asymp k$ variables, where the smallest deterministic finite state automaton has size $\asymp \exp(k^2)$. This makes heavy use of the $\lcm$ trick. This could be used to construct elements in $V$ of large order, maybe $\exp(n^2)$??
%\end{rem}

\subsection{Distortion of infinite cyclic subgroups}

Membership problem in cyclic subgroups is also decidable in co-context-free groups. (This is another case of the Knapsack Problem, which is decidable in co-context-free groups \cite{Knapsack_Heis}.) We make this result quantitative.
\begin{thm}
    Let $G$ be a co-context-free group, and $g\in G$ an element of infinite order. Then the distortion of the cyclic subgroup $\la g\ra$ satisfies
    \[ \dist^G_{\la g\ra}(n):=\max\left\{k\ge 0 \;\big|\; \|g^k\|_S\le n\right\}\preceq \exp(n^2). \]
\end{thm}
Note that cyclic subgroups of $V$ are undistorted \cite[Theorem 1.9]{Burillo_Cleary_Röver_2017}, so the true answer might be $\dist_{\la g\ra}^G(n)\asymp n$.
\begin{proof}
	Fix a element $g\in G$ of infinite order. The language $\coWP(G,S\cup\{g\})$ is context-free, suppose it is produced by a grammar in Chomsky normal form with $m$ variables. For each word $w\in S^n$, the context-free language $\coWP(G,S\cup\{g\}) \cap w\{g\}^*$ is recognised by a grammar with $m(n+1)^2+1$ variables. It follows that the unary language 
	\begin{align*}
	    \Lc_{g,w}
        & = \{1^0,1^1,\ldots,1^{n-1}\} \sqcup u\bigl(\coWP(G,S\cup\{g\})\cap w\{g\}^*\bigr) \\ 
        & = \{1^0,1^1,\ldots,1^{n-1}\} \sqcup \{1^{n+x}\mid x\ge 0,\, wg^x\ne 1_G\}
	\end{align*}
	is produced by a NFA with $\abs{\mathsf{Q}}\le 2^{2m(n+1)^2+1}+n$ states. On the other hand, any automaton recognising $\{1^x\mid x\ne k\}$ has size $\ge \sqrt{k}$ \cite[Theorem 7.2]{NFA_for_cofinite}. If there exists $k$ such that $wg^k=1$, we conclude that
	\[ k\le n+k \le \abs{\mathsf{Q}}^2 \le (2^{2m(n+1)^2+1}+n)^2 \asymp \exp(n^2).\]
	Therefore, $\dist^G_{\la g\ra}(n)=\max\left\{k\ge 0 \;\big|\; \|g^k\|_S\le n\right\}\preceq \exp(n^2)$.
\end{proof}
% Can we improve the exponent of $n$ in the size of the grammar recognising $\coWP \cap w\{g\}*$ ?
\section{Houghton groups \texorpdfstring{$H_m$}{}} \label{sec:Houghton}

We study the period growth of Houghton groups $H_m$, defined in \cite{houghton}.
\begin{defi}[Landau's function]
    Let $n\ge 0$ be an integer. We define
    \begin{align*}
        g(n)
        & = \max\!\left\{\ord(\sigma) \;\big|\; \sigma\in \Sym(n) \right\} \\
        & = \max\!\left\{ \lcm(m_1,m_2,\ldots,m_k) \;\big|\; m_1+m_2+\ldots+m_k=n \right\}.
    \end{align*}
    More generally $g^D(n)=\max\!\left\{\ord(\sigma) \;\big|\; \sigma\in\Sym(n) \text{ and }\ord(\sigma)\in D\right\}$.
\end{defi}

\begin{thm}[Landau \cite{Landau}]
    $g(n) \asymp \exp\!\big(\sqrt{n\log(n)}\big)$.
\end{thm}
%\begin{proof}[Intuition]
%    Take powers of primes up to $\sqrt n$.
%\end{proof}
\begin{prop}
    If $D_Q\subseteq D\subseteq D_k$ with $Q\subset \P$ and $\abs Q=k$, then $$g^D(n)\asymp n^k$$
\end{prop}
\begin{proof}
    Let $Q=\{q_1,\ldots,q_k\}$. We take $m_i=q_i^{\ell_i}$ between $\frac{n}{q_ik}$ and $\frac{n}{k}$, consider
    \[ \sigma = (1\,\ldots\,m_1)(m_1+1\,\ldots\,m_1+m_2)\ldots \in \Sym(n).\]
    This gives $g^D(n)\ge n^{k}\cdot \prod (q_ik)^{-1}$.

    For the matching upper bound, the order of a $D$-torsion element $\sigma$ is bounded by the product of the length of $k$ cycles in its decomposition (one with largest $q_i$ exponent for each $q_i$). It follows that $\ord(g)\le n^{k}$.
\end{proof}
\bigskip

\begin{thm} \label{thm:period_H}
    Let $H_2=\FSym(\Z)\rtimes\Z$ and $D\subseteq\N$. Then $p_{H_2}^D(n) \asymp g^D(n)$.
\end{thm}
\begin{proof} We fix integers $m_1,m_2,\ldots,m_k\in\Z_{>0}$ such that
    \[ m_1+\ldots+m_k=n \quad\text{and}\quad \lcm(m_1,\ldots,m_k)=g^D(n). \]
    For the lower bound, we consider the element
    \begin{align*}
    g & = (ta)^{m_1-1}t(ta)^{m_2-1}t\ldots t(ta)^{m_k-1}t^{-m_1-\ldots-m_k} \\
    & = (1\,2\,\ldots\,m_1)(m_1+1\,m_1+2\,\ldots\, m_1+m_2)\ldots
    \end{align*}
    with $a=(0\,1)$ and $t\colon x\mapsto x+1$. This element has length $\norm g\le 3n$, and order $\ord(g)=g^D(n)$.
    
    For the upper bound, if $g\in H_2$ satisfies $\ord(g)\in D$, then $g\in\FSym(\Z)$. If moreover $\norm g\le n$, then $\abs{\supp(g)}\preceq n$ and therefore $\ord(g)\preceq g^D(n)$.
\end{proof}
\begin{rem}
    This can be easily generalised to other Houghton groups $H_m$ (with $m\ge 2$), and $\FSym(G)\rtimes G$ with $G$ infinite torsion-free. Indeed, the lower bound follows from the fact that $H_2$ is a subgroup and by \cref{lem:sub-over-group}. For the upper bound, we observe that for an element to be torsion and of length at most $n$, its support is of cardinality $\preceq n$.
\end{rem}
\section{Length estimates in Thompson group \texorpdfstring{$V$}{V}}
\label{sec:metric}

\subsection{Background}

We fix symmetric generating sets $S_F\subset S_T\subset S_V$ for $F$, $T$ and $V$ respectively, and denote by $\norm{\vardot}_F,\norm{\vardot}_T,\norm{\vardot}_V$ the corresponding word lengths.  Moreover, we use right composition for permutations (so $(\sigma\tau)(x)=((x)\sigma)\tau$).

\medskip

We recall a few results of Birget and Burillo--Cleary--Stein--Taback:
\begin{thm} \label{thm:metric_all}
Consider $g\in V$ whose reduced diagram has $n$ leaves.
\begin{enumerate}[leftmargin=8mm, label={\normalfont(\alph*)}]
    \item If $g\in F$, then $\norm{g}_F\asymp n$. %$\frac1Cn\le\norm{g}_F\le Cn$
    \hfill{ \normalfont \cite{ThompsonV_metric}}
    \item If $g\in T$, then $\norm{g}_T\asymp n$. \hfill{\normalfont\cite{ThompsonT_metric}}
    \item If $g\in V$, then $n\preceq \norm{g}_V\preceq n\log(n)$. Moreover, there exists $c>0$ such that elements satisfying $\norm{g}_V\ge cn\log(n)$ are exponentially generic. \hfill{\normalfont\cite{ThompsonV_metric}}
\end{enumerate}
\end{thm}

One of the key ingredients in Birget's estimates (c) is that the top and bottom trees do not matter much in the length of an element of $V$:
\begin{defi}[Induced length on $\Sym(n)$]
    Given $\sigma\in\Sym(n)$, we define
    \[ L(\sigma) \coloneqq \min\left\{\norm{(\Uc,\sigma,\Vc)}_V \;\big|\; \Uc,\Vc\text{ are trees with $n$ leaves}\right\}. \]
\end{defi}
\begin{lemma} \label{lem:induced_length_composition}
There exists $C>0$ such that, for all positive integer $n$, all permutations $\sigma,\tau\in\Sym(n)$ and all trees $\Uc,\Vc$ with $n$ leaves, we have
\begin{enumerate}[leftmargin=8mm, label={\normalfont(\alph*)}]
    \item $\norm{(\Uc,\sigma,\Vc)}_V\le L(\sigma)+2Cn$, and
    \item $L(\sigma\tau)\le L(\sigma)+L(\tau)+Cn$.  
\end{enumerate}
\end{lemma}
\begin{proof}
    (a) Suppose that $L(\sigma)=\norm{(\Uc_0,\sigma,\Uc_0)}_V$, then
    \begin{align*} 
    \norm{(\Uc,\sigma,\Vc)}_V
    & \le \norm{(\Uc,\id,\Uc_0)}_V + \norm{(\Uc_0,\sigma,\Vc_0)}_V + \norm{(\Vc_0,\id,\Vc)}_V \\
    & \le Cn + L(\sigma)+Cn,
    \end{align*}
    using \cref{thm:metric_all}(a) for $(\Uc,\id,\Uc_0),(\Vc_0,\id,\Vc)\in F$.
    
    (b) Suppose that $L(\sigma)=\norm{(\Uc_0,\sigma,\Vc_0)}_V$ and $L(\tau)=\norm{(\Uc_1,\tau,\Vc_1)}_V$, then
    \begin{align*}
    L(\sigma\tau)
    & \le\norm{(\Uc_0,\sigma\tau,\Vc_1)}_V \\
    & \le \norm{(\Uc_0,\sigma,\Vc_0)}_V+\norm{(\Vc_0,\id,\Uc_1)}_V+\norm{(\Uc_1,\tau,\Vc_1)}_V \\
    & \le L(\sigma)+Cn+L(\tau). \qedhere
    \end{align*}
\end{proof}

\subsection{Riffle-shuffle permutations}

\begin{defi}
Given a permutation $\sigma\in\Sym(n)$, we define two statistics:
\begin{itemize}[leftmargin=7mm]
    \item The longest decreasing sequence of $\sigma$ is
    \[ \LDS(\sigma)= \max\bigl\{ \abs A : A\subseteq\{1,2\ldots,n\bigr\} \text{ and } \sigma|_A\text{ is decreasing}\}. \]
    %Equivalently, such that there exists $A\subseteq\{1,2,\ldots,n\}$ with $\abs A=p$ such that $\sigma|_A$ is decreasing.
    Equivalently, Dilworth's theorem \cite{Dilworth} implies that $\LDS(\sigma)$ coincides with the minimum number of parts in a partition
\[ \{1,2,\ldots,n\}=C_1\sqcup C_2\sqcup \ldots \sqcup C_q \]
such that $\sigma|_{C_i}$ is increasing.

    \item The rising number is $r(\sigma)=\#\{i\in\{1,\ldots,n-1\}:\sigma(i)>\sigma(i+1)\}+1$. Equivalently, $r(\sigma)$ is the smallest number of parts in a partition
    \[ \{1,2,\ldots,n\}=B_1\sqcup B_2\sqcup \ldots \sqcup B_q \]
    where $B_i$ consists of \emph{consecutive} integers and $\sigma|_{B_i}$ is increasing.
\end{itemize}
\end{defi}

\begin{defi}
$\sigma\in\Sym(n)$ is a \emph{riffle-shuffle permutation} if $r(\sigma)\le 2$.
\end{defi}

\begin{center}
    \begin{tikzpicture}[xscale=.9]
        \draw[thick, red] (1,1) -- (2,0);
        \draw[thick, red!66!yellow] (2,1) -- (4,0);
        \draw[thick, red!33!yellow] (3,1) -- (6,0);
        
        \draw[line width=2pt, white] (4,1) -- (1,0);
        \draw[thick, blue] (4,1) -- (1,0);
        \draw[line width=2pt, white] (5,1) -- (3,0);
        \draw[thick, blue!70!red] (5,1) -- (3,0);
        \draw[line width=2pt, white] (6,1) -- (5,0);
        \draw[thick, blue!40!red] (6,1) -- (5,0);

        \foreach \x in {1,...,6}{
            \draw[fill=black] (\x,1) ellipse (1pt and .9pt);
            \draw[fill=black] (\x,0) circle (1pt and .9pt);
        }
    \end{tikzpicture} \hspace{15mm}
    \begin{tikzpicture}[xscale=.6]
        \draw[thick, red] (1,1) -- (2,0);
        \draw[thick, red] (2,1) -- (3,0);
        \draw[thick, red!66!yellow] (3,1) -- (6,0);
        \draw[thick, red!33!yellow] (4,1) -- (8,0);
        \draw[thick, red!33!yellow] (5,1) -- (9,0);
        \draw[thick, red!33!yellow] (6,1) -- (10,0);
        
        \draw[line width=2pt, white] (7,1) -- (1,0);
        \draw[thick, blue] (7,1) -- (1,0);
        \draw[line width=2pt, white] (8,1) -- (4,0);
        \draw[thick, blue!70!red] (8,1) -- (4,0);
        \draw[line width=2pt, white] (9,1) -- (5,0);
        \draw[thick, blue!70!red] (9,1) -- (5,0);
        \draw[line width=2pt, white] (10,1) -- (7,0);
        \draw[thick, blue!40!red] (10,1) -- (7,0);

        \foreach \x in {1,...,10}{
            \draw[fill=black] (\x,1) circle (1.5pt and .9pt);
            \draw[fill=black] (\x,0) circle (1.5pt and .9pt);
        }
    \end{tikzpicture}
    
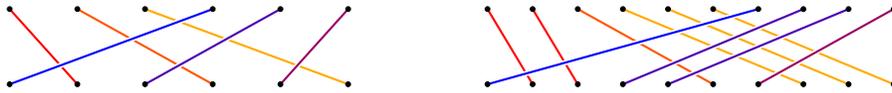
\captionof{figure}{Riffle-shuffle permutations}
    \label{fig:riffle_shuffle_example}
\end{center}

\begin{prop} \label{prop:riffle_length}
    The subset $R\subseteq\Sym(n)$ of riffle-shuffle permutations is a generating set of the symmetric group, and 
    \[ \forall\sigma\in\Sym(n),\quad \norm{\sigma}_R = \left\lceil\log_2(r(\sigma))\right\rceil. \]
\end{prop}
\begin{proof}
    Random walks on the Cayley graph of $\Sym(n)$ with respect to $R$ have been studied extensively. In particular, celebrated work of Bayer and Diaconis gives the exact formula
    \[ p^{(m)}(\sigma) = \frac1{2^{mn}}\binom{2^m+n-r(\sigma)}{n} \]
    for the probability that $m$ random shuffles compose to $\sigma$ \cite[Theorem 1]{BayerDiaconis}. We conclude using that $\norm{\sigma}_R= \min\bigl\{m\ge 0: p^{(m)}(\sigma)>0\bigr\}$.
\end{proof}
\begin{rem}
    The set $R$ is not symmetrical, hence $\norm{\sigma^{-1}}_R\ne \norm{\sigma}_R$.
\end{rem}
\subsection{Bounds on length}

We improve upon \cref*{thm:metric_all}(c), giving a linear upper bound in many cases.

\begin{thm} \label{thm:upper_bound_LDS}
    Consider an element $g\in V$ given by a diagram $(\Uc,\sigma,\Vc)$ with $n$ leaves. Its word length is bounded by $\norm{g}_V \preceq n\cdot \log_2(\LDS(\sigma))$.
\end{thm}
This does not contradict \cref{thm:metric_all} since we have $\E[\LDS(\sigma)] \sim 2\sqrt n$ for permutations $\sigma\in\Sym(n)$ picked uniformly at random \cite{Kerov_Vershik}, and even $\P[\LDS(\sigma)\le c\sqrt n]\to 0$ exponentially fast as $n\to\infty$ for any $c<2$ \cite{large_deviation}.
\bigskip

\begin{lemma} \label{lem:induced_length_riffle}
    $L(\sigma) \preceq n$ for all riffle-shuffle permutations $\sigma\in\Sym(n)$.
\end{lemma}
\begin{proof} We observe that riffle-shuffle permutations $\sigma$ can be obtained by \say{cloning} some strands in a permutation $\tau$ from one of the four families described in \cref*{fig:riffle_shuffle_families}. For instance, the permutation on the right of \cref{fig:riffle_shuffle_example} can be obtained from the permutation on the left (also $\tau_3$ in \cref*{fig:riffle_shuffle_families}), as suggested by the colors. This argument is explained in \cite[\S1]{burillo_note}.
    \begin{center}
        \begin{tikzpicture}[xscale=.9, yscale=1]

        \node at (-1.25,.5) {$\tau_n=$};
        \draw[thick, red] (1,1) -- (2,0);
        \draw[thick, red] (2,1) -- (4,0);
        \draw[thick, red] (3,1) -- (6,0);
        
        \draw[line width=2pt, white] (4,1) -- (1,0);
        \draw[thick, blue] (4,1) -- (1,0);
        \draw[line width=2pt, white] (5,1) -- (3,0);
        \draw[thick, blue] (5,1) -- (3,0);
        \draw[line width=2pt, white] (6,1) -- (5,0);
        \draw[thick, blue] (6,1) -- (5,0);

        \foreach \x in {1,...,6}{
            \draw[fill=black] (\x,1) ellipse (1pt and .9pt);
            \draw[fill=black] (\x,0) circle (1pt and .9pt);
        }

        \begin{scope}[shift={(0,-1.5)}]
        \node at (-1.25,.5) {$\tau'_n=$};
        \draw[thick, red] (1,1) -- (2,0);
        \draw[thick, red] (2,1) -- (4,0);
        \draw[thick, red] (3,1) -- (6,0);
        
        \draw[line width=2pt, white] (4,1) -- (1,0);
        \draw[thick, blue] (4,1) -- (1,0);
        \draw[line width=2pt, white] (5,1) -- (3,0);
        \draw[thick, blue] (5,1) -- (3,0);
        \draw[line width=2pt, white] (6,1) -- (5,0);
        \draw[thick, blue] (6,1) -- (5,0);
        \draw[line width=2pt, white] (7,1) -- (7,0);
        \draw[thick, blue] (7,1) -- (7,0);
        
        \foreach \x in {1,...,7}{
            \draw[fill=black] (\x,1) ellipse (1pt and .9pt);
            \draw[fill=black] (\x,0) circle (1pt and .9pt);}
        \end{scope}

        \begin{scope}[shift={(0,-3)}]
        \node at (-1.25,.5) {$\tau''_n=$};
        \draw[thick, red] (0,1) -- (0,0);
        \draw[thick, red] (1,1) -- (2,0);
        \draw[thick, red] (2,1) -- (4,0);
        \draw[thick, red] (3,1) -- (6,0);
        
        \draw[line width=2pt, white] (4,1) -- (1,0);
        \draw[thick, blue] (4,1) -- (1,0);
        \draw[line width=2pt, white] (5,1) -- (3,0);
        \draw[thick, blue] (5,1) -- (3,0);
        \draw[line width=2pt, white] (6,1) -- (5,0);
        \draw[thick, blue] (6,1) -- (5,0);
        
        \foreach \x in {0,...,6}{
            \draw[fill=black] (\x,1) ellipse (1pt and .9pt);
            \draw[fill=black] (\x,0) circle (1pt and .9pt);}
        \end{scope}
        
        \begin{scope}[shift={(0,-4.5)}]
        \node at (-1.25,.5) {$\tau'''_n=$};
        \draw[thick, red] (0,1) -- (0,0);
        \draw[thick, red] (1,1) -- (2,0);
        \draw[thick, red] (2,1) -- (4,0);
        \draw[thick, red] (3,1) -- (6,0);
        
        \draw[line width=2pt, white] (4,1) -- (1,0);
        \draw[thick, blue] (4,1) -- (1,0);
        \draw[line width=2pt, white] (5,1) -- (3,0);
        \draw[thick, blue] (5,1) -- (3,0);
        \draw[line width=2pt, white] (6,1) -- (5,0);
        \draw[thick, blue] (6,1) -- (5,0);
        \draw[line width=2pt, white] (7,1) -- (7,0);
        \draw[thick, blue] (7,1) -- (7,0);
        
        \foreach \x in {0,...,7}{
            \draw[fill=black] (\x,1) ellipse (1pt and .9pt);
            \draw[fill=black] (\x,0) circle (1pt and .9pt);}
        \end{scope}
    \end{tikzpicture}
        \captionsetup{font=small}
        
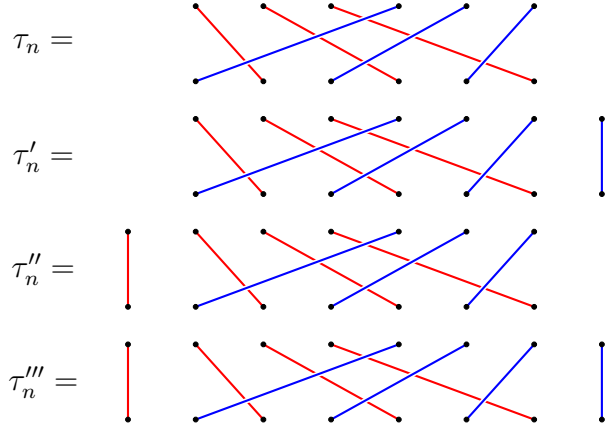
\captionof{figure}{The permutations $\tau$, $\tau'_n$, $\tau''_n$ and $\tau'''_n$, here for $n=3$.}
        \label{fig:riffle_shuffle_families}
    \end{center}
We prove that $L(\sigma)\le L(\tau)$. Consider $g=(\Uc,\tau,\Vc)$ such that $\norm g_V=L(\tau)$. We just have to observe that $g$ is also represented by a (reducible) tree diagram with permutation $\sigma$ (see \cref*{fig:cloning}), hence $L(\sigma)\le\norm{g}_V=L(\tau)$. Therefore, we only have to prove the statement for these four families. 
\begin{center}
    \begin{tikzpicture}[scale=.9]

        \node at (0.25,.75) {$g=$};
        \begin{scope}[very thick]
            % \tau
            \draw[red] (1,1.5) -- (2,0);
            \draw[red!66!yellow] (2,1.5) -- (4,0);
            \draw[red!33!yellow] (3,1.5) -- (6,0);
            
            \draw[line width=2pt, white] (4,1.5) -- (1,0);
            \draw[blue] (4,1.5) -- (1,0);
            \draw[line width=2pt, white] (5,1.5) -- (3,0);
            \draw[blue!70!red] (5,1.5) -- (3,0);
            \draw[line width=2pt, white] (6,1.5) -- (5,0);
            \draw[blue!40!red] (6,1.5) -- (5,0);

            % U
            %\draw (1,1.5) -- (3.5,4);
            \draw (2,1.5) -- (2.5,2);
            \draw (3,1.5) -- (2,2.5);
            \draw (4,1.5) -- (5,2.5);
            \draw (5,1.5) -- (5.5,2);
            \draw (6,1.5) -- (3.5,4) -- (1,1.5);
            
            % V
            %\draw (1,0) -- (3.5,-2.5);
            \draw (2,0) -- (1.5,-.5);
            \draw (3,0) -- (4.5,-1.5);
            \draw (4,0) -- (4.5,-.5);
            \draw (5,0) -- (4,-1);
            \draw (6,0) -- (3.5,-2.5) -- (1,0);
            
        \node at (6.75,.75) {$=$};
    
        \begin{scope}[shift={(6.5,0)}]
            % \sigma
            \draw[red] (.85,1.35) -- (1.85,.15);
            \draw[red] (1.15,1.35) -- (2.15,.15);
            \draw[red!66!yellow] (2,1.5) -- (4,0);
            \draw[red!33!yellow] (2.7,1.3) -- (5.7,0.2);
            \draw[red!33!yellow] (3,1.3) -- (6,0.2);
            \draw[red!33!yellow] (3.3,1.3) -- (6.3,0.2);
            
            \draw[line width=2pt, white] (4,1.5) -- (1,0);
            \draw[blue] (4,1.5) -- (1,0);
            \draw[line width=2pt, white] (4.85,1.35) -- (2.85,0.15);
            \draw[blue!70!red] (4.85,1.35) -- (2.85,0.15);
            \draw[line width=2pt, white] (5.15,1.35) -- (3.15,0.15);
            \draw[blue!70!red] (5.15,1.35) -- (3.15,0.15);
            \draw[line width=2pt, white] (6,1.5) -- (5,0);
            \draw[blue!40!red] (6,1.5) -- (5,0);
            
            % U'
            %\draw (1,1.5) -- (3.5,4);
            \draw[red] (0.85,1.35) -- (1,1.5) -- (1.15,1.35);
            \draw (2,1.5) -- (2.5,2);
            \draw (3,1.5) -- (2,2.5);
            \draw[red!33!yellow] (2.7,1.3) -- (3,1.5) -- (3.3,1.3);
            \draw[red!33!yellow] (3.15,1.4) -- (3,1.3);
            \draw (4,1.5) -- (5,2.5);
            \draw (5,1.5) -- (5.5,2);
            \draw[blue!70!red] (4.85,1.35) -- (5,1.5) -- (5.15,1.35);
            \draw (6,1.5) -- (3.5,4) -- (1,1.5);
            
            % V'
            %\draw (1,0) -- (3.5,-2.5);
            \draw (2,0) -- (1.5,-.5);
            \draw[red] (1.85,.15) -- (2,0) -- (2.15,.15);
            \draw (3,0) -- (4.5,-1.5);
            \draw[blue!70!red] (2.85,.15) -- (3,0) -- (3.15,.15);
            \draw (4,0) -- (4.5,-.5);
            \draw (5,0) -- (4,-1);
            \draw(6,0) -- (3.5,-2.5) -- (1,0);
            \draw[red!33!yellow] (5.7,.2) -- (6,0) -- (6.3,.2);
            \draw[red!33!yellow] (6.15,0.1) -- (6,.2);
        \end{scope}
        \end{scope}
    \end{tikzpicture}
    \captionsetup{font=small}
    
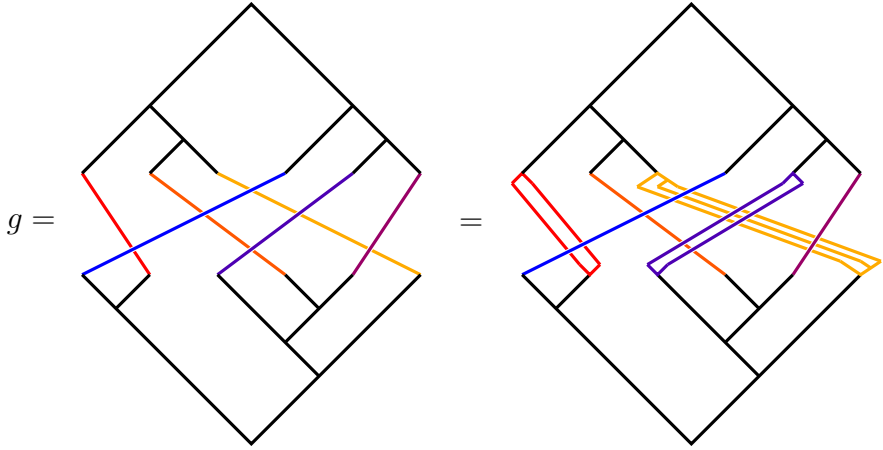
\captionof{figure}{$g$ given by two tree diagrams, with permutations $\sigma$ and $\tau$.}
    \label{fig:cloning}
\end{center}
We first treat the case of permutations $\tau_n\in\Sym(2n)$ defined by
\[
\tau_n(i)=2i \quad\text{and}\quad \tau_n(n+i)=2i-1 \quad\text{for all }1\le i\le n.
\]
To bound $L(\tau_n)$, we estimate the length of $t_n=(\Rc_n\wedge\Rc_n,\tau_n,\Rc_{2n})$, where $\Rc_n$ is the right comb with $n$ leaves and  \say{$\wedge$} is the concatenation (see \cref*{fig:riffle_shuffle_rec}).
\begin{center}
    \begin{tikzpicture}[scale = .85]
    \begin{scope}[very thick, shift={(0,0)}]
        \node at (1.7,2) {$t_3$};

        \draw[red] (1,0) -- (2,-1);
        \draw[red] (2,0) -- (4,-1);
        \draw[red] (3,0) -- (6,-1);
        \draw[line width=4pt, white] (4,0) -- (1,-1);
        \draw[blue] (4,0) -- (1,-1);
        \draw[line width=4pt, white] (5,0) -- (3,-1);
        \draw[blue] (5,0) -- (3,-1);
        \draw[line width=4pt, white] (6,0) -- (5,-1);
        \draw[blue] (6,0) -- (5,-1);
        
        \draw (1,0) -- (3.5,2.5);
        \draw (2,0) -- (2.5,.5);
        \draw (3,0) -- (2,1);
        \draw (4,0) -- (5,1);
        \draw (5,0) -- (5.5,.5);
        \draw (6,0) -- (3.5,2.5);
        
        \draw (1,-1) -- (3.5,-3.5);
        \draw (2,-1) -- (4,-3);
        \draw (3,-1) -- (4.5,-2.5);
        \draw (4,-1) -- (5,-2);
        \draw (5,-1) -- (5.5,-1.5);
        \draw (6,-1) -- (3.5,-3.5);
    \end{scope}

    \node at (6.75, -.5) {$=$};
    
    \begin{scope}[very thick, shift={(10.5,0)}]
        \node at (1.5,2) {$x_0^{-1}fx_0$};

        \draw[blue] (1,0) -- (1,-1);
        \draw[red] (2,0) -- (2,-1);
        \draw[blue] (3,0) -- (3,-1);
        \draw[blue] (4,0) -- (5,-1);
        \draw[line width=4pt, white] (5,0) -- (4,-1);
        \draw[very thick, red] (5,0) -- (4,-1);
        \draw[red] (6,0) -- (6,-1);
        
        \draw (1,0) -- (3.5,2.5);
        \draw (2,0) -- (4,2);
        \draw (3,0) -- (4.5,1.5);
        \draw (4,0) -- (3.5,.5);
        \draw (5,0) -- (5.5,.5);
        \draw (6,0) -- (3.5,2.5);
        
        \draw (1,-1) -- (3.5,-3.5);
        \draw (2,-1) -- (4,-3);
        \draw (3,-1) -- (4.5,-2.5);
        \draw (4,-1) -- (5,-2);
        \draw (5,-1) -- (5.5,-1.5);
        \draw (6,-1) -- (3.5,-3.5);
    \end{scope}

    \node at (11, -.5) {$\cdot$};
    
    \begin{scope}[very thick, shift={(6.5,0)}]
        \node at (1.4,1.2) {$t_2$};

        \draw[red] (1,0) -- (2,-1);
        \draw[red] (2,0) -- (4,-1);
        \draw[line width=4pt, white] (3,0) -- (1,-1);
        \draw[blue] (3,0) -- (1,-1);
        \draw[line width=4pt, white] (4,0) -- (3,-1);
        \draw[blue] (4,0) -- (3,-1);
        
        \draw (1,0) -- (2.5,1.5);
        \draw (2,0) -- (1.5,.5);
        \draw (3,0) -- (3.5,.5);
        \draw (4,0) -- (2.5,1.5);
        
        \draw (1,-1) -- (2.5,-2.5);
        \draw (2,-1) -- (3,-2);
        \draw (3,-1) -- (3.5,-1.5);
        \draw (4,-1) -- (2.5,-2.5);
    \end{scope}
\end{tikzpicture}
    
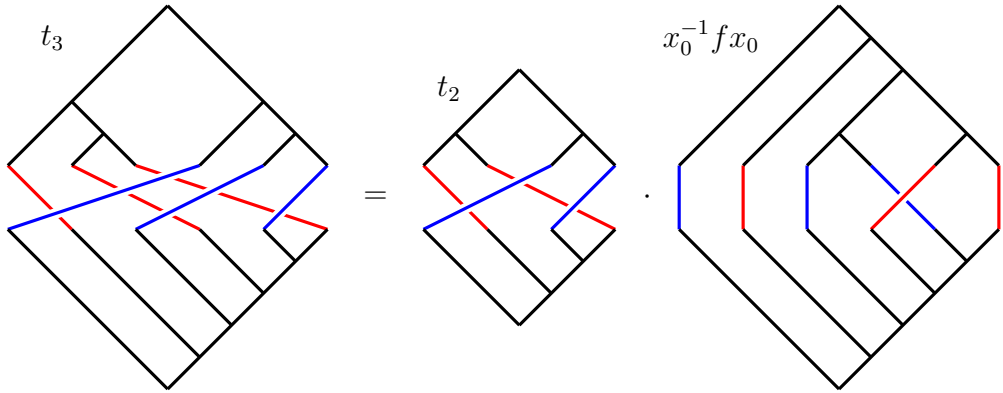
\captionof{figure}{Recurrence relation between $t_n$ and $t_{n-1}$.}
    \label{fig:riffle_shuffle_rec}
\end{center}
We have $t_n= t_{n-1}\cdot x_0^{5-2n} f x_0^{2n-5}$ with $x_0$, $f$ as in \cref{fig:random_V_elements}. It follows that % \todo{I get something different here: I think the right $t_n$ has upper tree consisting of two right combs $R_n$ and lower tree $R_{2n}$}
\[
t_n 
= t_2 \cdot x_0^{-1}fx_0 \cdot \ldots \cdot x_0^{5-2n}f x_0^{2n-5}
= t_2 \cdot x_0^{-1}f(x_0^{-2}f)^{n-3}x_0^{2n-5},
\]
therefore $L(\tau_n)\le \norm{t_n}_V\preceq n$.

Next, we consider the family $\tau'_n$. Consider once again $g=(\Uc,\tau_n,\Vc)$ such that $\norm{g}_V=L(\tau_n)$. We define an injective endomorphism $h\mapsto h@0$, where the element $h@0\in V$ is defined by its action on the Cantor set $\mathfrak C$:
\[
(h@0)(0w)=0h(w) \quad\text{and}\quad (h@0)(1w)=1w.
\]
In a sentence, $h@0$ acts as $h$ below $0$, and trivially elsewhere. We observe that $g@0=(\hspace{1pt}\Uc',\tau'_n,\Vc')$ for some trees $\Uc',\Vc'$, hence
\[
L(\tau'_n)\le \norm{g@0}_V \preceq \norm{g@0}_{V@0} = \norm{g}_V = L(\tau_n) \preceq n.
\]
Similar arguments apply for the two other families.
\end{proof}
\begin{center}
    \begin{tikzpicture}[scale = 1.34, very thick]
    \node at (.2,.77) {$x_0$};
    
    \draw (0,0) -- (1,1) -- (2,0);
    \draw (0.5,.5) -- (1,0);
    
    \draw[blue] (0,0) -- (0,-.8);
    \draw[blue] (1,0) -- (1,-.8);
    \draw[blue] (2,0) -- (2,-.8);

    \draw (0,-.8) -- (1,-1.8) -- (2,-.8);
    \draw (1.5,-1.3) -- (1,-.8);
\end{tikzpicture} \qquad
\begin{tikzpicture}[scale = 1.34, very thick]
    \node at (.2,.77) {$y$};

    % \sigma
    \draw[blue] (0,0) -- (0,-.8);
    \draw[blue] (2,0) -- (1,-.8);
    \draw[line width=4pt, white] (1,0) -- (2,-.8);
    \draw[blue] (1,0) -- (2,-.8);

    % \Uc
    \draw (0,0) -- (1,1) -- (2,0);
    \draw (1.5,.5) -- (1,0);

    % \Vc
    \draw (0,-.8) -- (1,-1.8) -- (2,-.8);
    \draw (1.5,-1.3) -- (1,-.8);
\end{tikzpicture} \qquad
\begin{tikzpicture}[scale = .75, very thick]
        \node at (1.4,1.5) {$f$};

        \draw[blue] (1,0) -- (1,-1);
        \draw[blue] (2,0) -- (2,-1);
        \draw[blue] (3,0) -- (4,-1);
        \draw[line width=4pt, white] (4,0) -- (3,-1);
        \draw[blue] (4,0) -- (3,-1);
        \draw[blue] (5,0) -- (5,-1);
        
        %\draw (1,0) -- (3,2);
        \draw (2,0) -- (3.5,1.5);
        \draw (3,0) -- (2.5,.5);
        \draw (4,0) -- (4.5,.5);
        \draw (5,0) -- (3,2) -- (1,0);
        
        %\draw (1,-1) -- (3,-3);
        \draw (2,-1) -- (3.5,-2.5);
        \draw (3,-1) -- (4,-2);
        \draw (4,-1) -- (4.5,-1.5);
        \draw (5,-1) -- (3,-3) -- (1,-1);
\end{tikzpicture}
    
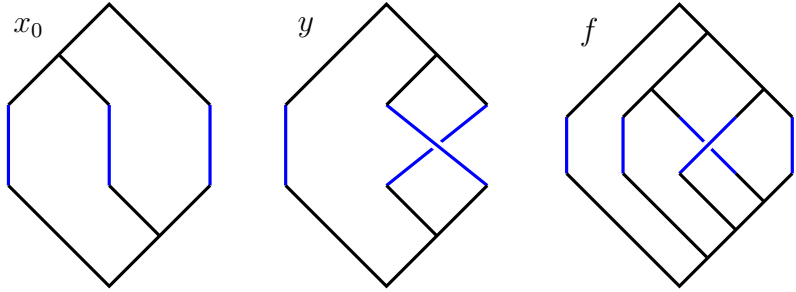
\captionof{figure}{Three elements of $V$}
    \label{fig:random_V_elements}
\end{center}
\begin{proof}[Proof of \cref{thm:upper_bound_LDS}]
Combining Lemma \ref{lem:induced_length_composition}, \ref{lem:induced_length_riffle} and Proposition \ref{prop:riffle_length} gives
\[ \forall \alpha\in\Sym(n),\quad  L(\alpha)\preceq n\cdot\log_2(r(\alpha)).\]
We conclude by observing that any permutation $\sigma$ can be written as $\alpha^{-1}\beta$ with $r(\alpha),r(\beta)\le \LDS(\sigma)$ (see \cref*{fig:rising_vs_LDS}), so that $g=(\hspace{1pt}\Uc,\sigma,\Vc)$ has length
\[ \norm{g}_V \le \norm{(\Uc,\alpha,\Uc)^{-1}}_V + \norm{(\Uc,\beta,\Vc)}_V \preceq n\cdot \log_2(\LDS(\sigma)). \vspace{2mm} \qedhere\]
\end{proof}
\begin{center}
    \begin{tikzpicture}[xscale=.6]
        \draw[thick, red] (1,1.5) -- (2,0);
        \draw[thick, red] (2,1.5) -- (8,0);
        \draw[thick, red] (9,1.5) -- (9,0);
        
        \draw[line width=2pt, white] (3,1.5) -- (5,0);
        \draw[thick, red!33!yellow] (3,1.5) -- (5,0);
        \draw[line width=2pt, white] (6,1.5) -- (10,0);
        \draw[thick, red!33!yellow] (6,1.5) -- (10,0);

        \draw[line width=2pt, white] (4,1.5) -- (4,0);
        \draw[thick, blue] (4,1.5) -- (4,0);
        
        \draw[line width=2pt, white] (5,1.5) -- (1,0);
        \draw[thick, ForestGreen] (5,1.5) -- (1,0);
        \draw[line width=2pt, white] (7,1.5) -- (3,0);
        \draw[thick, ForestGreen] (7,1.5) -- (3,0);
        \draw[line width=2pt, white] (8,1.5) -- (6,0);
        \draw[thick, ForestGreen] (8,1.5) -- (6,0);
        \draw[line width=2pt, white] (10,1.5) -- (7,0);
        \draw[thick, ForestGreen] (10,1.5) -- (7,0);

        \foreach \x in {1,...,10}{
            \draw[fill=black] (\x,1.5) circle (1.5pt and .9pt);
            \draw[fill=black] (\x,0) circle (1.5pt and .9pt);
        }
    \end{tikzpicture} \hspace{10mm}
    \begin{tikzpicture}[xscale=.6]
        \draw[thick, red] (1,1.5) -- (1,.75) -- (2,0);
        \draw[thick, red] (2,1.5) -- (2,.75) -- (8,0);
        \draw[thick, red] (9,1.5) -- (3,.75) -- (9,0);   
        
        \draw[line width=2pt, white] (3,1.5) -- (4,.75) -- (5,0);
        \draw[thick, red!33!yellow] (3,1.5) -- (4,.75) -- (5,0);
        \draw[line width=2pt, white] (6,1.5) -- (5,.75) -- (10,0);
        \draw[thick, red!33!yellow] (6,1.5) -- (5,.75) -- (10,0);

        \draw[line width=2pt, white] (4,1.5) -- (6,.75) -- (4,0);
        \draw[thick, blue] (4,1.5) -- (6,.75) -- (4,0); 
        
        \draw[line width=2pt, white] (5,1.5) -- (7,.75) -- (1,0);
        \draw[thick, ForestGreen] (5,1.5) -- (7,.75) -- (1,0);
        \draw[line width=2pt, white] (7,1.5) -- (8,.75) -- (3,0);
        \draw[thick, ForestGreen] (7,1.5) -- (8,.75) -- (3,0);
        \draw[line width=2pt, white] (8,1.5) -- (9,.75) -- (6,0);
        \draw[thick, ForestGreen] (8,1.5) -- (9,.75) -- (6,0);
        \draw[line width=2pt, white] (10,1.5) -- (10,.75) -- (7,0);
        \draw[thick, ForestGreen] (10,1.5) -- (10,.75) -- (7,0);

        \foreach \x in {1,...,10}{
            \draw[fill=black] (\x,1.5) circle (1.5pt and .9pt);
            \draw[fill=black] (\x,0.75) circle (1.5pt and .9pt);
            \draw[fill=black] (\x,0) circle (1.5pt and .9pt);}

        \node at (11.35,1.125) {$\alpha^{-1}$};
        \node at (11,0.375) {$\beta$};
    \end{tikzpicture}
    \captionsetup{font=small, margin=11mm}
    
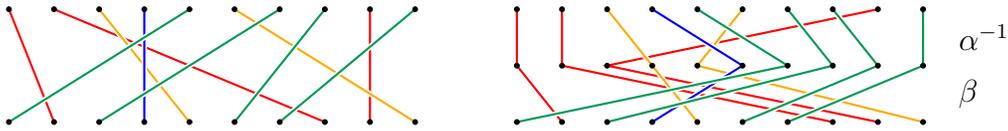
\captionof{figure}{Decomposing $\sigma=\alpha^{-1}\beta$. The sets $C_i$ are indicated in different colors, and $\alpha^{-1}$ rearrange them as sets of consecutive integers.}
    \label{fig:rising_vs_LDS}
\end{center}

We conclude by a complementary proposition, giving the same upper bound in the opposite case where $\LDS(\sigma)=n$ is maximal.
\begin{prop} \label{prop:flip_length}
The permutation $\varphi_n\in\Sym(n)$ defined by $$\varphi_n(i)=n+1-i \quad\text{for all }1\le i\le n$$ satisfies $L(\varphi_n)\preceq n$.    
\end{prop}
\begin{proof}
We consider the element $g_n=(\Rc_n,\varphi_n,\Lc_n)$ where $\Lc_n$ (resp.\ $\Rc_n$) is the left (resp.\ right) comb of size $n$ (see \cref*{fig:flip_length_rec}).
    \begin{center}
    \begin{tikzpicture}[scale = .9]
    \begin{scope}[very thick, shift={(0,0)}]
        \node at (.5,1.5) {$g_5$};
        
        \draw (0,0) -- (2,2) -- (4,0);
        \draw (1,0) -- (2.5,1.5);
        \draw (2,0) -- (3,1);
        \draw (3,0) -- (3.5,.5);

        \draw (0,-.8) -- (2,-2.8) -- (4,-.8);
        \draw (1,-.8) -- (.5,-1.3);
        \draw (2,-.8) -- (1,-1.8);
        \draw (3,-.8) -- (1.5,-2.3);

        \draw[blue] (0,0) -- (4,-.8);
        \draw[blue] (1,0) -- (3,-.8);
        \draw[blue] (2,0) -- (2,-.8);
        \draw[blue] (3,0) -- (1,-.8);
        \draw[blue] (4,0) -- (0,-.8);
    \end{scope}

    \node at (4.75, -.4) {$=$};
    
    \begin{scope}[very thick, shift={(5.5,0)}]
        \node at (.2,1.5) {$x_0^{-2}yx_0^{2}$};
        
        \draw (0,0) -- (2,2) -- (4,0);
        \draw (1,0) -- (2.5,1.5);
        \draw (2,0) -- (3,1);
        \draw (3,0) -- (3.5,.5);

        \draw (0,-.8) -- (2,-2.8) -- (4,-.8);
        \draw (1,-.8) -- (2.5,-2.3);
        \draw (2,-.8) -- (3,-1.8);
        \draw (3,-.8) -- (3.5,-1.3);

        \draw[blue] (0,0) -- (0,-.8);
        \draw[blue] (1,0) -- (1,-.8);
        \draw[blue] (2,0) -- (2,-.8);
        \begin{pgfonlayer}{background}
            \draw[very thick, blue] (4,0) -- (3,-.8);
            \draw[line width=4pt, white] (3,0) -- (4,-.8);
            \draw[very thick, blue] (3,0) -- (4,-.8);
        \end{pgfonlayer}
    \end{scope}

    \node at (10, -.4) {$\cdot$};
    
    \begin{scope}[very thick, shift={(10.5,0)}]
        \node at (.4,1.2) {$g_4$};
        
        \draw (0,0) -- (1.5,1.5) -- (3,0);
        \draw (1,0) -- (2,1);
        \draw (2,0) -- (2.5,.5);

        \draw (0,-.8) -- (1.5,-2.3) -- (3,-.8);
        \draw (1,-.8) -- (.5,-1.3);
        \draw (2,-.8) -- (1,-1.8);

        \draw[blue] (0,0) -- (3,-.8);
        \draw[blue] (1,0) -- (2,-.8);
        \draw[blue] (2,0) -- (1,-.8);
        \draw[blue] (3,0) -- (0,-.8);
    \end{scope}
\end{tikzpicture}
    
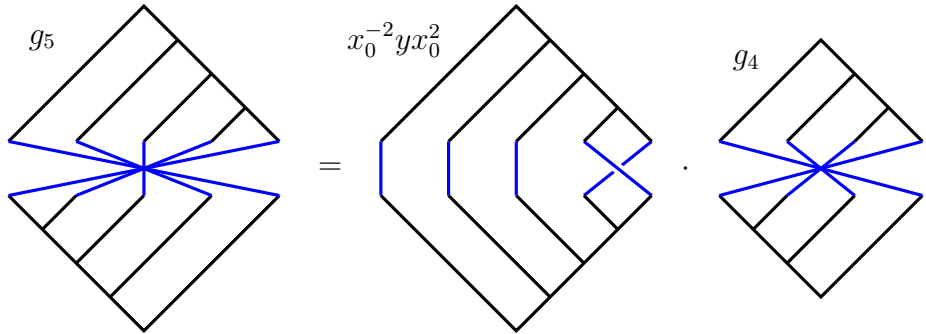
\captionof{figure}{Recurrence relation between $g_n$ and $g_{n-1}$.} \label{fig:flip_length_rec}
\end{center}

We have $g_n= x_0^{-(n-3)} y x_0^{n-3} \cdot g_{n-1}$ with $x_0$, $y$ as in \cref{fig:random_V_elements}. It follows that 
\[
g_n 
= x_0^{3-n}y x_0^{n-3} \cdot x_0^{4-n}y x_0^{n-4} \cdot\ldots\cdot y \cdot g_1 
 = x_0^{3-n}y(x_0y)^{n-3}g_2,
\]
therefore $L(\varphi_n)\le \norm{g_n}_V\preceq n$. \qedhere
\end{proof}
\section{Computing orders in Thompson groups} \label{sec5_diagrams}%\texorpdfstring{$T$}{T} and } 

%In this section, we study period growth for Thompson's groups $T$ and $V$. We use generalised tree diagrams, inspired by Thurston's formalism of train tracks \cite{Thurston}, as used by Calegari in \cite{calegari2007denominator}.

\subsection{Marked strand diagrams}

Elements of Thompson group $V$ can be represented via marked strand diagrams, which are a slight generalisation of abstract strand diagrams (see \cite{strand_diagrams}).

A rotation system on a graph is an assignment of a circular order to the edges incident on each vertex.
% \begin{defi}
% An abstract strand diagram is a finite acyclic digraph together with a
% rotation system, having the following properties:
% \begin{enumerate}
%     \item The graph has a single univalent source and a single univalent sink.
%     \item Every other vertex is either a split or a merge.
% \end{enumerate}
% Two abstract strand diagrams are considered equal if there exists a digraph isomorphism
% between them that is compatible with the corresponding rotation systems.
% \end{defi}

% In order to compose elements more easily (in particular, take powers), we allows more general diagrams.

\begin{defi}
A marked strand diagram is a finite sequence of finite acyclic digraphs together with a rotation system, having the following properties:
\begin{enumerate}[leftmargin=7mm]
    \item Each digraph has univalent source and sink vertices, the number $n$ of source vertices is equal to the number of sink vertices.
    \item There are two ordered $n$-partitions of the Cantor set $\mathfrak C$ and two order-preserving bijective marking maps, one from source vertices to the first partitions and one from sink vertices to the second. 
    \item Every other vertices are a split (\Vsplit) or a merge (\Vmerge).
\end{enumerate}
Two marked strand diagrams are considered equal if there exist digraph isomorphisms between the connected components that are compatible with the corresponding rotation systems and the marking maps.
\end{defi}

\begin{center}
    \begin{tikzpicture}[scale=.8,
        d/.style={circle, draw=black, fill=yellow, inner sep=0mm, line width=.8pt, minimum size = 4pt},
        s/.style={circle, draw=black, fill=black, inner sep=0mm, line width=.8pt, minimum size = 3pt}]

        % Left drawing
        \begin{scope}[very thick, shift={(-5,3.8)}]
            \draw (1,0) -- (2.5,1.5);
            \draw (2,0) -- (3,1);
            \draw (3,0) -- (3.5,.5);
            \draw (4,0) -- (2.5,1.5);

            \draw[very thick] (2.5,1.5) -- (2.5,2.2);

            \draw[out=-135, in=135] (1,0) to (3,-1); 
            \draw[out=-135, in=135] (2,0) to (1,-1);
            \draw[out=-135, in=45] (3,0) to (2,-1);
            \draw[out=-45, in=45] (4,0) to (4,-1);

            \node[label=above:$\varepsilon$] (source) at (2.5,2.2) [s] {};
            \node (eps) at (3,1) [d] {};
            \node (1) at (2.5,1.5) [d] {};
            \node (11) at (3.5,.5) [d] {};

            %\node[] at (2.5,2.5) {$\varepsilon$};
        \end{scope}

        \begin{scope}[very thick, shift={(-5,2.8)}, yscale=-1]
            \draw (1,0) -- (2.5,1.5);
            \draw (2,0) -- (1.5,.5);
            \draw (3,0) -- (3.5,.5);
            \draw (4,0) -- (2.5,1.5);

            \draw[very thick] (2.5,1.5) -- (2.5,2.2);
           
            \node[label=below:$\varepsilon$] (sink) at (2.5,2.2) [s] {};
            \node (eps-) at (2.5,1.5) [d] {};
            \node (0-) at (1.5,.5) [d] {};
            \node (1-) at (3.5,.5) [d] {};

            %\node[] at (2.5,2.5) {$\varepsilon$};
        \end{scope}

        % Right drawing
        \begin{scope}[very thick, shift={(0,4.8)}]
            \draw (3.5,.5) -- (3.5,1.2);
            \draw (3,0) -- (3.5,.5);
            \draw (4,0) -- (3.5,.5);

            \draw[out=-90, in=90] (1.5,1.2) to (2.5,-4.2); 
            \draw[out=-90, in=135] (2.5,1.2) to (1,-3);
            \draw[out=-135, in=45] (3,0) to (2,-3);
            \draw[out=-45, in=90] (4,0) to (3.5,-4.2);

            \node[label=above:$0$] (0) at (1.5,1.2) [s] {};
            \node[label=above:$10$] (10) at (2.5,1.2) [s] {};
            \node[label=above:$11$] (11) at (3.5,1.2) [s] {};
            \node (11split) at (3.5,.5) [d] {};

            %\node[] at (1.5,1.5) {$0$};
            %\node[] at (2.5,1.5) {$10$};
            %\node[] at (3.5,1.5) {$11$};
        \end{scope}

        \begin{scope}[very thick, shift={(0,1.8)},yscale=-1]
            \draw (1.5,.5) -- (1.5,1.2);
            \draw (1,0) -- (1.5,.5);
            \draw (2,0) -- (1.5,.5);

            \node[label=below:$0$] (0) at (1.5,1.2) [s] {};
            \node[label=below:$10$] (10) at (2.5,1.2) [s] {};
            \node[label=below:$11$] (11) at (3.5,1.2) [s] {};
            \node (0split) at (1.5,.5) [d] {};

            %\node[] at (1.5,1.5) {$0$};
            %\node[] at (2.5,1.5) {$10$};
            %\node[] at (3.5,1.5) {$11$};
        \end{scope}
        \end{tikzpicture}
        \captionsetup{font=small}
        \hspace*{10mm}
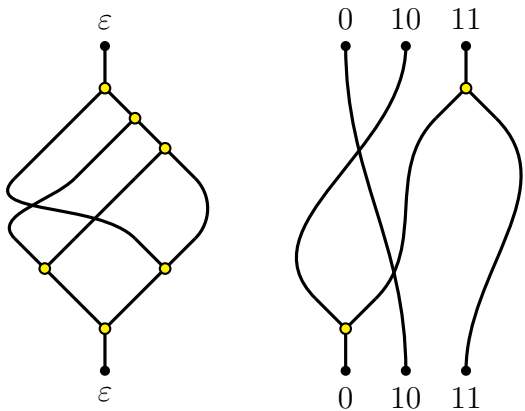
\captionof{figure}{Two marked strand diagrams representing the same element in $V$, one is obtained from the other by applying multiple times the rules in \cref{fig:thompson_rules}.}
        \label{fig:marked_diagram}
    \end{center}

Given a marked strand diagram with sink partition equal to the source partition of another, one can compose them by gluing the sink vertices of the first to the source vertices of the second, according to the markings, and then eliminating the resulting bivalent vertices. It is straightforward that the composition of marked strand diagrams representing elements  $g$ and $h$ in $V$ is a marked strand diagram representing $gh$.\

Starting from a marked strand diagram one can obtain another one representing the same element of $V$, by applying the rules in \cref*{fig:thompson_rules}.

\begin{center}
    \begin{tikzpicture}[scale=.95,
        d/.style={circle, draw=black, fill=yellow, inner sep=0mm, line width=.8pt, minimum size = 4pt},
        s/.style={circle, draw=black, fill=black, inner sep=0mm, line width=.8pt, minimum size = 3pt}]

    % Rule I
    \begin{scope}[very thick, shift={(0,0)}]
        \begin{scope}[shift={(0,0)},
        decoration={markings, mark=at position 0.65 with {\arrow{stealth}}}]
            \node (a) at (0,1) [d] {};
            \node (z) at (0,0) [d] {};

            \draw[postaction=decorate] (-.5,1.5) -- (a);
            \draw[postaction=decorate] (.5,1.5) -- (a);
            \draw[postaction=decorate] (a) -- (z);
            \draw[postaction=decorate] (z) -- (-.5,-.5);
            \draw[postaction=decorate] (z) -- (.5,-.5);
        \end{scope}
        \node at (1.25,.8) {\footnotesize Rule~I};
        \draw[thick, latex-latex] (.5,.5) -- (2,.5);
        \begin{scope}[shift={(2.7,0)},
        decoration={markings, mark=at position 0.6 with {\arrow{stealth}}}]
            \draw[postaction=decorate, out=-65, in=65] (-.5,1.5) to (-.5,-.5);
            \draw[postaction=decorate, out=-115, in=115] (.5,1.5) to (.5,-.5);
        \end{scope}
    \end{scope}

    % Rule II
    \begin{scope}[very thick, shift={(0,-2.75)}]
        \begin{scope}[shift={(0,0)},
        decoration={markings, mark=at position 0.65 with {\arrow{stealth}}}]
            \node (a) at (0,1) [d] {};
            \node (z) at (0,0) [d] {};

            \draw[postaction=decorate] (0,1.5) -- (a);
            \draw[postaction=decorate, out=-55, in=55] (a) to (z);
            \draw[postaction=decorate, out=-125, in=125] (a) to (z);
            \draw[postaction=decorate] (z) -- (0,-.5);
        \end{scope}
        \node at (1.25,.8) {\footnotesize Rule~II};
        \draw[thick, latex-latex] (.5,.5) -- (2,.5);
        \begin{scope}[shift={(2.7,0)},
        decoration={markings, mark=at position 0.6 with {\arrow{stealth}}}]
            \draw[postaction=decorate] (0,1.5) -- (0,-.5);
        \end{scope}
    \end{scope}

    % Rule III top
    \begin{scope}[very thick, shift={(4.5,1.2)}]
        \begin{scope}[shift={(0,0)},
        decoration={markings, mark=at position 0.65 with {\arrow{stealth}}}]]
            \draw[thick] (0,0) -- (1.5,0);
            \node[label=above:\small$\alpha$] (a) at (.75,0) [s] {};
            \node (n) at (.75,-.8) [d] {};

            \draw[postaction=decorate] (a) -- (n);
            \draw[postaction=decorate] (n) -- (.3,-1.2);
            \draw[postaction=decorate] (n) -- (1.2,-1.2);
        \end{scope}
        \node at (2.1,-.4) {\footnotesize Rule~III};
        \draw[thick, latex-latex] (1.5,-.7) -- (2.7,-.7);
        \begin{scope}[shift={(2.8,0)},
        decoration={markings, mark=at position 0.55 with {\arrow{stealth}}}]]
            \draw[thick] (0,0) -- (1.5,0);
            \node[label=above:\small$\alpha0$] (a0) at (.4,0) [s] {};
            \node[label=above:\small$\alpha1$] (a1) at (1.1,0) [s] {};

            \draw[postaction=decorate, out=-90, in=70] (a0) to (.3,-1.2);
            \draw[postaction=decorate, out=-90, in=110] (a1) to (1.2,-1.2);
        \end{scope}
    \end{scope}
    
    % Rule III bottom
    \begin{scope}[very thick, very thick, shift={(4.5,-2.75)}, yscale=-1]
        \begin{scope}[shift={(0,0)},
        decoration={markings, mark=at position 0.65 with {\arrow{stealth}}}]]
            \draw[thick] (0,0) -- (1.5,0);
            \node[label=below:\small$\alpha$] (a) at (.75,0) [s] {};
            \node (n) at (.75,-.8) [d] {};

            \draw[postaction=decorate] (n) -- (a);
            \draw[postaction=decorate] (.3,-1.2) -- (n);
            \draw[postaction=decorate] (1.2,-1.2) -- (n);
        \end{scope}
        \node at (2.1,-.8) {\footnotesize Rule~III'};
        \draw[thick, latex-latex] (1.5,-.5) -- (2.7,-.5);
        \begin{scope}[shift={(2.8,0)},
        decoration={markings, mark=at position 0.55 with {\arrow{stealth}}}]]
            \draw[thick] (0,0) -- (1.5,0);
            \node[label=below:\small$\alpha0$] (a0) at (.4,0) [s] {};
            \node[label=below:\small$\alpha1$] (a1) at (1.1,0) [s] {};

            \draw[postaction=decorate, in=-90, out=70] (.3,-1.2) to (a0);
            \draw[postaction=decorate, in=-90, out=110] (1.2,-1.2) to (a1);
        \end{scope}
    \end{scope}

    % Rule IV
    \begin{scope}[very thick, very thick, shift={(9.5,-.25)}]
        \begin{scope}[shift={(0,0)},
        decoration={markings, mark=at position 0.65 with {\arrow{stealth}}}]]
            \draw[thick] (0,0) -- (1.5,0);
            \draw[thick] (0,-1.2) -- (1.5,-1.2);
            \node[label=above:\small$\alpha$] (a) at (.75,0) [s] {};
            \node[label=below:\small$\beta$] (b) at (.75,-1.2) [s] {};
            \draw[postaction=decorate] (a) -- (b);
        \end{scope}
        \node at (2.1,-.8) {\footnotesize Rule~IV};
        \draw[thick, latex-latex] (1.5,-.5) -- (2.7,-.5);
        \begin{scope}[shift={(2.8,0)},
        decoration={markings, mark=at position 0.65 with {\arrow{stealth}}}]]
            \draw[thick] (0,0) -- (1.5,0);
            \draw[thick] (0,-1.2) -- (1.5,-1.2);
            \node[label=above:\small$\alpha0$] (a0) at (.4,0) [s] {};
            \node[label=above:\small$\alpha1$] (a1) at (1.1,0) [s] {};
            \node[label=below:\small$\beta0$] (b0) at (.4,-1.2) [s] {};
            \node[label=below:\small$\beta1$] (b1) at (1.1,-1.2) [s] {};
            \draw[postaction=decorate] (a0) -- (b0);
            \draw[postaction=decorate] (a1) -- (b1);
        \end{scope}
    \end{scope}
\end{tikzpicture}
    \captionsetup{font=small}
    
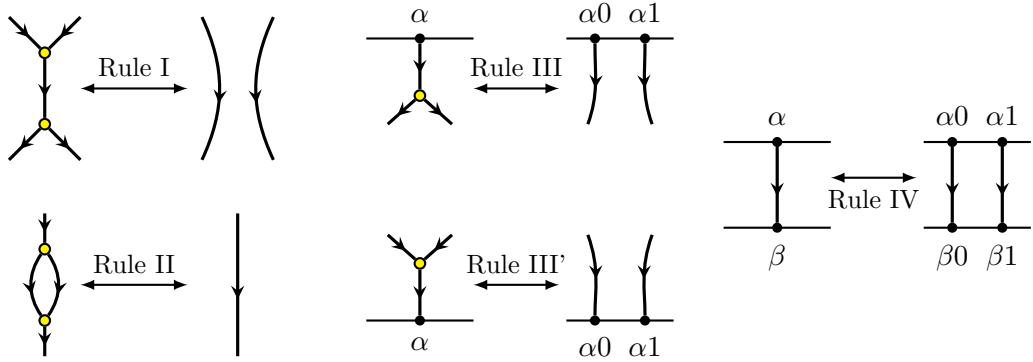
\captionof{figure}{The four rewriting rules for marked strand diagrams.}
    \label{fig:thompson_rules}
\end{center}

\subsection{Detecting torsion and computing order}

We describe an algorithm to decide if an element of Thompson group $V$ has finite order, and in this case to compute efficiently its order.

\subsubsection{Algorithm}\label{subsection:algorithm}

\begin{enumerate}[leftmargin=7mm, label=(\arabic*)]
    \item We write the element on a cylinder and write a label $1$ on the edge crossing the \say{seam} (see left of \cref*{fig:algo_step1}) and label $0$ on all other edges.
    \begin{center}
    \begin{tikzpicture}[scale=.8, d/.style={circle, draw=black, fill=yellow, inner sep=0mm, line width=.8pt, minimum size = 4pt}]
        \begin{scope}[very thick, shift={(-8,3.8)}]
            \draw (1,0) -- (3.5,2.5);
            \draw (2,0) -- (2.5,.5);
            \draw (3,0) -- (2,1);
            \draw (4,0) -- (5,1);
            \draw (5,0) -- (5.5,.5);
            \draw (6,0) -- (3.5,2.5);

            \draw[very thick] (3.5,2.5) -- (3.5,3.2);
            \draw[thick, dotted, red] (1,3.2) -- (6,3.2);
            \node[circle, draw=red, fill=white, inner sep=2pt] at (3.5,3.2) {{$\color{red}1$}};

            \draw[ForestGreen, out=-135, in=135] (1,0) to (3,-1); 
            \draw[ForestGreen, out=-135, in=135] (2,0) to (1,-1);
            \draw[ForestGreen, out=-45, in=45] (3,0) to (6,-1);
            \draw[ForestGreen, out=-135, in=45] (4,0) to (2,-1);
            \draw[ForestGreen, out=-135, in=135] (5,0) to (4,-1);
            \draw[ForestGreen, out=-45, in=45] (6,0) to (5,-1);
            \node[ForestGreen] at (.5,-.5) {$\sigma$};

            \node (eps) at (3.5,2.5) [d] {};
            \node (0) at (2,1) [d] {};
            \node (1) at (5,1) [d] {};
            \node (01) at (2.5,.5) [d] {};
            \node (11) at (5.5,.5) [d] {};
        \end{scope}

        \begin{scope}[very thick, shift={(-8,2.8)}, yscale=-1]
            \draw (1,0) -- (3.5,2.5);
            \draw (2,0) -- (1.5,.5);
            \draw (3,0) -- (4.5,1.5);
            \draw (4,0) -- (4.5,.5);
            \draw (5,0) -- (4,1);
            \draw (6,0) -- (3.5,2.5);

            \draw[very thick] (3.5,2.5) -- (3.5,2.8);
            \draw[thick, dotted, red] (1,2.8) -- (6,2.8);

            \node (eps-) at (3.5,2.5) [d] {};
            \node (0-) at (1.5,.5) [d] {};
            \node (1-) at (4.5,1.5) [d] {};
            \node (10) at (4,1) [d] {};
            \node (101) at (4.5,.5) [d] {};
        \end{scope}

        \begin{scope}[very thick, shift={(0,0)}]
            \draw (1,0) -- (3.5,2.5);
            \draw (2,0) -- (2.5,.5);
            \draw (3,0) -- (2,1);
            \draw (4,0) -- (5,1);
            \draw (5,0) -- (5.5,.5);
            \draw (6,0) -- (3.5,2.5);

            \node at (1,-.4) {\footnotesize$1$};
            \node at (2,-.4) {\footnotesize$2$};
            \node at (6,-.4) {\footnotesize$n$};

            \node (0) at (2,1) [d] {};
            \node (1) at (5,1) [d] {};
            \node (01) at (2.5,.5) [d] {};
            \node (11) at (5.5,.5) [d] {};
        \end{scope}

        \begin{scope}[very thick, shift={(0,7)}, yscale=-1]
            \draw (1,0) -- (3.5,2.5);
            \draw (2,0) -- (1.5,.5);
            \draw (3,0) -- (4.5,1.5);
            \draw (4,0) -- (4.5,.5);
            \draw (5,0) -- (4,1);
            \draw (6,0) -- (3.5,2.5);

            \node at (1,-.4) {\scriptsize$\sigma^{-1}(1)$};
            \node at (6,-.4) {\scriptsize$\sigma^{-1}(n)$};
            
            \node (0-) at (1.5,.5) [d] {};
            \node (1-) at (4.5,1.5) [d] {};
            \node (10) at (4,1) [d] {};
            \node (101) at (4.5,.5) [d] {};
        \end{scope}

        \draw[very thick] (3.5,2.5) -- (3.5,4.5);
        \node[circle, draw=red, very thick, fill=white, inner sep=2pt] at (3.5,3.5) {{$\color{red}1$}};
        \node (eps) at (3.5,2.5) [d] {};
        \node (eps-) at (3.5,4.5) [d] {};
        \end{tikzpicture}
        \captionsetup{font=small}
        \hspace*{10mm}
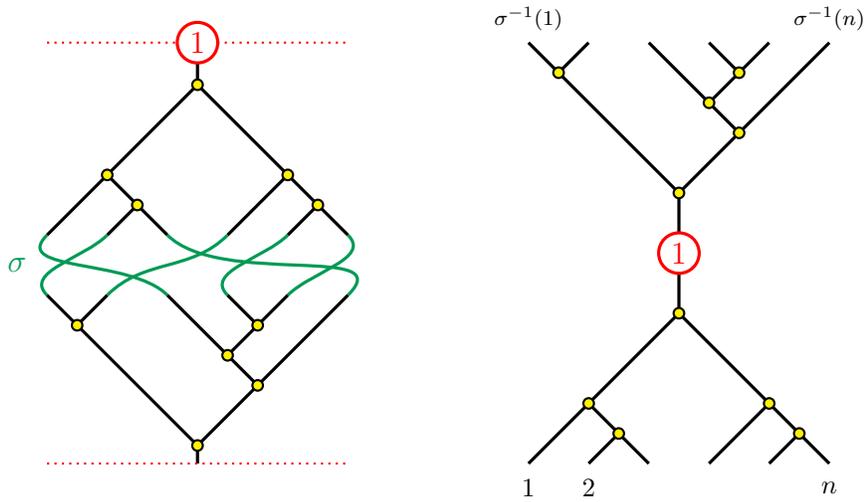
\captionof{figure}{The ideal picture on a cylinder, and the pratical picture.}
        \label{fig:algo_step1}
    \end{center}
    In practice, only the strand diagram and the edge labels matter (not the embedding in the cylinder). We draw the element with the leaves to the exterior and number the leaves using $\sigma$ (see right of \cref*{fig:algo_step1}) and allow ourselves to glue back leaves with the same numbers when doable.
    
    \item If we find a merge-split pair of successive carets, we apply Rule~I, duplicate the central label and sum all labels that end up on each edge. Note that some of the incoming and outgoing edges may coincide, in which case we get slightly different pictures, as explained in \cref*{fig:algo_step2} and \ref*{fig:algo_step2b}.
    \begin{center}
    \begin{minipage}{.45\linewidth}
        \centering
        \begin{tikzpicture}[
	d/.style={circle, draw=black, fill=yellow, inner sep=0mm, line width=.8pt, minimum size = 4pt},
	e/.style={inner sep=0mm, minimum size = 0pt}]
        \clip (-.8,-1.1) rectangle (4.7,1.6);
        \begin{scope}[very thick, decoration={markings, mark=at position 0.65 with {\arrow{stealth}}}]
            \node (a) at (0,1) [d] {};
            \node (z) at (0,-.5) [d] {};
        
            \draw[postaction=decorate] (-.5,1.5) -- (a);
            \draw[postaction=decorate] (.5,1.5) -- (a);
            \draw[postaction=decorate] (a) -- (z);
            \draw[postaction=decorate] (z) -- (-.5,-1);
            \draw[postaction=decorate] (z) -- (.5,-1);

            \node at (-.5,1.05) {$x_1$};
            \node at (.5,1.05) {$x_2$};
            \node at (-.3,.25) {$y$};
            \node at (-.5,-.5) {$z_1$};
            \node at (.5,-.5) {$z_2$};
        \end{scope}

        \begin{scope}[very thick, shift={(3.7,0)},
            decoration={markings, mark=at position 0.57 with {\arrow{stealth}}}]
            \draw[postaction=decorate, out=-65, in=65] (-.5,1.5) to (-.5,-1);
            \draw[postaction=decorate, out=-115, in=115] (.5,1.5) to (.5,-1);

            \node at (-.7,1.05) {$x_1$};
            \node at (-.7,.65) {$+$};
            \node at (-.7,.25) {$y$};
            \node at (-.7,-.15) {$+$};
            \node at (-.7,-.55) {$z_1$};
            
            \node at (.7,1.05) {$x_2$};
            \node at (.7,.65) {$+$};
            \node at (.7,.25) {$y$};
            \node at (.7,-.15) {$+$};
            \node at (.7,-.55) {$z_2$};
        \end{scope}
        \draw[thick, -latex] (.9,.25) -- (2.4,.25);
    \end{tikzpicture}
    \captionsetup{font=small, margin=3mm}
    \captionof{figure}{Rule~I with labels if all $4$ edges are distinct} \label{fig:algo_step2} 
    \end{minipage}\hspace*{10mm}
    \begin{minipage}{.45\linewidth}
    \begin{tikzpicture}[
	d/.style={circle, draw=black, fill=yellow, inner sep=0mm, line width=.8pt, minimum size = 4pt},
	e/.style={inner sep=0mm, minimum size = 0pt}]
        \clip (-.8,-1.1) rectangle (4.7,1.6);
        \begin{scope}[very thick, decoration={markings, mark=at position 0.65 with {\arrow{stealth}}}]
            \node (a) at (0,1) [d] {};
            \node (z) at (0,-.5) [d] {};
        
            \draw[postaction=decorate] (-.5,1.5) -- (a);
            \draw[postaction=decorate] (.5,1.5) -- (a);
            \draw[postaction=decorate] (a) -- (z);
            \draw[postaction=decorate] (z) -- (-.5,-1);
            \draw[postaction=decorate] (z) -- (.5,-1);

            \draw[dotted] (-.5,-1) to[out=10, in=-70] (.5,1.5);
            \node at (-.5,1.05) {$x_1$};
            \node at (.5,1.05) {$x_2$};
            \node at (-.3,.25) {$y$};
            \node at (.5,-.5) {$z_2$};
        \end{scope}

        \begin{scope}[very thick, shift={(3.7,0)},
            decoration={markings, mark=at position 0.57 with {\arrow{stealth}}}]
            \draw[postaction=decorate, out=-65, in=65] (-.5,1.5) to (-.5,-1);
            \draw[postaction=decorate, out=-115, in=115] (.5,1.5) to (.5,-1);
            \draw[dotted] (-.5,-1) to[out=10, in=-70] (.5,1.5);
            
            \node at (-.97,1.05) {\footnotesize$x_1+x_2$};
            \node at (-1.11,.65) {\footnotesize$+2y+z_2$};
        \end{scope}
        \draw[thick, -latex] (.9,.25) -- (2.4,.25);
    \end{tikzpicture}
    \captionsetup{font=small}
    \captionof{figure}{Rule~I with labels if top-right and bottom-left coincide} \label{fig:algo_step2b} 
    \end{minipage}
    \end{center}

    We iterate until there is no more merge-split pairs of successive carets. %In practice, we glue leaves with identical number when doable.
    \item At the end, there are two possibilities:
    \begin{itemize}[leftmargin=5mm]
        \item Either the diagram is a finite union of circles. Then the element is torsion, and its order is the least common multiple of the labels.
        \item Or some carets remains. Then the action of the element on $\mathfrak C$ admits attractor and repellor points, hence the order is infinite. Moreover, the length of finite orbits can be read as the sum of labels along cycles.
        \begin{center}
        \begin{tikzpicture}
        \clip (-2.5,-1.2) rectangle (2.5,1.2);
        \begin{scope}[every node/.style={
            circle, draw=black, fill=yellow,
            inner sep=0mm, line width=.8pt, minimum size = 4pt}]
            \node (1) at (-1.5,.5) {};                    
            \node (2) at (-1,-.5) {};
            \node (3) at (1,-.5) {};
            \node (4) at (1.5,.5) {};
        \end{scope}

        \begin{scope}[very thick]
            \draw[decoration={markings, mark=at position 0.35 with {\arrow{stealth}}}, postaction=decorate, out=-110, in=-90, looseness=2.5] (1) to (-2.3,.4);
            \draw[out=90, in=80, looseness=2] (-2.3,.4) to (1);
            \draw[decoration={markings, mark=at position 0.52 with {\arrow{stealth}}}, postaction=decorate, out=-70, in=110] (1) to (2);
            
            \draw[decoration={markings, mark=at position 0.10 with {\arrow{stealth}}}, postaction=decorate, out=-80, in=110, looseness=2] (2) to (3);
            \draw[decoration={markings, mark=at position 0.10 with {\arrow{stealth}}}, postaction=decorate, out=-100, in=70, looseness=2] (3) to (2);
            
            \draw[decoration={markings, mark=at position 0.35 with {\arrow{stealth}}}, postaction=decorate, out=-70, in=-90, looseness=2.5] (4) to (2.3,.4);
            \draw[out=90, in=100, looseness=2] (2.3,.4) to (4);
            \draw[decoration={markings, mark=at position 0.52 with {\arrow{stealth}}}, postaction=decorate, out=-110, in=70] (4) to (3);
            \end{scope}
\end{tikzpicture}
        \captionsetup{font=small, margin=18mm}
        
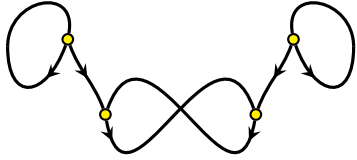
\captionof{figure}{A diagram with remaining merge and split carets, but all merge carets are \say{downstream} from split carets.}
        \end{center}
    \end{itemize}
\end{enumerate}
\begin{exa} \label{exa:adding_machine}
We consider the element $g_n$ from \cref{prop:flip_length}, which can be alternatively defined as follows:
\[ g_n(1^k0w)=0^k1w\;\text{ for }0\le k< n-1, \quad\text{and}\quad g_n(1^{n-1}w)=0^{n-1}w. \]
When reading the first $n$ digits in reverse, $g_n$ adds $1$ modulo $2^{n-1}$, therefore $g_n$ has order $2^{n-1}$. We can verify this using our algorithm, here for $n=3$:
\end{exa}

\begin{center}
    \begin{tikzpicture}[d/.style={circle, draw=black, fill=yellow, inner sep=0mm, line width=.8pt, minimum size = 4pt}]
        % First drawing
        \begin{scope}[very thick, shift={(0,0)}]
            \draw (1,0) -- (2,1);
            \draw (2,0) -- (2.5,.5);
            \draw (3,0) -- (2,1);

            \node at (1,-.4) {\footnotesize$1$};
            \node at (2,-.4) {\footnotesize$2$};
            \node at (3,-.4) {\footnotesize$3$};
        \end{scope}

        \begin{scope}[very thick, shift={(0,4)}, yscale=-1]
            \draw (1,0) -- (2,1);
            \draw (2,0) -- (1.5,.5);
            \draw (3,0) -- (2,1);

            \node at (1,-.4) {\scriptsize$3$};
            \node at (2,-.4) {\scriptsize$2$};
            \node at (3,-.4) {\scriptsize$1$};
        \end{scope}
        \draw[very thick] (2,1) -- (2,3);

        \node (a) at (1.5,3.5) [d] {};
        \node (b) at (2,3) [d] {};

        \node (a) at (2.5,.5) [d] {};
        \node (b) at (2,1) [d] {};
        
        \node[circle, draw=red, very thick, fill=white, inner sep=2pt] at (2,2) {{$\color{red}1$}};

        \draw[thick, -latex] (3,2) -- (4.5,2);
        \node at (3.7,2.3) {\small Rule~I};

        % Second drawing
        \begin{scope}[shift={(4,0)}]
            \begin{scope}[very thick, shift={(0,0)}]
            \draw (1.5,0) -- (1.5,3.5);
            \draw (2,0) -- (2.5,.5);
            \draw (3,0) -- (2.5,.5);

            \node at (1.5,-.4) {\footnotesize$1$};
            \node at (2,-.4) {\footnotesize$2$};
            \node at (3,-.4) {\footnotesize$3$};
        \end{scope}

        \begin{scope}[very thick, shift={(0,4)}, yscale=-1]
            \draw (1,0) -- (1.5,.5);
            \draw (2,0) -- (1.5,.5);
            \draw (2.5,0) -- (2.5,3.5);

            \node at (1,-.4) {\scriptsize$3$};
            \node at (2,-.4) {\scriptsize$2$};
            \node at (2.5,-.4) {\scriptsize$1$};
        \end{scope}

        \node (a) at (1.5,3.5) [d] {};
        %\node (b) at (1.5,2.75) [d] {};
        %\node (c) at (1.5,1.25) [d] {};

        %\node (a1) at (2.5,2.75) [d] {};
        %\node (b1) at (2.5,1.25) [d] {};
        \node (c1) at (2.5,.5) [d] {};

        \node[circle, draw=red, very thick, fill=white, inner sep=2pt] at (1.5,2) {{$\color{red}1$}};
        \node[circle, draw=red, very thick, fill=white, inner sep=2pt] at (2.5,2) {{$\color{red}1$}};
        \end{scope}

        \draw[thick, -latex] (7.3,2) -- (8.3,2);

        % Third drawing 
        \begin{scope}[shift={(7.5,0)}]
            \begin{scope}[very thick, shift={(0,0)}]
            %\draw (1.5,0) -- (1.5,3.5);
            \draw (1,0) -- (1.5,.5);
            \draw (2,0) -- (1.5,.5);

            %\node at (1.5,-.4) {\footnotesize$1$};
            \node at (1,-.4) {\footnotesize$2$};
            \node at (2,-.4) {\footnotesize$3$};
        \end{scope}

        \begin{scope}[very thick, shift={(0,4)}, yscale=-1]
            \draw (1,0) -- (1.5,.5);
            \draw (2,0) -- (1.5,.5);
            \draw (1.5,.5) -- (1.5,3.5);

            \node at (1,-.4) {\scriptsize$3$};
            \node at (2,-.4) {\scriptsize$2$};
            %\node at (2.5,-.4) {\scriptsize$1$};
        \end{scope}

        \node (a) at (1.5,3.5) [d] {};
        \node (c) at (1.5,.5) [d] {};
        
        \node[circle, draw=red, very thick, fill=white, inner sep=2pt] at (1.5,2) {{$\color{red}2$}};
        \end{scope}

        \draw[thick, -latex] (9.5,2) -- (10.5,2);
        \node at (10,2.3) {\footnotesize Rule~I};

        % Fourth drawing
        \begin{scope}[shift={(10.5,0)}]
            \begin{scope}[very thick, shift={(0,0)}]
            \draw (.5,0) -- (.5,4);
            % \draw (2,0) -- (2.5,.5);
            % \draw (3,0) -- (2.5,.5);

            \node at (.5,-.4) {\footnotesize$2$};
             \node at (1.5,-.4) {\footnotesize$3$};
            % \node at (3,-.4) {\footnotesize$3$};
        \end{scope}

        \begin{scope}[very thick, shift={(0,4)}, yscale=-1]
            % \draw (1,0) -- (1.5,.5);
            % \draw (2,0) -- (1.5,.5);
            \draw (1.5,0) -- (1.5,4);

            \node at (.5,-.4) {\scriptsize$3$};
            \node at (1.5,-.4) {\scriptsize$2$};
            % \node at (2.5,-.4) {\scriptsize$1$};
        \end{scope}

        \node (a) at (.5,3.5) [d] {};
        %\node (c) at (.5,.5) [d] {};
        
        %\node (a1) at (1.5,3.5) [d] {};
        %\node (c1) at (1.5,.5) [d] {};
        
        \node[circle, draw=red, very thick, fill=white, inner sep=2pt] at (.5,2) {{$\color{red}2$}};
        \node[circle, draw=red, very thick, fill=white, inner sep=2pt] at (1.5,2) {{$\color{red}2$}};
        \end{scope}

        \draw[thick, -latex] (12.5,2) -- (12.8,2);

        % Fifth drawing
        \begin{scope}[very thick,shift={(13,0)}]
        \draw[] (.5,2) circle (.5cm);

        \node (a) at (0,2) [d] {};
        
        \node[circle, draw=red, very thick, fill=white, inner sep=2pt] at (1,2) {{$\color{red}4$}};
        \end{scope}
        \end{tikzpicture}
    \captionsetup{font=small, margin=18mm}
    
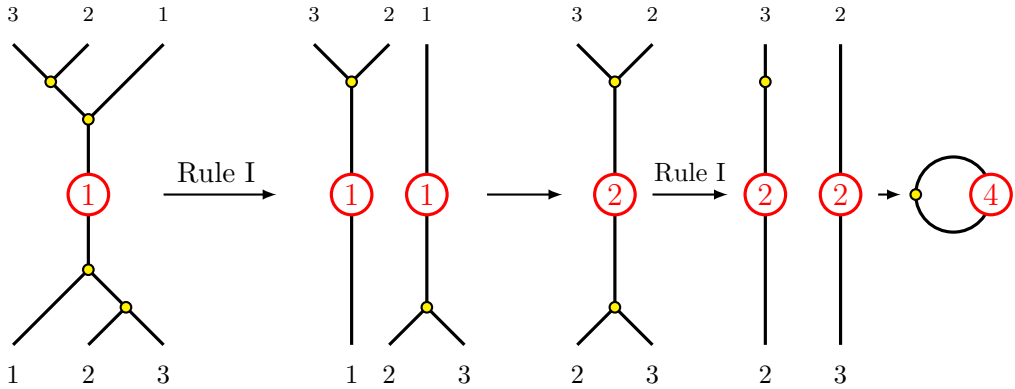
\captionof{figure}{Algorithm applied to $g_3$.}
\end{center}

\begin{exa}
    Consider $g$ such that $g @ 0 \colon00w \leftrightarrow 01w$ and $g @ 1= x_0^{-1}$. This element is clearly non-torsion and it has an attractor-repellor behavior on the right subtree, while it describes an orbit of period 2 on the left subtree, this can be verified by the algorithm (see \cref{fig:non-torsion}).
\end{exa}

\begin{center}
    \begin{tikzpicture}[scale=.95,
    d/.style={circle, draw=black, fill=yellow, inner sep=0mm, line width=.8pt, minimum size = 4pt}]
    % Element
    \begin{scope}[very thick, shift={(-5,2.5)}]
        \draw (0,0) -- (2,2) -- (4,0);
        \draw (1,0) -- (.5,.5);
        \draw (2,0) -- (3,1);
        \draw (3,0) -- (3.5,.5);

        \begin{scope}[rotate=180,shift={(-4,1)}]
            \draw (0,0) -- (2,2) -- (4,0);
            \draw (1,0) -- (1.5,.5);
            \draw (2,0) -- (1,1);
            \draw (3,0) -- (3.5,.5);
        \end{scope}

        \draw[blue] (2,0) -- (2,-1);
        \draw[blue] (3,0) -- (3,-1);
        \draw[blue] (4,0) -- (4,-1);
        \draw[purple] (1,0) -- (0,-1);
        \draw[line width=4pt,white] (0.3,-.3) -- (0.7,-.7);
        \draw[purple] (0,0) -- (1,-1);
    \end{scope}

    % Start of algorithm
        \begin{scope}[very thick, shift={(0,-1)}]
            \draw (0,0) -- (2,2) -- (4,0);
            \draw (1,0) -- (.5,.5);
            \draw (2,0) -- (3,1);
            \draw (3,0) -- (3.5,.5);

            \node at (0,-.4) {\footnotesize$1$};
            \node at (1,-.4) {\footnotesize$2$};
            \node at (2,-.4) {\footnotesize$3$};
            \node at (3,-.4) {\footnotesize$4$};
            \node at (4,-.4) {\footnotesize$5$};
        \end{scope}

        \begin{scope}[very thick, shift={(0,5)}, yscale=-1]
            \draw (0,0) -- (2,2) -- (4,0);
            \draw (1,0) -- (.5,.5);
            \draw (2,0) -- (3,1);
            \draw (3,0) -- (2.5,.5);

            \node at (0,-.4) {\scriptsize$2$};
            \node at (1,-.4) {\scriptsize$1$};
            \node at (2,-.4) {\scriptsize$3$};
            \node at (3,-.4) {\scriptsize$4$};
            \node at (4,-.4) {\scriptsize$5$};
        \end{scope}
        
        \draw[very thick] (2,1) -- (2,3);

        \node at (2.5,4.5) [d] {};
        \node at (3,4) [d] {};
        \node at (2,3) [d] {};
        \node at (.5,4.5) [d] {};

        \node at (3,0) [d] {}; %1
        \node at (2,1) [d] {}; %eps
        \node at (3.5,-.5) [d] {}; %11
        \node at (.5,-.5) [d] {}; %0
        
        \node[circle, draw=red, very thick, fill=white, inner sep=2pt] at (2,2) {{$\color{red}1$}};

        % End of algorithm
        \begin{scope}[shift={(4.5,0)}]
            \begin{scope}[very thick, shift={(0,-1)}]
            \draw (1.5,0) -- (1.5,5.5);
            \draw (2,0) -- (2.5,.5);
            \draw (3,0) -- (2.5,.5);

            \node at (1.5,-.4) {\footnotesize$3$};
            \node at (2,-.4) {\footnotesize$4$};
            \node at (3,-.4) {\footnotesize$5$};
        \end{scope}

        \begin{scope}[very thick, shift={(0,5)}, yscale=-1]
            \draw (1,0) -- (1.5,.5);
            \draw (2,0) -- (1.5,.5);
            \draw (2.5,0) -- (2.5,5.5);

            \node at (1,-.4) {\scriptsize$3$};
            \node at (2,-.4) {\scriptsize$4$};
            \node at (2.5,-.4) {\scriptsize$5$};
        \end{scope}

        \node (a) at (1.5,4.5) [d] {}; %a
        %\node (b) at (1.5,3) [d] {}; %b
        %\node (c) at (1.5,1) [d] {}; %c

        %\node (a1) at (2.5,3) [d] {}; %a1
        %\node (b1) at (2.5,1) [d] {}; %b1
        \node (c1) at (2.5,-.5) [d] {}; %c1

        \node[circle, draw=red, very thick, fill=white, inner sep=2pt] at (1.5,2) {{$\color{red}1$}};
        \node[circle, draw=red, very thick, fill=white, inner sep=2pt] at (2.5,2) {{$\color{red}1$}};
        \end{scope}

    \begin{scope}[very thick,shift={(4,0)}]
        \draw (.5,2) circle (.5cm);

        \node at (0,2) [d] {};
        
        \node[circle, draw=red, fill=white, inner sep=2pt] at (1,2) {{$\color{red}2$}};
    \end{scope}
\end{tikzpicture}
    \captionsetup{font=small}%, margin=18mm}
    
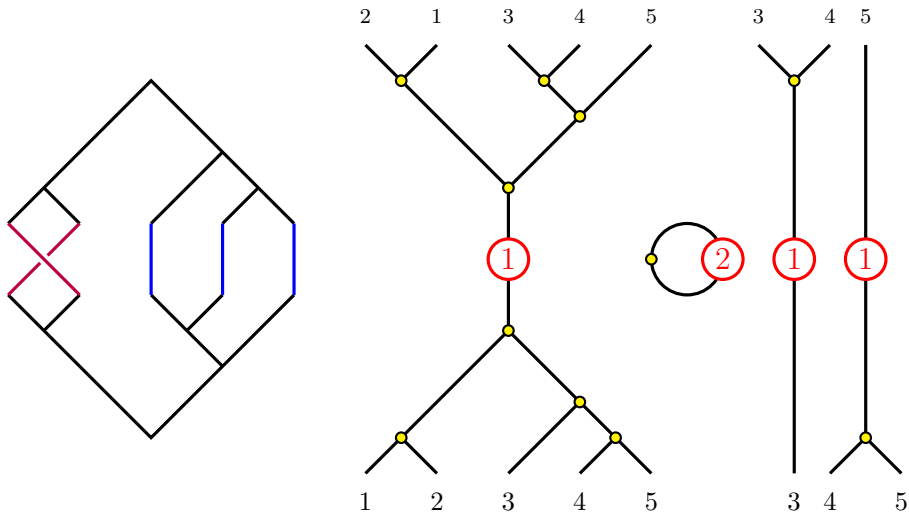
\captionof{figure}{An element of $V$ together with the first (center) and the last (right) step of algorithm. In this case there is a periodic orbit (with the period highlighted) and two strands representing the attractor-repellor behavior.}
    \label{fig:non-torsion}
\end{center}

\begin{prop}\label{prop:components&labels}
Starting with an element represented by a pair of trees with $n$ leaves (hence $n-1$ carets) each, the end diagram satisfies
\begin{itemize}[leftmargin=7mm]
    \item The number of connected components is bounded by $n$.
    \item The sum of all labels is bounded by $2^{n-1}$.
\end{itemize}
\end{prop}
\begin{proof}We first observe that the algorithm terminates after at most $n-1$ applications of Rule~I, since the number of carets start at $2n-2$ and decreases by $2$ at each operation. We observe that the number of connected components starts at $1$ and increase at most by $1$ at each application of Rule~I. Similarly, the sum of labels starts at $1$ and is at most doubled at each application.
\end{proof}

\begin{rem}[Complexity]
   The first two steps of the algorithm run in $\mathrm{DSPACE}(n^2)$ and $\mathrm{DTIME}(n^2)$. Indeed,
   \begin{itemize}[leftmargin=7mm]
       \item One needs to keep track of a graph with at most $2n$ vertices, their type (merge or split), and for each of them the target vertex of the $1$ or $2$ outgoing edges and their labels (which are $n$-digit integers).
       \item At each step, one needs to scan for a merge-split pair and do some $n$-digits additions, both operation taking $O(n)$ time.
   \end{itemize}
   The most expensive part of the algorithm seems to be to compute the least common multiple in part (3). Using recursively the formula
   \[ \lcm(l_1,l_2,\ldots,l_k)=\lcm\left(\lcm\bigl(l_1,\ldots,l_{[\frac k2]}\bigr),\lcm\bigl(l_{[\frac k2]+1},\ldots,l_k\bigr)\right), \]
    we need to compute $\frac n{2^{i+1}}$ least common multiples of two $2^in$-digit integers, for $i=0,\ldots,\log_2n$. This takes $O\bigl((n\log n)^2\bigr)$ operations using state-of-the-art algorithms for computing the product, quotient and $\gcd$ of integers \cite{compute_multiplication, compute_gcd}.
\end{rem}
%\todo{Is this Alex-good? Or just Alex-fine? (Complexity)}

\subsubsection{Correctness} \label{subsection:correctness}

In the previous subsection, we explain how to obtain reduced strand diagrams (meaning without merge-split pairs) which are even reduced \say{in the cylinder}. This reduction can be done using only marked strand diagrams all representing the same element of $V$, using rules III, III' and sometimes Rule~IV for simplifying merge-split pairs, depending how many times the edge linking the pair wraps around the cylinder (see \cref{fig:splitted-diagram}).
\begin{center}
    \begin{tikzpicture}[very thick,
	   d/.style={circle, draw=black, fill=yellow, inner sep=0mm, line width=.8pt, minimum size = 4pt},
	   s/.style={circle, draw=black, fill=black, inner sep=0mm, line width=.8pt, minimum size = 3pt}]
       
    \node at (-.8,1.5) {$g=$};
    \begin{scope}[shift={(0,0)}]
        \draw[thick] (0,0) rectangle (2,3);
        
        \node[label=above:$\varepsilon$] (source) at (1,3) [s] {};
        \node[label=below:$\varepsilon$] (sink) at (1,0) [s] {};
        \node (eps) at (1,2.4) [d] {};
        \node (1) at (1.4,2) [d] {};
        \node (eps') at (1,.6) [d] {};
        \node (1') at (1.4,1) [d] {};

        \draw (source) to (eps) to (1);
        \draw (1') to (eps') to (sink);
        \draw (1) to [out=-45, in = 135, looseness=2] (eps');
        \draw[white, line width=4pt] (1) to[out=-135, in=45, looseness=2] (1');
        \draw (1) to[out=-135, in=45, looseness=2] (1');
        \draw[white, line width=4pt] (eps) to [out=-135, in=135, looseness=2] (1');
        \draw (eps) to [out=-135, in=135, looseness=2] (1');
    \end{scope}
    \node at (2.75,1.5) {$=$};
    \begin{scope}[shift={(3.5,0)}]
        \draw[thick] (0,0) rectangle (2,3);
        
        \node[label=above:$0$] (source0) at (.5,3) [s] {};
        \node[label=above:$1$] (source1) at (1.5,3) [s] {};
        \node[label=below:$0$] (sink0) at (.5,0) [s] {};
        \node[label=below:$1$] (sink1) at (1.5,0) [s] {};
        \node (1) at (1.4,2) [d] {};
        \node (1') at (1.4,1) [d] {};

        \draw (source1) to[out=-90, in=90] (1);
        \draw (1) to [out=-45, in = 90] (sink0);
        
        \draw[white, line width=4pt]  (source0) to [out=-90, in=135] (1');
        \draw (source0) to [out=-90, in=135] (1');
        \draw[white, line width=4pt] (1') to[out=-135, in=90] (sink1);
        \draw (1') to[out=-90, in=90] (sink1);
        \draw[white, line width=4pt] (1) to[out=-135, in=45, looseness=2] (1');
        \draw (1) to[out=-135, in=45, looseness=2] (1');
    \end{scope} 
    \node at (6.25,1.5) {$=$};
    \begin{scope}[shift={(7,0)}]
        \draw[thick] (0,0) rectangle (2,3);
        
        \node[label=above:$0$] (source0) at (.5,3) [s] {};
        \node[label=above:$10$] (source10) at (1.25,3) [s] {};
        \node[label=above:$11$] (source11) at (1.75,3) [s] {};
        \node[label=below:$0$] (sink0) at (.5,0) [s] {};
        \node[label=below:$10$] (sink10) at (1.25,0) [s] {};
        \node[label=below:$11$] (sink11) at (1.75,0) [s] {};

        \draw (source11) to[out=-90, in=90] (sink0);
        \draw[white, line width=4pt] (source0) to[out=-90, in=90] (sink10);
        \draw (source0) to[out=-90, in=90] (sink10);
        \draw[white, line width=4pt] (source10) to[out=-90, in=90] (sink11);
        \draw (source10) to[out=-90, in=90] (sink11);
    \end{scope}  
    \end{tikzpicture}
    \captionsetup{font=small}
        \hspace*{10mm}
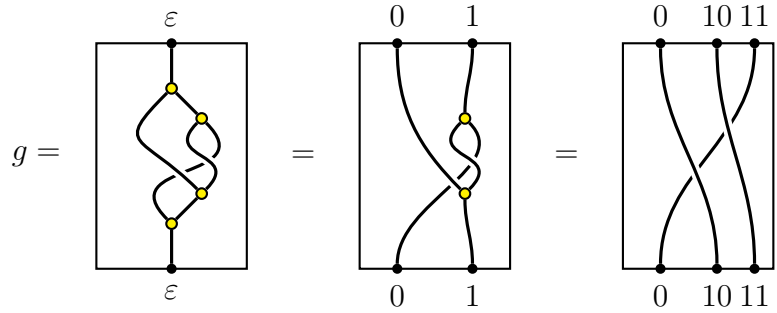
\captionof{figure}{The reduction to a diagram without merge-splits for an element $g \in T$.}
    \label{fig:splitted-diagram}
\end{center}
Once we have these reduced marked strand diagrams representing the element, understanding the dynamic of $g\acts \mathfrak C$ is much easier.
\begin{center}
    \begin{tikzpicture}[very thick,
	   d/.style={circle, draw=black, fill=yellow, inner sep=0mm, line width=.8pt, minimum size = 4pt},
	   s/.style={circle, draw=black, fill=black, inner sep=0mm, line width=.8pt, minimum size = 3pt}]
    \begin{scope}[shift={(0,0)}]
        \node at (-.4,1.5) {$g$};
        \draw[thick, dotted, red] (0,0) -- (2,0);
        \draw[thick, dotted, red] (0,3) -- (2,3);
        
        \node[label=above:$0$] (source0) at (.5,3) [s] {};
        \node[label=above:$10$] (source10) at (1.25,3) [s] {};
        \node[label=above:$11$] (source11) at (1.75,3) [s] {};
        \node (sink0) at (.5,0) [s] {};
        \node (sink10) at (1.25,0) [s] {};
        \node (sink11) at (1.75,0) [s] {};

        \draw (source11) to[out=-90, in=90] (sink0);
        \draw[white, line width=4pt] (source0) to[out=-90, in=90] (sink10);
        \draw (source0) to[out=-90, in=90] (sink10);
        \draw[white, line width=4pt] (source10) to[out=-90, in=90] (sink11);
        \draw (source10) to[out=-90, in=90] (sink11);

        \node[circle, draw=black, fill=red, inner sep=0mm, line width=.8pt, minimum size = 4pt] at (.502,2.8) {};

        \node[circle, draw=red, fill=white, inner sep=2pt] at (.6,2.3) {{\small$\color{red}3$}};
    \end{scope}
    \begin{scope}[shift={(0,-3)}]
        \node at (-.4,1.5) {$g$};
        \draw[thick, dotted, red] (0,0) -- (2,0);
        \draw[thick, dotted, red] (0,3) -- (2,3);
        
        \node (source0) at (.5,3) [s] {};
        \node (source10) at (1.25,3) [s] {};
        \node (source11) at (1.75,3) [s] {};
        \node (sink0) at (.5,0) [s] {};
        \node (sink10) at (1.25,0) [s] {};
        \node (sink11) at (1.75,0) [s] {};

        \draw (source11) to[out=-90, in=90] (sink0);
        \draw[white, line width=4pt] (source0) to[out=-90, in=90] (sink10);
        \draw (source0) to[out=-90, in=90] (sink10);
        \draw[white, line width=4pt] (source10) to[out=-90, in=90] (sink11);
        \draw (source10) to[out=-90, in=90] (sink11);

        \node[circle, draw=black, fill=red, inner sep=0mm, line width=.8pt, minimum size = 4pt] at (.502,2.8) {};
    \end{scope}  
    \begin{scope}[shift={(0,-6)}]
        \node at (-.4,1.5) {$g$};
        \draw[thick, dotted, red] (0,0) -- (2,0);
        \draw[thick, dotted, red] (0,3) -- (2,3);
        
        \node (source0) at (.5,3) [s] {};
        \node (source10) at (1.25,3) [s] {};
        \node (source11) at (1.75,3) [s] {};
        \node[circle, draw=black, fill=red, inner sep=0mm, line width=.8pt, minimum size = 4pt] (node) at (.502,-.2) {};
        \draw (.5,0) -- (node);
        \node[label=below:$0$] (sink0) at (.5,0) [s] {};
        \node[label=below:$10$] (sink10) at (1.25,0) [s] {};
        \node[label=below:$11$] (sink11) at (1.75,0) [s] {};

        \draw (source11) to[out=-90, in=90] (sink0);
        \draw[white, line width=4pt] (source0) to[out=-90, in=90] (sink10);
        \draw (source0) to[out=-90, in=90] (sink10);
        \draw[white, line width=4pt] (source10) to[out=-90, in=90] (sink11);
        \draw (source10) to[out=-90, in=90] (sink11);

        \node[circle, draw=black, fill=red, inner sep=0mm, line width=.8pt, minimum size = 4pt] at (.502,2.8) {};
    \end{scope}
    \end{tikzpicture} \hspace*{10mm}
    \begin{tikzpicture}[scale=.83,
        d/.style={circle, draw=black, fill=yellow, inner sep=0mm, line width=.8pt, minimum size = 4pt},
        s/.style={circle, draw=black, fill=black, inner sep=0mm, line width=.8pt, minimum size = 3pt}]
        % Right drawing
        \begin{scope}[very thick, shift={(0,4.8)}]
            \draw[thick, dotted, red] (1,1.2) -- (4,1.2);
        
            \draw (3.5,.5) -- (3.5,1.2);
            \draw (3,0) -- (3.5,.5);
            \draw (4,0) -- (3.5,.5);

            \draw[out=-90, in=90] (1.5,1.2) to (2.5,-4.2); 
            \draw[out=-90, in=135] (2.5,1.2) to (1,-3);
            \draw[out=-135, in=45] (3,0) to (2,-3);
            \draw[out=-45, in=90] (4,0) to (3.5,-4.2);

            \node[label=above:$0$] (0) at (1.5,1.2) [s] {};
            \node[label=above:$10$] (10) at (2.5,1.2) [s] {};
            \node[label=above:$11$] (11) at (3.5,1.2) [s] {};
            \node (11split) at (3.5,.5) [d] {};

            \node[circle, draw=red, fill=white, inner sep=2pt] at (1.55,.6) {{\small$\color{red}2$}};
            \node[circle, draw=red, fill=white, inner sep=2pt] at (4.15,-.2) {{\small$\color{red}1$}};
            \node[circle, draw=red, fill=white, inner sep=2pt] at (2.85,-.2) {{\small$\color{red}0$}};
        \end{scope}

        \begin{scope}[very thick, shift={(0,1.8)},yscale=-1]
            \draw (1.5,.5) -- (1.5,1.2);
            \draw (1,0) -- (1.5,.5);
            \draw (2,0) -- (1.5,.5);

            \node (0) at (1.5,1.2) [s] {};
            \node (10) at (2.5,1.2) [s] {};
            \node (11) at (3.5,1.2) [s] {};
            \node (0split) at (1.5,.5) [d] {};
        \end{scope}
        % Right drawing
        \begin{scope}[shift={(0,-5.4)}]
        \begin{scope}[very thick, shift={(0,4.8)}]
            \draw[thick, dotted, red] (1,1.2) -- (4,1.2);
            \draw (3.5,.5) -- (3.5,1.2);
            \draw (3,0) -- (3.5,.5);
            \draw (4,0) -- (3.5,.5);

            \draw[out=-90, in=90] (1.5,1.2) to (2.5,-4.2); 
            \draw[out=-90, in=135] (2.5,1.2) to (1,-3);
            \draw[out=-135, in=45] (3,0) to (2,-3);
            \draw[out=-45, in=90] (4,0) to (3.5,-4.2);

            \node (0) at (1.5,1.2) [s] {};
            \node (10) at (2.5,1.2) [s] {};
            \node (11) at (3.5,1.2) [s] {};
            \node (11split) at (3.5,.5) [d] {};

            %\node[] at (1.5,1.5) {$0$};
            %\node[] at (2.5,1.5) {$10$};
            %\node[] at (3.5,1.5) {$11$};
        \end{scope}

        \begin{scope}[very thick, shift={(0,1.8)},yscale=-1]
            \draw[thick, dotted, red] (1,1.2) -- (4,1.2);
            \draw (1.5,.5) -- (1.5,1.2);
            \draw (1,0) -- (1.5,.5);
            \draw (2,0) -- (1.5,.5);

            \node[label=below:$0$] (0) at (1.5,1.2) [s] {};
            \node[label=below:$10$] (10) at (2.5,1.2) [s] {};
            \node[label=below:$11$] (11) at (3.5,1.2) [s] {};
            \node (0split) at (1.5,.5) [d] {};

            %\node[] at (1.5,1.5) {$0$};
            %\node[] at (2.5,1.5) {$10$};
            %\node[] at (3.5,1.5) {$11$};
        \end{scope}
        \end{scope}
    \end{tikzpicture}
    \captionsetup{font=small}
        \hspace*{10mm}
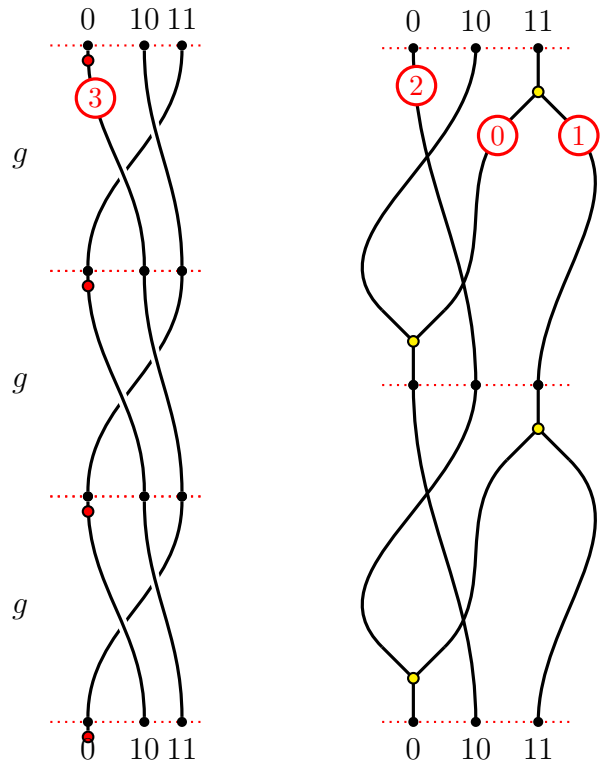
\captionof{figure}{On the left, an element belonging to $T$ of period $3$. On the right, a non-torsion element of $V$.}
        \label{fig:orbits}
\end{center}
We note that the label on each edge is exactly the number of horizontal separations it crosses. (For \say{circle} components, we add a dummy vertex in red to mark the end of the edge.) We distinguish two cases:

\begin{itemize}[leftmargin=7mm]
    \item If there is no merge and split left, then the element is just a permutation of a partition of $\mathfrak C$. For each part $\alpha$ of the partition, the label $l_i$ on the associated edge is the smallest integer such that $g^l|_\alpha=\id$. Therefore, the order of $g$ is the least common multiple of the labels $l_i$. See the left side of \cref{fig:orbits}.

    \item If there are some merges and splits left, then starting from any merge and going down, one can only meet other merges. Since there are only finitely many vertices, some merge must repeat: this is an attractor cycle, in particular the order is infinite. This attractor contains a point with a finite orbit, whose length is the sum of the labels on the cycle. See the right side of \cref{fig:orbits}.
\end{itemize}
%%%%%%%%%%%%%%%%%%

\subsection{Rotation number in Thompson \texorpdfstring{$T$}{T}}
\label{sec:compute-rotation}

Let $g\in\Homeo^+(\bS)$ be a homeomorphism of the unit circle and let $\hat{g}$ be a lift of $g$ to $\Homeo^+(\mathbb{R})$. The rotation number of $g$ is defined as
$$\rot(g)=\lim_{i \to\infty} \dfrac{\hat{g}^i(x)}{i} \mod \Z.$$
This number is independent of the choice of $x\in\R$, and also of the choice of lift thanks to the congruence $\mod \Z$. We recall that $g$ admits a fixed point if and only if $\rot(g)=0$. In particular, if $\rot(g)\in\Q$, then its denominator $k$ is the smallest natural number such that $g^k$ admits a fixed point. Indeed
$$\rot(g^k)=k\cdot \rot(g)=0.$$

A remarkable result of Ghys and Sergiescu states that rotation numbers of elements of Thompson $T$ are always rational numbers \cite{Ghys_rotation}, see also \cite{Liousse_rotation}. In this subsection, we explain a variation on our main algorithm which computes $\rot(g)$ efficiently for $g\in T$. Note that the main algorithm already gave the length of finite orbits of $g\acts\mathbb S$, i.e., the denominator of $\rot(g)$. We need one more set of labels to keep track of the numerator.

\subsubsection{Algorithm}

We describe an algorithm which computes efficiently the rotation number $\rot(g)\in\Q/\Z$ of an element $g\in T$.

\begin{enumerate}[leftmargin=7mm, label=(\arabic*)]
    \item  We write the element on a torus and write a red label $1$ on each edge crossing the red \say{seam}, and a blue label $1$ on each edge crossing the blue \say{seam}, as in \cref*{fig:algo_rot_step1}.

\begin{center}
    \begin{tikzpicture}
        % Element
        \begin{scope}[very thick, shift={(-6,3.3)}]
            \draw (1,0) -- (2.5,1.5);
            \draw (2,0) -- (3,1);
            \draw (3,0) -- (2.5,.5);
            \draw (4,0) -- (2.5,1.5);

            \draw[very thick] (2.5,1.5) -- (2.5,2.2);
            \draw[thick, dotted, red] (.5,2.2) -- (4.5,2.2);
            \node[circle, draw=red, fill=white, inner sep=2pt] at (2.5,2.2) {{$\color{red}1$}};

            \begin{scope}
                \clip (.5,-3.3) rectangle (4.5,2.2);
                \draw[ForestGreen] (1,0) -- (3,-1.5);
                \draw[ForestGreen] (2,0) -- (4,-1.5);
                \draw[ForestGreen] (3,0) -- (5,-1.5);
                \draw[ForestGreen] (4,0) -- (6,-1.5);
                \draw[ForestGreen] (-1,0) -- (1,-1.5);
                \draw[ForestGreen] (0,0) -- (2,-1.5);
            \end{scope}

            \draw[thick, dotted, blue] (.5,2.2) -- (.5,-3.3);
            \draw[thick, dotted, blue] (4.5,2.2) -- (4.5,-3.3);
            \node[circle, draw=blue, fill=white, inner sep=2pt] at (4.5,-0.375) {{\footnotesize$\color{blue}1$}};
            \node[circle, draw=blue, fill=white, inner sep=2pt] at (4.5,-1.125) {{\footnotesize$\color{blue}1$}};     
        \end{scope}

        \begin{scope}[very thick, shift={(-6,1.8)}, yscale=-1]
            \draw (1,0) -- (2.5,1.5);
            \draw (2,0) -- (1.5,.5);
            \draw (3,0) -- (3.5,.5);
            \draw (4,0) -- (2.5,1.5);
            
            \draw[very thick] (2.5,1.5) -- (2.5,1.8);
            \draw[thick, dotted, red] (.5,1.8) -- (4.5,1.8);
        \end{scope}

        % Start of algorithm
        \begin{scope}[very thick, shift={(0,0)}]
            \draw (1,0) -- (2.5,1.5);
            \draw (2,0) -- (3,1);
            \draw (3,0) -- (2.5,.5);
            \draw (4,0) -- (2.5,1.5);

            \node at (1,-.4) {\footnotesize$1$};
            \node at (2,-.4) {\footnotesize$2$};
            \node at (3,-.4) {\footnotesize$3$};
            \node at (4,-.4) {\footnotesize$4$};

            \node[circle, draw=blue, fill=white, inner sep=1.5pt] at (2.75,.25) {{\footnotesize$\color{blue}1$}};
            \node[circle, draw=blue, fill=white, inner sep=1.5pt] at (3.5,.5) {{\footnotesize$\color{blue}1$}};
        \end{scope}

        \begin{scope}[very thick, shift={(0,5.5)}, yscale=-1]
            \draw (1,0) -- (2.5,1.5);
            \draw (2,0) -- (1.5,.5);
            \draw (3,0) -- (3.5,.5);
            \draw (4,0) -- (2.5,1.5);

            \node at (1,-.4) {\footnotesize$3$};
            \node at (2,-.4) {\footnotesize$4$};
            \node at (3,-.4) {\footnotesize$1$};
            \node at (4,-.4) {\footnotesize$2$};
        \end{scope}

        \draw[very thick] (2.5,1.5) -- (2.5,4);
        \node[circle, draw=red, very thick, fill=white, inner sep=2pt] at (2.5,2.75) {{$\color{red}1$}};
        \end{tikzpicture}
        \captionsetup{font=small}
        \hspace*{10mm}
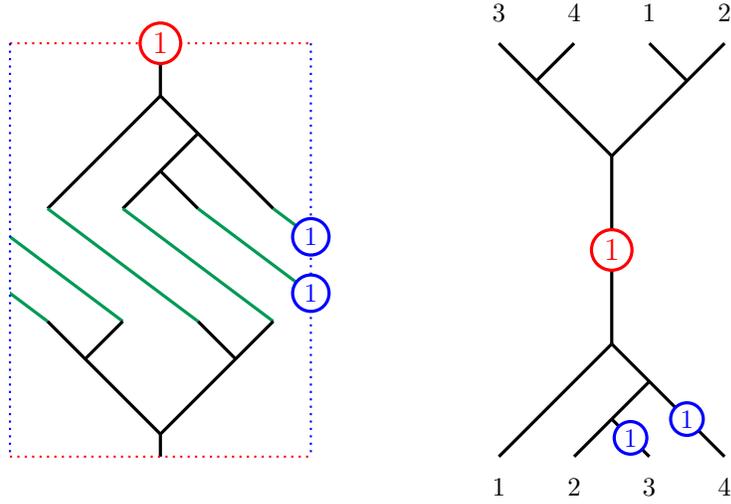
\captionof{figure}{The ideal picture on a torus, and the pratical picture.}
        \label{fig:algo_rot_step1}
    \end{center}
    \item If we find a merge-split pair of successive carets, we apply Rule~I as in the algorithm described in \S\ref{subsection:algorithm} for each color. Repeat until there are no more merge-split pairs of successive carets.
    \item At the end, the rotation number is $\frac ab$ where $a$ (resp.\ $b$) is the sum of blue (resp.\ red) labels along any oriented cycle of the final diagram.
\end{enumerate}

\begin{exa}
    Let us finish the example in \cref{fig:algo_rot_step1}.
    \begin{center}
    \begin{tikzpicture}
        % Start of algorithm
        \begin{scope}[very thick, shift={(0,0)}]
            \draw (1,0) -- (2.5,1.5);
            \draw (2,0) -- (3,1);
            \draw (3,0) -- (2.5,.5);
            \draw (4,0) -- (2.5,1.5);

            \node at (1,-.4) {\footnotesize$1$};
            \node at (2,-.4) {\footnotesize$2$};
            \node at (3,-.4) {\footnotesize$3$};
            \node at (4,-.4) {\footnotesize$4$};

            \node[circle, draw=blue, fill=white, inner sep=1.5pt] at (2.75,.25) {{\footnotesize$\color{blue}1$}};
            \node[circle, draw=blue, fill=white, inner sep=1.5pt] at (3.5,.5) {{\footnotesize$\color{blue}1$}};

            \draw[very thick] (2.5,1.5) -- (2.5,4);
            \node[circle, draw=red, very thick, fill=white, inner sep=2pt] at (2.5,2.75) {{$\color{red}1$}};
        \end{scope}

        \begin{scope}[very thick, shift={(0,5.5)}, yscale=-1]
            \draw (1,0) -- (2.5,1.5);
            \draw (2,0) -- (1.5,.5);
            \draw (3,0) -- (3.5,.5);
            \draw (4,0) -- (2.5,1.5);

            \node at (1,-.4) {\footnotesize$3$};
            \node at (2,-.4) {\footnotesize$4$};
            \node at (3,-.4) {\footnotesize$1$};
            \node at (4,-.4) {\footnotesize$2$};
        \end{scope}

        \draw[very thick, -latex] (4,2.75) -- (5.25,2.75);
        % Second step
        \begin{scope}[very thick, shift={(4.5,0)}]
            \draw (1,0) -- (2,1);
            \draw (2,0) -- (1.5,.5);
            \draw (3,0) -- (2,1);

            \node at (1,-.4) {\footnotesize$2$};
            \node at (2,-.4) {\footnotesize$3$};
            \node at (3,-.4) {\footnotesize$4$};

            \node[circle, draw=blue, fill=white, inner sep=1.5pt] at (1.75,.25) {{\footnotesize$\color{blue}1$}};
            \node[circle, draw=blue, fill=white, inner sep=1.5pt] at (2.5,.5) {{\footnotesize$\color{blue}1$}};

            \draw[very thick] (2,1) -- (2,3.5);
            \node[circle, draw=red, very thick, fill=white, inner sep=2pt] at (2,2.75) {{$\color{red}1$}};
        \end{scope}

        \begin{scope}[very thick, shift={(4.5,5.5)}, yscale=-1]
            \draw (1,0) -- (1.5,.5) -- (2,0);
            \draw (1.5,0.5) -- (1.5,1.5);
            \draw (1.5,1.5) -- (2,2) -- (2.5,1.5);
            
            \node at (1,-.4) {\footnotesize$3$};
            \node at (2,-.4) {\footnotesize$4$};
            \node at (2.5,1.1) {\footnotesize$2$};

            \node[circle, draw=red, very thick, fill=white, inner sep=2pt] at (1.5,1) {{\footnotesize$\color{red}1$}};
        \end{scope}

        \draw[very thick, -latex] (7.5,2.75) -- (8.75,2.75);
        % Third step
        \begin{scope}[very thick, shift={(8,0)}]
            \draw (1,0) -- (1.5,.5);
            \draw (2,0) -- (1.5,.5);

            \node at (1,-.4) {\footnotesize$2$};
            \node at (2,-.4) {\footnotesize$3$};

            \node[circle, draw=blue, fill=white, inner sep=1.5pt] at (1.75,.25) {{\footnotesize$\color{blue}1$}};

            \draw[very thick] (1.5,.5) -- (1.5,3.5);
            \node[circle, draw=red, very thick, fill=white, inner sep=2pt] at (1.5,2.75) {{$\color{red}2$}};
        \end{scope}

        \begin{scope}[very thick, shift={(8,5.5)}, yscale=-1]
            \draw (1,1.5) -- (1.5,2) -- (2,1.5) -- (2,0);
            
            \node at (1,1.1) {\footnotesize$3$};
            \node at (2,-.4) {\footnotesize$2$};

            \node[circle, draw=red, very thick, fill=white, inner sep=2pt] at (2,.45) {{\footnotesize$\color{red}1$}};
            \node[circle, draw=blue, fill=white, inner sep=1.5pt] at (2,1.05) {{\footnotesize$\color{blue}1$}};
        \end{scope}

        \draw[very thick, -latex] (10.125,2.75) -- (11.375,2.75);
        % Third step
        \begin{scope}[very thick, shift={(11,0)}]
            \node at (1,-.4) {\footnotesize$2$};
            \node at (1,5.9) {\footnotesize$2$};

            \draw[very thick] (1,0) -- (1,5.5);
            \node[circle, draw=red, very thick, fill=white, inner sep=2pt] at (1,2.75) {{$\color{red}5$}};
            \node[circle, draw=blue, very thick, fill=white, inner sep=2pt] at (1,2) {{$\color{blue}2$}};
        \end{scope}
    \end{tikzpicture}
    \captionsetup{font=small, margin=10mm}
    \captionof{figure}{The algorithm applied to $\tilde a_1\in T$ from \cref*{thm:period_T}.}
\end{center}
Therefore, this element has rotation number $\frac25$.
\end{exa}

\subsubsection{Correctness}

Just as in \cref{subsection:correctness}, we can get to a reduced diagram using only marked strand diagrams representing the same element $g\in T$. Moreover, we can concatenate horizontally copies of this diagram to get a marked diagram for a lift $\hat g\in\Homeo^+(\R)$ (see \cref*{fig:algo_rot_corr}). The labels on an edge are exactly the number of vertical and horizontal separations it crosses.
\begin{center}
    \import{Pictures/}{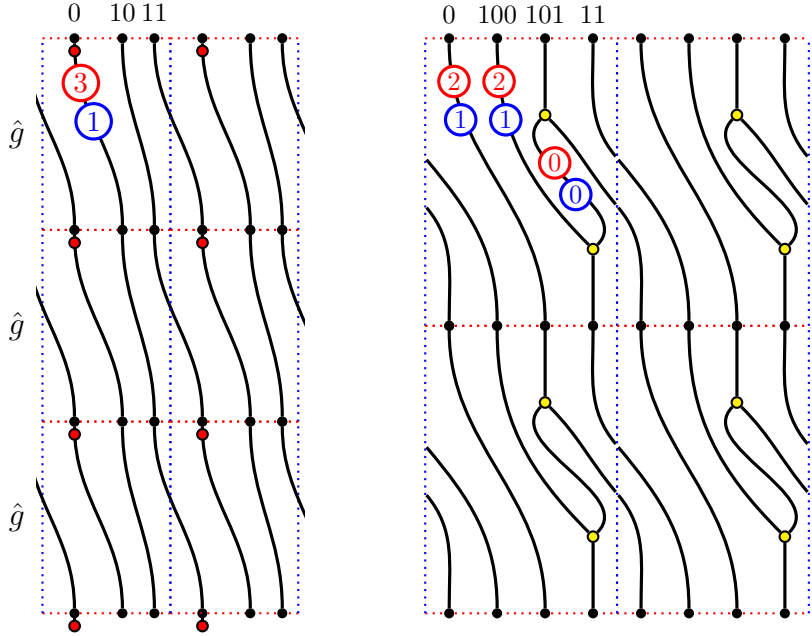}
    \captionsetup{font=small, margin=10mm}
    \captionof{figure}{The reduced diagrams of two elements of $T$, with the corresponding labels, concatenated to get diagrams of $\hat g\in\Homeo^+(\R)$. The rotation numbers are $\frac13$ and $\frac12$ respectively.}
    \label{fig:algo_rot_corr}
\end{center}
We observe that every oriented path between two copies of the same vertex gives a point $x\in \bS$ such that $g^b(x)=x$, where $b$ is the sum of the red labels on the path. Therefore, we have a point $\hat x\in\R$ such that $\hat g^b(\hat x)=\hat x+a$ for some integer $a\in\Z$, more specifically $a$ is the sum of the blue labels on the path. We conclude that $\rot(g)=\frac ab$.

\section{Thompson groups}
\subsection{Thompson group \texorpdfstring{$T$}{}}
We start by computing the period growth of Thompson group $T$.
\begin{thm} \label{thm:period_T}
The period growth of Thompson $T$ satisfies $p_T(n)\asymp \exp(n)$.
\end{thm}
\begin{proof}
A torsion element $g\in T$ gives a number of connected components with the same label at the end of the algorithm described in \S\ref{subsection:algorithm}. Indeed, if a cycle has length $l$, the element $g^l$ fixes a point, hence it belongs to $F$ up to conjugation. It is therefore the trivial element, and all the cycles need to have length $l$. If $g$ has $n$ leaves (hence by \cref{thm:metric_all} of length $\asymp n$), its order is the label on each of the components, hence is bounded by $2^{n-1}$ by \cref{prop:components&labels}. This gives the upper bound.
\medskip

To prove that $p_T(n)$ is exponential, we observe that, starting from any torsion element $g$ satisfying $g(\frac12)=0$ (i.e., the vertices sent back to the left are exactly the vertices below $1$) with $n$ leaves, we can construct an element $\tilde{g}$ with $n+2$ leaves and $\ord(\tilde{g})\geq 2 \ord(g)$, as in \cref*{fig:Tlow}. Moreover, $\tilde g^{-1}(\frac12)=0$.
\begin{center}
    \begin{tikzpicture}[scale = .85]   
    
    \begin{scope}[very thick, shift={(10.5,0)}]
     \begin{scope}[shift={(0,.5)}]
        \node at (.2,1.5) {$\tilde{g}$};
        \node[purple] at (4.5,.2) {\small A};
        \node[blue] at (2.5,.2) {\small B};
        
        \draw (0,0) -- (2.5,2.5) -- (4.5,.5);
        \draw (1,0) -- (3,2);
        \draw (2.5,.5) -- (2,1);
        \draw[blue] (2,0) -- (2.5,.5);
        \draw[blue] (2,0) -- (3,0);
        \draw[blue] (3,0) -- (2.5,.5);
        \draw[purple] (4,0) -- (4.5,.5);
        \draw[purple] (4,0) -- (5,0);
        \draw[purple] (4.5,.5) -- (5,0);
    \end{scope}

        \draw[blue] (2,.5) -- (0,-1.3);
        \draw[blue] (3,.5) -- (1,-1.3);
        \draw[purple] (4,.5) -- (2,-1.3);
        \draw[purple] (5,.5) -- (3,-1.3);
        \draw[line width=3pt,white](1.444,.3) -- (4.8,-1.2);
        \draw[line width=3pt,white](0.444,.3) -- (3.8,-1.2);
        \draw[ForestGreen] (1,.5) -- (5,-1.3);
        \draw[ForestGreen] (0,.5) -- (4,-1.3);

    \begin{scope}[rotate=180,shift={(-5,1.3)}]
        \node[] at (3.5,.3) {\small C};

        \draw (0,0) -- (2.5,2.5) -- (5,0);
        \draw[] (1,0) -- (.5,.5);
        \draw[] (2,0) -- (3.5,1.5);
        \draw[] (2,0) -- (5,0);
    \end{scope}

    \end{scope}

    \begin{scope}[very thick, shift={(0,0)}]
        \node at (.4,1.2) {$g$};
        \node[purple] at (.5,.2) {\small A};
        \node[blue] at (2.5,.2) {\small B};
        \node[] at (1.5,-1.1) {\small C};

        \draw[purple] (0,0) -- (.5,.5);
        \draw[blue] (2.5,.5) -- (3,0);
        \draw (.5,.5) -- (1.5,1.5) -- (2.5,.5);
        \draw[purple] (1,0) -- (.5,.5);
        \draw[purple] (0,0) -- (1,0);
        \draw[blue] (2,0) -- (2.5,.5);
        \draw[blue] (2,0) -- (3,0);

        \draw[] (0,-.8) -- (1.5,-2.3) -- (3,-.8) -- (0,-.8);
        
        \draw[blue] (2,0) -- (0,-.8);
         \draw[blue] (3,0) -- (1,-.8);
           \draw[line width=3pt,white](1.5,-.2) -- (1,-.4);
         \draw[line width=4pt,white](2,-.4) -- (1.5,-.6) ;
        \draw[line width=4pt,white](0.3,-.12) -- (1.7,-.68);
        \draw[line width=4pt,white](1.3,-.12) -- (2.7,-.68);
         \draw[purple] (0,0) -- (2,-.8);
         \draw[purple] (1,0) -- (3,-.8);
         
    \end{scope}

\begin{scope}[very thick, shift={(5.5,-2.5)}]
        \node[purple] at (.5,.2) {\small A};
        \node[blue] at (2.5,.2) {\small B};

        \draw[purple] (0,0) -- (.5,.5);
        \draw[blue] (2.5,.5) -- (3,0);
        \draw (.5,.5) -- (1.5,1.5) -- (2.5,.5);
        \draw[purple] (1,0) -- (.5,.5);
        \draw[purple] (0,0) -- (1,0);
        \draw[blue] (2,0) -- (2.5,.5);
        \draw[blue] (2,0) -- (3,0);

\begin{scope}[very thick, shift={(0,5)}]
        \draw[] (0,-.8) -- (1.5,-2.3) -- (3,-.8) -- (0,-.8);
        \draw (1.5,-2.3) -- (1.5,-3.5);
        \node[] at (1.5,-1.4) {\small C};
 \end{scope}       

 \node[circle, draw=red, very thick, fill=white, inner sep=2pt] at (1.5,2) {{$\color{red}x$}};

 \node[circle, draw=red, very thick, fill=white, inner sep=2pt] at (1,1) {{$\color{red}y$}};
    
\end{scope}

\end{tikzpicture}
    \captionsetup{font=small, margin=5mm}
    
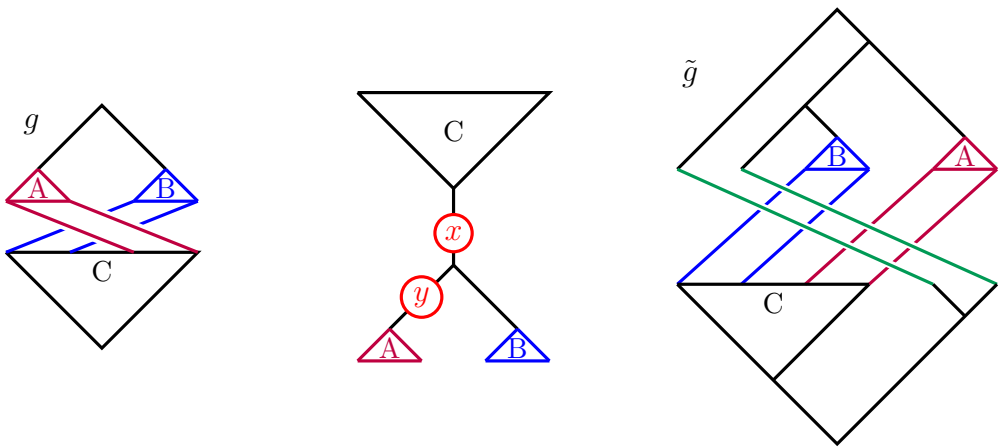
\captionof{figure}{At the first stage of the algorithm the element $g$ (on the left) provides the central diagram with $x=1$ and $y=0$ while, after two applications of Rule~I to the element $\tilde{g}$, the diagram is the same with $x=2$ and $y=1$.}
    \label{fig:Tlow}
\end{center}

Starting from $a_1\colon 0w\leftrightarrow 1w$ (with $2$ leaves and $\ord(a_1)=2$), we define recursively $a_{n+1}=\tilde a_n^{-1}$. By induction, the element $a_n$ has $2n$ leaves (hence $\norm{a_n}_T\asymp n$ by \cref{thm:metric_all}(b)) and $\ord(a_n)\ge 2^n$.
\end{proof}

\subsection{A question of D. Calegari} \label{subsection:Calegari}

Another natural question in Thompson $T$ is to bound the denominator of $\rot(g)$, for elements that are not necessarily torsion. This question is partially addressed in an article of D.\ Calegari \cite{calegari2007denominator}. For a reduced element $g=(\hspace{1pt}\Uc,\sigma,\Vc)\in T$, he proved that
\[ \text{denominator of }\rot(g) \le m(g)\cdot 2^{m(g)} \]
where $m(g)\coloneqq 2^{d(g)}$ and $d(g)$ is the maximal depth between $\Uc$ and $\Vc$. He provided examples where the denominator is proportional to $m^2$, and asked whether a polynomial (or even quadratic) upper bound in $m$ holds in full generality \cite[Question 4.4]{calegari2007denominator}.

We answer negatively to Calegari's question, using another family of torsion elements. In particular, the denominator of $\rot(g)$ coincides with $\ord(g)$ (see the proof of \cref{thm:period_T}). We now optimise the order for a given depth instead of a given number of leaves/word length.
\begin{center}
    \begin{tikzpicture}[scale = .9]
    \begin{scope}[very thick, shift={(0,0)}]
        \node at (.4,1.2) {\large$g$};
        \node at (1.5,.6) {$\Uc$};
        \node at (1.5,-1.4) {$\Vc$};
        \node at (.5,.25) {$\ell$};
        \draw[latex-latex] (3.5,0) -- (3.5,1.5);
        \node at (3.75,.8) {$d$};
        
        \draw[blue] (2,0) -- (0,-.8);
        \draw[blue] (3,0) -- (1,-.8);
        \draw[line width=4pt,white] (0,0) -- (2,-.8);
        \draw[line width=4pt,white] (1,0) -- (3,-.8);

        \draw (0,0) --  (1.5,1.5) -- (3,0) -- cycle;
        \draw (0,-.8) -- (1.5,-2.3) -- (3,-.8) -- cycle;
        
        \draw[purple] (0,0) -- (2,-.8);
        \draw[purple] (1,0) -- (3,-.8);
        \draw[purple] (0,0) -- (1,0);
        \draw[purple] (2,-.8) -- (3,-.8);
        \draw[blue] (2,0) -- (3,0);
        \draw[blue] (0,-.8) -- (1,-.8);  
    \end{scope}
    
    \begin{scope}[very thick, shift={(5.5,-.1)}, scale=.6]
        \draw[blue] (10,.5) -- (0,-1.3);
        \draw[blue] (11,.5) -- (1,-1.3);
        \draw[purple] (12,.5) -- (2,-1.3);
        \draw[purple] (14,.5) -- (3,-1.3);

        \draw[line width=3pt,white] (0,.5) -- (4,-1.3);
        \draw[ForestGreen] (0,.5) -- (4,-1.3);
        \draw[line width=3pt,white] (4,.5) -- (4.25,-1.3);
        \draw[ForestGreen] (4,.5) -- (4.25,-1.3);
        \draw[line width=3pt,white] (4.25,.5) -- (4.5,-1.3);
        \draw[ForestGreen] (4.25,.5) -- (4.5,-1.3);
        \draw[line width=3pt,white] (4.5,.5) -- (4.75,-1.3);
        \draw[ForestGreen] (4.5,.5) -- (4.75,-1.3);
        \draw[white] (5.2,-.4) circle (2.5pt);
        \draw[ForestGreen] (5.2,-.4) circle (1pt);
        \draw[white] (5.5,-.4) circle (2.5pt);
        \draw[ForestGreen] (5.5,-.4) circle (1pt);
        \draw[white] (5.8,-.4) circle (2.5pt);
        \draw[ForestGreen] (5.8,-.4) circle (1pt);
        \draw[line width=3pt,white] (6.25,.5) -- (6.5,-1.3);
        \draw[ForestGreen] (6.25,.5) -- (6.5,-1.3);
        \draw[line width=3pt,white] (6.5,.5) -- (6.75,-1.3);
        \draw[ForestGreen] (6.5,.5) -- (6.75,-1.3);
        \draw[line width=3pt,white] (6.75,.5) -- (7,-1.3);
        \draw[ForestGreen] (6.75,.5) -- (7,-1.3);
        \draw[line width=3pt,white] (7,.5) -- (8,-1.3);
        \draw[ForestGreen] (7,.5) -- (8,-1.3);
        \draw[line width=3pt,white] (8,.5) -- (12,-1.3);
        \draw[purple] (8,.5) -- (12,-1.3);
        \draw[line width=3pt,white] (9,.5) -- (14,-1.3);
        \draw[purple] (9,.5) -- (14,-1.3);
        
     \begin{scope}[shift={(0,.5)}]
        \node at (2,4.5) {\large $\tilde g$};
        \node at (5.5,.6) {$\Tc_{2^d}$};
        \node at (9.5,.6) {$\Uc$};
        \node at (13,.4) {\footnotesize$\Tc_\ell$};
        
        \draw (0,0) -- (7,7) -- (14,0);
        \draw (7,0) -- (3.5,3.5);
        \draw (4,0) -- (5.5,1.5) -- (7,0) -- (4,0);
        \draw (8,0) -- (9.5,1.5) -- (11,0) -- (8,0);
        \draw (9.5,1.5) -- (11,3);
        \draw (12,0) -- (13,1) -- (14,0) -- (12,0);
    \end{scope}

    \begin{scope}[yscale=-1, shift={(0,1.3)}]
        \node at (1.5,.5) {$\Vc$};
        \node at (5.5,.5) {$\Tc_{2^d}$};
        \node at (13,.4) {\footnotesize$\Tc_\ell$};
        
        \draw (0,0) -- (7,7) -- (14,0);
        \draw (0,0) -- (1.5,1.5) -- (3,0) -- (0,0);
        \draw (7,0) -- (3.5,3.5);
        \draw (4,0) -- (5.5,1.5) -- (7,0) -- (4,0);
        \draw (8,0) -- (11,3);
        \draw (12,0) -- (13,1) -- (14,0) -- (12,0);
    \end{scope}

    \end{scope}

\end{tikzpicture} \\
    
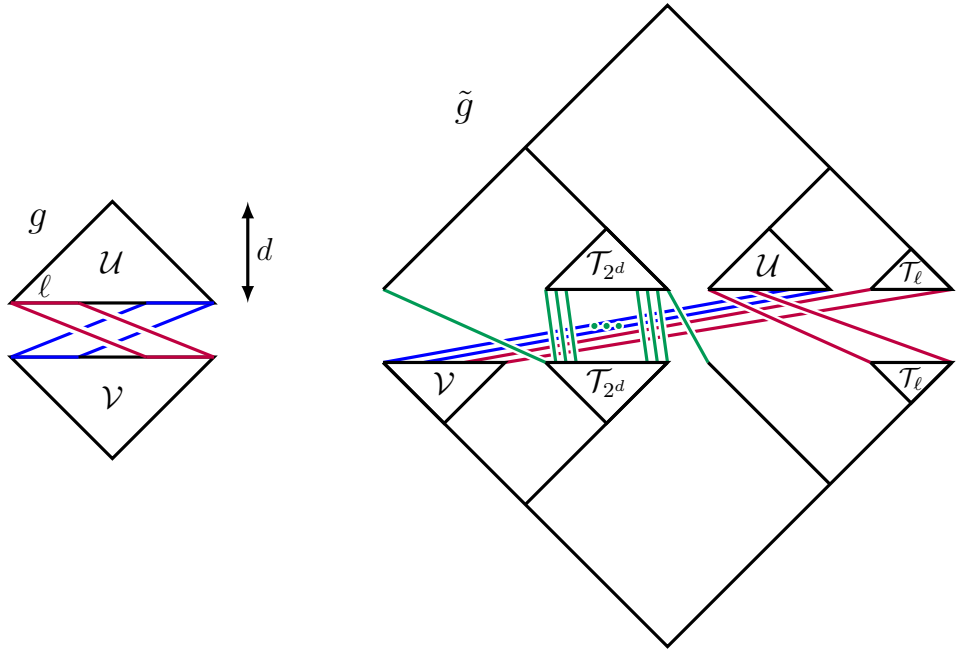
\captionof{figure}{Alternate iterative construction, optimising for depth.}
    \label{fig:TlowCalegari}
\end{center}
From each (torsion) element $g\in T$, we associate another (torsion) element $\tilde g\in T$, as described in \cref*{fig:TlowCalegari}. Here $\Tc_n$ denotes any tree with $n$ leaves and depth $\lceil\log_2(n)\rceil$. We control the different parameters of the new element $\tilde g$:
\begin{itemize}[leftmargin=7mm]
    \item The depths of $\Tc_{2^d}$ and $\Tc_\ell$ are smaller than the depths of $\Uc$ and $\Vc$, so we get an equality $d(\tilde g) = d(g)+2$.
    \item The element $\tilde g$ is still torsion, and its order satisfies $\ord(\tilde g)\ge 2^{d(g)}\cdot\ord(g)$. %\coco{Maybe add a picture for this. I mean this will be the 10th example, so maybe not.}
\end{itemize}
If we start with the element $b_0\colon 0w\leftrightarrow 1w$ (with $d(b_0)=1$ and $\ord(b_0)=2$) and iterate the construction, we get a sequence $(b_k)_{k\ge 0}$ with depth $d(b_k)=2k+1$. Moreover, $\ord(b_{k+1})\ge 2^{d(b_k)}\cdot \ord(b_k)$ inductively gives 
\[ \ord(b_k)\ge 2^{k^2+1} \asymp m(b_k)^{\frac14\log_2 m(b_k)}. \]
This proves that $\ord(g)$ (hence the denominator of $\rot(g)$) does not satisfy any polynomial inequality in terms of the \say{height} $m(g)$.
%\begin{center}
%    \import{Pictures/}{tikz_T_cubic.tex}
%    \captionof{figure}{Elements of $T$ with $n$ leaves and period $n^3$.} \label{fig:Tcubic}
%\end{center}

\subsection{Thompson group \texorpdfstring{$V$}{V}}
Now we focus on Thompson group $V$ both in the case $D=\mathbb{N}$ and in a more general case. We must assume that $D$ is not too sparse.
\begin{thm} \label{thm:period_V} $\,$
    \begin{enumerate}[leftmargin=7mm, label={\normalfont(\alph*)}]
        \item If $D\subseteq D_k$ satisfies a\say{logarithmic Bertrand postulate}, i.e., $\exists C\ge 1$, $\forall M\in \N$, $\exists m\in D$ such that $M\le m\le M^C$, then $p_V^D(n)\asymp \exp(n)$.
        \item $p_V(n)\asymp \exp(n^2)$.
    \end{enumerate}
\end{thm}

\begin{proof}
    (a) For the upper bound, let $g\in V$ have $\ord(g)\in D\subseteq D_k$, so that $\ord(g)=p_1^{e_1}\dots p_k^{e_k}$. Let $\{l_1,\dots, l_m\}$ be the set of labels we get at the end of the algorithm described in \S\ref{subsection:algorithm}. We have that $p_i^{e_i}\mid \lcm(l_1,\ldots,l_m)$, hence $$p_i^{e_i}\leq \max\limits_{i=1}^m l_i\leq \sum_{i=1}^m l_i.$$ Thus, if $\norm{g}_V\leq n$, the order satisfies $\ord(g)\leq (\sum\limits_{i=1}^m  l_i)^k\preceq \exp(n)$ by \cref{prop:components&labels}.

Let us now treat the lower bound. By assumption there exists $m\in D $ such that $(C \exp(Cn))^{\frac{1}{C}}\leq m\leq C \exp(Cn)$. We will show that, for every $m\in\N$, there exists $g\in V$ with length $\asymp\log(m)$ such that $\ord(g)=m$. This implies that $p^D_V(n)\geq (C \exp(Cn))^{\frac{1}{C}}\asymp \exp(n)$. 

Let $m\in [2^{n-1}, 2^n -1]$, and write it as $a_0+a_1 \cdot 2+\ldots +a_{n-2}\cdot 2^{n-2}+2^{n-1}$ for some $a_i\in \{0,1\}$. Let $\{i\mid a_i=1\} = \{i_1,i_2,\ldots,i_d\}$ with $i_1<\ldots<i_d$. Consider the element $c_m=(\Lc_n\wedge\Rc_d,\phi_m ,\Rc_n\wedge\Rc_d) $ where 
\[ \phi_m\colon\left\{ \begin{array}{rcll}
1 & \mapsto & n, & \\
i+2 & \mapsto & n+1-(i+2) & \text{ if }0\leq i\leq n-2 \text{ and } a_i=0,\\
i_j +2 & \mapsto & n+j, & \\
n+j & \mapsto & n+1-(i_j +2). & 
\end{array}\right.\]
The element has $n+d$ leaves with $d\leq n-1$, and $\phi_m=\alpha \beta$, where $\alpha\in~\Sym(n+d)$ satisfies $\LDS(\alpha)\le 3$ and $\beta\in\Sym(n+d)$ is represented by a non-reduced diagram of $g_n@0=g_2g_{n+1}$ (with $g_n$ defined in \cref{prop:flip_length}). Hence $\norm{c_m}_V\preceq n$ by \cref{thm:upper_bound_LDS} and \cref{prop:flip_length}. The example with $m=45$ is depicted in \cref*{fig:lower_bd_V_D}, and \cref*{fig:lower_bd_V_D_algorithm} describes the algorithm.

\begin{center}
    \begin{tikzpicture}[scale=.75]
    \begin{scope}[very thick]
        \draw[line width=4pt,white] (7,2) -- (5,0);
        \draw[ForestGreen] (7,2) -- (5,0);
        \draw[line width=4pt,white] (8,2) -- (3,0);
        \draw[ForestGreen] (8,2) -- (3,0);
        \draw[line width=4pt,white] (9,2) -- (2,0);
        \draw[ForestGreen] (9,2) -- (2,0);

        \draw[line width=4pt,white] (2,2) -- (7,0);
        \draw[blue] (2,2) -- (7,0);
        \draw[line width=4pt,white] (4,2) -- (8,0);
        \draw[blue] (4,2) -- (8,0);
        \draw[line width=4pt,white] (5,2) -- (9,0);
        \draw[blue] (5,2) -- (9,0);
        
        \draw[line width=4pt,white] (1,2) -- (6,0);
        \draw[orange] (1,2) -- (6,0);
        \draw[line width=4pt,white] (3,2) -- (4,0);
        \draw[purple] (3,2) -- (4,0);
        \draw[line width=4pt,white] (6,2) -- (1,0);
        \draw[purple] (6,2) -- (1,0);        
    \end{scope}
    
    \begin{scope}[very thick, shift={(0,2)}]
        \draw (1,0) -- (5,4);
        \draw (2,0) -- (1.5,.5);
        \draw (3,0) -- (2,1);
        \draw (4,0) -- (2.5,1.5);
        \draw (5,0) -- (3,2);
        \draw (6,0) -- (3.5,2.5);
        
        \draw (7,0) -- (8,1);
        \draw (8,0) -- (8.5,.5);
        \draw (9,0) -- (5,4);
    \end{scope}
    
    \begin{scope}[very thick, shift={(0,0)}, yscale=-1]
        \draw (1,0) -- (5,4);
        \draw (2,0) -- (4,2);
        \draw (3,0) -- (4.5,1.5);
        \draw (4,0) -- (5,1);
        \draw (5,0) -- (5.5,.5);
        \draw (6,0) -- (3.5,2.5);

        \draw (7,0) -- (8,1);
        \draw (8,0) -- (8.5,.5);
        \draw (9,0) -- (5,4);
    \end{scope}

    \begin{scope}[very thick, shift={(10,0)}, xscale=.7]
        \draw[line width=4pt,white] (7,2) -- (2,1);
        \draw[ForestGreen] (7,2) -- (2,1);
        \draw[line width=4pt,white] (8,2) -- (4,1);
        \draw[ForestGreen] (8,2) -- (4,1);
        \draw[line width=4pt,white] (9,2) -- (5,1);
        \draw[ForestGreen] (9,2) -- (5,1);

        \draw[line width=4pt,white] (2,2) -- (7,1);
        \draw[blue] (2,2) -- (7,1);
        \draw[line width=4pt,white] (4,2) -- (8,1);
        \draw[blue] (4,2) -- (8,1);
        \draw[line width=4pt,white] (5,2) -- (9,1);
        \draw[blue] (5,2) -- (9,1);
        
        \draw[line width=4pt,white] (1,2) -- (1,1);
        \draw[orange] (1,2) -- (1,1);
        \draw[line width=4pt,white] (3,2) -- (3,1);
        \draw[purple] (3,2) -- (3,1);
        \draw[line width=4pt,white] (6,2) -- (6,1);
        \draw[purple] (6,2) -- (6,1);

        \draw (1,0) -- (6,1);
        \draw (2,0) -- (5,1);
        \draw (3,0) -- (4,1);
        \draw (4,0) -- (3,1);
        \draw (5,0) -- (2,1);
        \draw (6,0) -- (1,1);
        
        \draw[gray] (7,0) -- (7,1);
        \draw[gray] (8,0) -- (8,1);
        \draw[gray] (9,0) -- (9,1);

        \foreach \x in {1,...,9}{
        \draw[fill=black] (\x,0) circle (1.28pt and .9pt);
        \draw[fill=black] (\x,1) circle (1.28pt and .9pt);
        \draw[fill=black] (\x,2) circle (1.28pt and .9pt);}
        
        \node at (10,1.5) {$\alpha$};
        \node at (10,0.5) {$\beta$};
    \end{scope}
\end{tikzpicture}
    \captionsetup{font=small, margin=7mm}
    
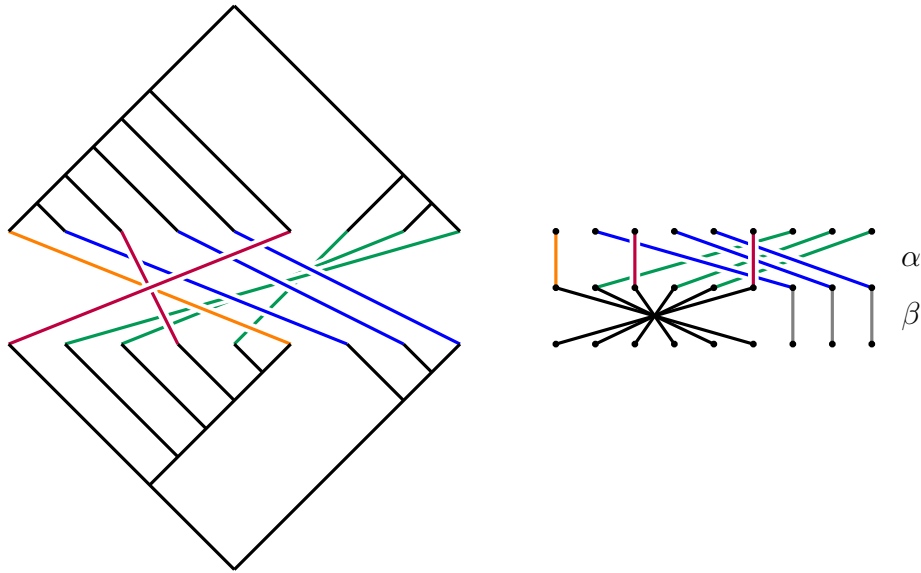
\captionof{figure}{Element for $m=({\color{blue}1}\cdot 2^0+{\color{purple}0}\cdot 2^1+{\color{blue}1}\cdot 2^2+{\color{blue}1}\cdot 2^3+{\color{purple}0}\cdot 2^4)+2^5$, and the associated decomposition $\phi_m=\alpha\beta$.}
    \label{fig:lower_bd_V_D}
\end{center}
After $d$ iterations of Rule~I, we get the tree on the left of \cref*{fig:lower_bd_V_D_algorithm} with $k=n$ and $x=1$. Next we repeat Rule~I for $n-1$ more times. Overall, every $a_i$ contributes $a_i 2^i$ to the order, while the label in the middle contributes $2^{n-1}$, giving a total of $m$.
\begin{center}
    \begin{tikzpicture}[scale=1]
    \begin{scope}[very thick, shift={(0,7)}, yscale=-1]
        \draw (1,0) -- (3.5,2.5);
        \draw (2,0) -- (4,2);
        %\draw (3,0) -- (4.5,1.5);
        \draw (4,0) -- (5,1);
        \draw (5,0) -- (5.5,.5);
        \draw (6,0) -- (3.5,2.5);
        
        \node at (1,-.4) {\footnotesize$k$};
        \node at (2,-.4) {\footnotesize$k-1$};
        \node at (5,-.4) {\footnotesize$2$};
        \node at (6,-.4) {\footnotesize$1$};

        \node[circle, draw=red, very thick, fill=white, inner sep=1.5pt] at (5.25,.25) {{\scriptsize$\color{red}a_0$}};
        \node[circle, draw=red, very thick, fill=white, inner sep=1.5pt] at (4.5,.5) {{\scriptsize$\color{red}a_1$}};
        \node[rotate=17] at (3.8,.7) {\color{red}$\cdots$};
        \node[circle, draw=red, very thick, fill=white, inner sep=1.5pt] at (3,1) {{\scriptsize$\color{red}a_{k-3}$}};
        \node[circle, draw=red, very thick, fill=white, inner sep=1.5pt] at (2.25,1.25) {{\scriptsize$\color{red}a_{k-2}$}};
    \end{scope}

    \begin{scope}[very thick, shift={(0,0)}]
        \draw (1,0) -- (3.5,2.5);
        \draw (2,0) -- (1.5,.5);
        \draw (3,0) -- (2,1);
        %\draw (4,0) -- (2.5,1.5);
        \draw (5,0) -- (3,2);
        \draw (6,0) -- (3.5,2.5);

        \node at (1,-.4) {\footnotesize$1$};
        \node at (2,-.4) {\footnotesize$2$};
        \node at (5,-.4) {\footnotesize$k-1$};
        \node at (6,-.4) {\footnotesize$k$};

        \draw[very thick] (3.5,2.5) -- (3.5,4.5);
        \node[circle, draw=red, very thick, fill=white, inner sep=2pt] at (3.5,3.5) {{$\color{red}x$}};
        \node[rotate=17] at (3.2,.7) {$\cdots$};
    \end{scope}
    
    \draw[very thick, -latex] (6,3.5) -- (8.5,3.5);
    \node at (7.25,4) {Rule~I};
    
    \begin{scope}[very thick, shift={(7.5,7)}, yscale=-1]
        \draw (1,0) -- (3,2);
        %\draw (2,0) -- (3.5,1.5);
        \draw (3,0) -- (4,1);
        \draw (4,0) -- (4.5,.5);
        \draw (5,0) -- (3,2);

        \node at (1,-.4) {\footnotesize$k-1$};
        \node at (4,-.4) {\footnotesize$2$};
        \node at (5,-.4) {\footnotesize$1$};

        \node[circle, draw=red, very thick, fill=white, inner sep=1.5pt] at (4.25,.25) {{\scriptsize$\color{red}a_0$}};
        \node[circle, draw=red, very thick, fill=white, inner sep=1.5pt] at (3.5,.5) {{\scriptsize$\color{red}a_1$}};
        \node[rotate=17] at (2.8,.7) {\color{red}$\cdots$};
        \node[circle, draw=red, very thick, fill=white, inner sep=1.5pt] at (2,1) {{\scriptsize$\color{red}a_{k-3}$}};
    \end{scope}

    \begin{scope}[very thick, shift={(7.5,0)}]
        \draw (1,0) -- (3,2);
        \draw (2,0) -- (1.5,.5);
        \draw (3,0) -- (2,1);
        %\draw (4,0) -- (2.5,1.5);
        \draw (5,0) -- (3,2);

        \node at (1,-.4) {\footnotesize$1$};
        \node at (2,-.4) {\footnotesize$2$};
        \node at (5,-.4) {\footnotesize$k-1$};

        \draw[very thick] (3,2) -- (3,5);
        \node[rectangle, rounded corners=3, draw=red, very thick, fill=white, inner sep=2pt] at (3,3.5) {{$\color{red}2x+a_{k-2}$}};
        \node[rotate=17] at (3.2,.7) {$\cdots$};
    \end{scope}
\end{tikzpicture}
    \captionsetup{font=small}
    
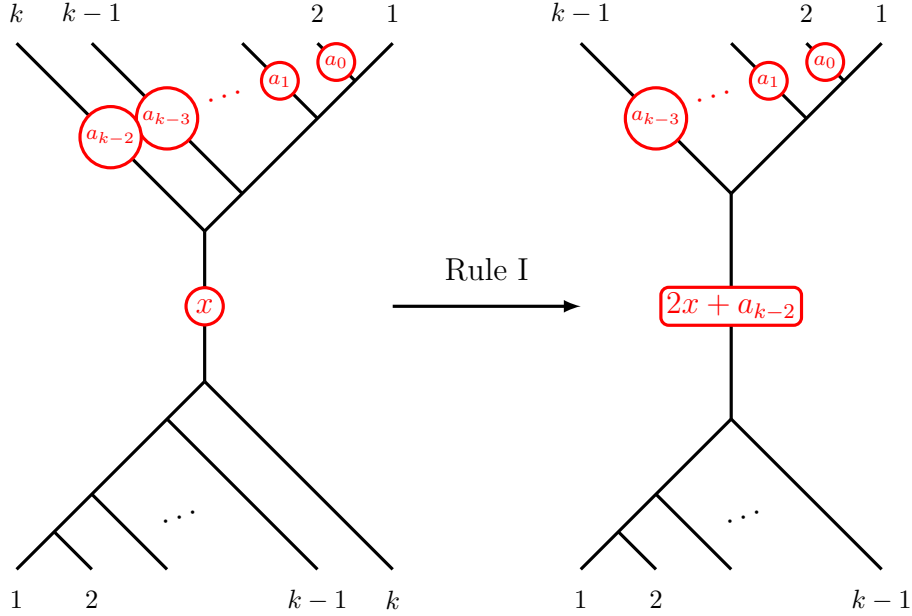
\captionof{figure}{One step of the algorithm applied to $c_m$.}
    \label{fig:lower_bd_V_D_algorithm}
\end{center}

(b) For the upper bound, an element $g\in V$ of length at most $n$ can be represented by a tree pair diagram with at most $Cn$ leaves by \cref{thm:metric_all}. It gives at most $Cn$ connected components at the end of the algorithm in \S\ref{subsection:algorithm}, and the label in each component is at most $2^{Cn-1}$ by \cref{prop:components&labels}. The order of $g$ is at most $(2^{Cn-1})^{Cn}\asymp2^{n^2}$.\\

For the lower bound, we consider $d_n=(\Sc_{2n+1}\wedge\Rc_{n+1},\,\psi_n,\,\Lc_{2n+1}\wedge \Rc_{n+1})$ where $\Sc_{2n+1}$ is a tree with $2n+1$ leaves highlighted in \cref{fig:Vlow} and $\psi_n\in\Sym(3n+2)$ is defined by
\[ \psi_n\colon \left\{ \begin{array}{rcll}
i & \mapsto & 2n-2i+2 & \text{ if }1\le i\le n \\
i & \mapsto & i+n+1 & \text{ if }n+1\le i\le 2n+1 \\
i & \mapsto & 2i-4n-3 & \text{ if }2n+2\le i\le 3n+2\end{array}\right.\]
We decompose $\psi_n=\alpha\beta$ (see the right of \cref{fig:Vlow}) where
\begin{itemize}[leftmargin=7mm]
    \item $\alpha$ is the permutation in a \emph{non-reduced} diagram for $g_n@0=g_{n+1}g_2$ (with $g_n$ defined in \cref{prop:flip_length}), hence $L(\alpha)\le \norm{g_n@0}_V\preceq n$, and
    \item $\beta\in\Sym(3n+2)$ is a riffle shuffle hence $L(\beta)\preceq n$ by \cref{lem:induced_length_riffle}. 
\end{itemize}

\begin{center}
    \begin{tikzpicture}[scale = .65]
    \pgfmathsetmacro{\n}{3}
    \pgfmathsetmacro{\np}{\n+1}

    \begin{scope}[very thick]
        \foreach \i in {1,...,\np}{
            \draw[ForestGreen] ({2*\n+\i+1},2) -- ({2*\i-1},0);}
        \foreach \i in {1,...,\n}{
            \draw[white, line width=4pt] (\i,2) -- ({2*\n+2-2*\i},0);}
        \foreach \i in {1,...,\np}{
            \draw[white, line width=4pt] ({\n+\i},2) -- ({2*\n+\i+1},0);}
            
        \fill[pink!20, rounded corners=7, shift={(0,2)}] (.4,-.2) -- (4,3.4) -- (7.6,-.2) -- cycle;
        
        \foreach \i in {1,...,\n}{
            \draw[purple] (\i,2) -- ({2*\n+2-2*\i},0);}
        \foreach \i in {1,...,\np}{
            \draw[blue] ({\n+\i},2) -- ({2*\n+\i+1},0);}
    \end{scope}

    \begin{scope}[very thick, shift={(0,2)}]
        \node at (2.8,3.3) {\large$d_3$};

        \draw (1,0) -- (6,5);
        \draw (2,0) -- (4,2);
        \draw (3,0) -- (4,1);
        \draw (4,0) -- (3.5,0.5);
        \draw (5,0) -- (3.5,1.5);
        \draw (6,0) -- (3.5,2.5);
        \draw (7,0) -- (4,3);
        
        \draw (8,0) -- (9.5,1.5);
        \draw (9,0) -- (10,1);
        \draw (10,0) -- (10.5,.5);
        \draw (11,0) -- (6,5);
    \end{scope}

    \begin{scope}[very thick, yscale=-1]
        \draw (1,0) -- (6,5);
        \draw (2,0) -- (1.5,.5);
        \draw (3,0) -- (2,1);
        \draw (4,0) -- (2.5,1.5);
        \draw (5,0) -- (3,2);
        \draw (6,0) -- (3.5,2.5);
        \draw (7,0) -- (4,3);
        
        \draw (8,0) -- (9.5,1.5);
        \draw (9,0) -- (10,1);
        \draw (10,0) -- (10.5,.5);
        \draw (11,0) -- (6,5);
    \end{scope}

    \begin{scope}[very thick, shift={(12,0)}, xscale=.7]
        \foreach \i in {1,...,\np}{
            \draw[gray] ({2*\n+\i+1},2) -- ({2*\n+\i+1},1);
            \draw[ForestGreen] ({2*\n+\i+1},1) -- ({2*\i-1},0);}
        \foreach \i in {1,...,\n}{
            \draw[white, line width=4pt] ({\n+1-\i},1) -- ({2*\n+2-2*\i},0);
            \draw[purple] (\i,2) -- ({\n+1-\i},1);
            \draw[purple] ({\n+1-\i},1) -- ({2*\n+2-2*\i},0);}
        \foreach \i in {1,...,\np}{
            \draw[white, line width=3pt] ({\n+\i},1) -- ({2*\n+\i+1},0);
            \draw[gray] ({\n+\i},2) -- ({\n+\i},1);
            \draw[blue] ({\n+\i},1) -- ({2*\n+\i+1},0);}
        \foreach \x in {1,...,11}{
            \draw[fill=black] (\x,0) circle (1.3pt and .9pt);
            \draw[fill=black] (\x,1) circle (1.3pt and .9pt);
            \draw[fill=black] (\x,2) circle (1.3pt and .9pt);}
        \node at (12,1.5) {$\alpha$};
        \node at (12,0.5) {$\beta$};
    \end{scope}
\end{tikzpicture}
    \captionsetup{font=small, margin=20mm}
    
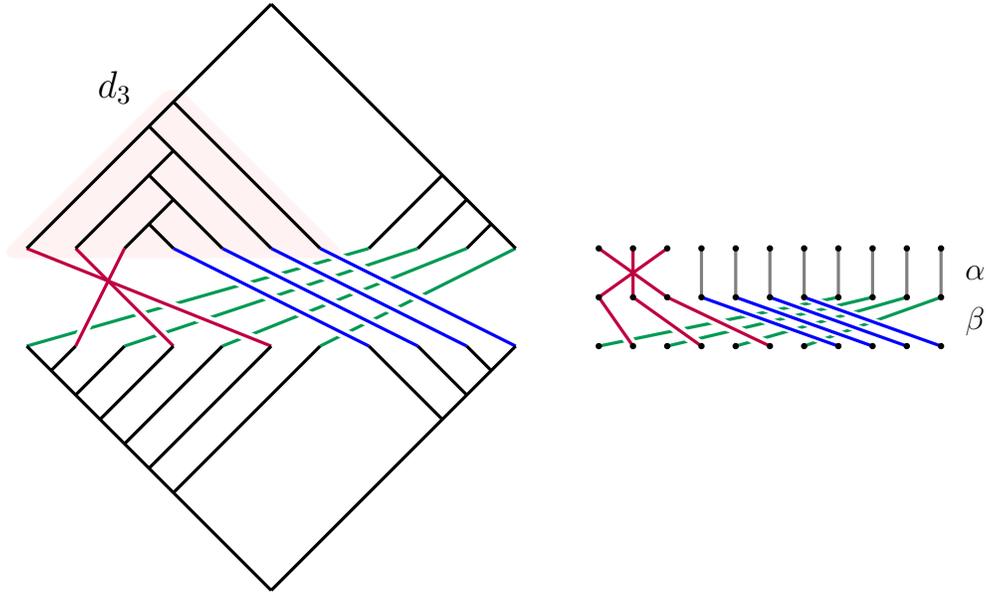
\captionof{figure}{Elements of the sequence $d_n$ realizing the lower bound, and the associated decomposition $\psi_n=\alpha\beta$.} \label{fig:Vlow}
\end{center}

Using \cref{lem:induced_length_composition}, we conclude that
\[ \norm{d_n}_V \le L(\psi_n)+O(n) \le L(\alpha)+L(\beta)+O(n) = O(n). \]
Finally, we need to compute the order of $d_n$. Recall that $d_n$ has $3n+2$ leaves. We first apply the algorithm (Rule~I) once to split the diagram into two parts, and then apply Rule~I another $n$ times to the right tree. This produces a diagram with $2n+1$ leaves, whose labels are exactly those shown in \cref{fig:algorithm-V} (first diagram on the left with $k=n$ and $x=1$). We then apply Rule~I twice, as illustrated in \cref{fig:algorithm-V}. This yields two connected components: one is a $(2n-1)$-leaves version of the $(2n+1)$-leaves diagram, and the other is a cycle, both having central label $2$.
\begin{center}
    \begin{tikzpicture}[scale=.555]
    \clip (0.7,-.7) rectangle (25,8.7);
    
    \begin{scope}[very thick, shift={(0,0)}]
        \draw (1,0) -- (4,3);
        \draw (2,0) -- (4,2);
        \draw (3,0) -- (4,1);
        \draw (4,0) -- (3.5,0.5);
        \draw (5,0) -- (3.5,1.5);
        \draw (6,0) -- (3.5,2.5);
        \draw (7,0) -- (4,3);

        \node at (1,-.4) {\scriptsize$2k$};
        \node at (2,-.4) {\scriptsize$\ldots$};
        \node at (3,-.4) {\scriptsize$2$};
        \node at (4,-.4) {\scriptsize$1$};
        \node at (5,-.4) {\scriptsize$3$};
        \node at (5.75,-.4) {\scriptsize$\cdots$};
        \node at (7,-.4) {\scriptsize$2k+1$};

        \node[circle, draw=red, very thick, fill=white, inner sep=1pt] at (3.75,.25) {{\scriptsize$\color{red}1$}};
        \node[circle, draw=red, very thick, fill=white, inner sep=1pt] at (4.5,.5) {{\scriptsize$\color{red}1$}};
        \node[circle, draw=red, very thick, fill=white, inner sep=1pt] at (5.25,.75) {{\scriptsize$\color{red}1$}};
        \node[circle, draw=red, very thick, fill=white, inner sep=1pt] at (6,1) {{\scriptsize$\color{red}1$}};

        \draw (4,3) -- (4,5);
        \node[circle, draw=red, very thick, fill=white, inner sep=1.5pt] at (4,4) {{\small$\color{red}x$}};
    \end{scope}

    \begin{scope}[very thick, shift={(0,8)}, yscale=-1]
        \draw (1,0) -- (4,3);
        \draw (2,0) -- (1.5,.5);
        \draw (3,0) -- (2,1);
        \draw (4,0) -- (2.5,1.5);
        \draw (5,0) -- (3,2);
        \draw (6,0) -- (3.5,2.5);
        \draw (7,0) -- (4,3);

        \node at (1,-.4) {\scriptsize$1$};
        \node at (2,-.4) {\scriptsize$2$};
        \node at (7,-.4) {\scriptsize$2k+1$};
    \end{scope}

    %%%%%%%%%%%%
    \node at (7.75,4.5) {\scriptsize Rule~I};
    \draw[very thick, -latex] (6,4) -- (9.5,4);
    %%%%%%%%%%%%
    
    \begin{scope}[very thick, shift={(8,0)}]
        \draw (1,0) -- (3.5,2.5);
        \draw (2,0) -- (4,2);
        \draw (3,0) -- (4,1);
        \draw (4,0) -- (3.5,0.5);
        \draw (5,0) -- (3.5,1.5);
        \draw (6,0) -- (3.5,2.5);
        \draw (7.5,0) -- (7.5,8);

        \node at (1,-.4) {\scriptsize$2k$};
        \node at (2,-.4) {\scriptsize$\ldots$};
        \node at (3,-.4) {\scriptsize$2$};
        \node at (4,-.4) {\scriptsize$1$};
        \node at (4.75,-.4) {\scriptsize$\ldots$};
        \node at (6,-.4) {\scriptsize$2k-1$};
        \node at (7.5,-.4) {\scriptsize$2k+1$};

        \node[circle, draw=red, very thick, fill=white, inner sep=1pt] at (3.75,.25) {{\scriptsize$\color{red}1$}};
        \node[circle, draw=red, very thick, fill=white, inner sep=1pt] at (4.5,.5) {{\scriptsize$\color{red}1$}};
        \node[circle, draw=red, very thick, fill=white, inner sep=1pt] at (5.25,.75) {{\scriptsize$\color{red}1$}};

        \draw (3.5,2.5) -- (3.5,5.5);
        \node[circle, draw=red, very thick, fill=white, inner sep=1.5pt] at (3.5,4) {{\small$\color{red}x$}};
        \node[rectangle, rounded corners=3, draw=red, very thick, fill=white, inner sep=2pt] at (7.5,4) {{\small$\color{red}x+1$}};
    \end{scope}

    \begin{scope}[very thick, shift={(8,8)}, yscale=-1]
        \draw (1,0) -- (3.5,2.5);
        \draw (2,0) -- (1.5,.5);
        \draw (3,0) -- (2,1);
        \draw (4,0) -- (2.5,1.5);
        \draw (5,0) -- (3,2);
        \draw (6,0) -- (3.5,2.5);

        \node at (1,-.4) {\scriptsize$1$};
        \node at (2,-.4) {\scriptsize$2$};
        \node at (6,-.4) {\scriptsize$2k$};
        \node at (7.5,-.4) {\scriptsize$2k+1$};
    \end{scope}

    %%%%%%%%%%%%
    \node at (18.25,4.5) {\scriptsize Rule~I};
    \draw[very thick, -latex] (17.5,4) -- (19,4);
    %%%%%%%%%%%%

    \begin{scope}[very thick, shift={(17.5,0)}]
        \draw (1,0) -- (3,2);
        \draw (2,0) -- (3,1);
        \draw (3,0) -- (2.5,0.5);
        \draw (4,0) -- (2.5,1.5);
        \draw (5,0) -- (3,2);
        \draw (6.5,0) -- (6.5,8);

        \node at (1,-.4) {\scriptsize$2k-2$};
        %\node at (2,-.4) {\scriptsize$\ldots$};
        \node at (2,-.4) {\scriptsize$2$};
        \node at (3,-.4) {\scriptsize$1$};
        \node at (3.75,-.4) {\scriptsize$\ldots$};
        \node at (5,-.4) {\scriptsize$2k-1$};
        \node at (6.5,-.4) {\scriptsize$2k+1$};

        \node[circle, draw=red, very thick, fill=white, inner sep=1pt] at (2.75,.25) {{\scriptsize$\color{red}1$}};
        \node[circle, draw=red, very thick, fill=white, inner sep=1pt] at (3.5,.5) {{\scriptsize$\color{red}1$}};
        \node[circle, draw=red, very thick, fill=white, inner sep=1pt] at (4.25,.75) {{\scriptsize$\color{red}1$}};

        \draw (3,2) -- (3,6);
        \node[rectangle, rounded corners=3, draw=red, very thick, fill=white, inner sep=2pt] at (3,4) {{\small$\color{red}2x$}};
        \node[rectangle, rounded corners=3, draw=red, very thick, fill=white, inner sep=2pt] at (6.5,4) {{\small$\color{red}x+1$}};
    \end{scope}

    \begin{scope}[very thick, shift={(17.5,8)}, yscale=-1]
        \draw (1,0) -- (3,2);
        \draw (2,0) -- (1.5,.5);
        \draw (3,0) -- (2,1);
        \draw (4,0) -- (2.5,1.5);
        \draw (5,0) -- (3,2);

        \node at (1,-.4) {\scriptsize$1$};
        \node at (2,-.4) {\scriptsize$2$};
        \node at (5,-.4) {\scriptsize$2k-1$};
        \node at (6.5,-.4) {\scriptsize$2k+1$};
    \end{scope}
\end{tikzpicture}
    \captionsetup{font=small, margin=10mm}
    
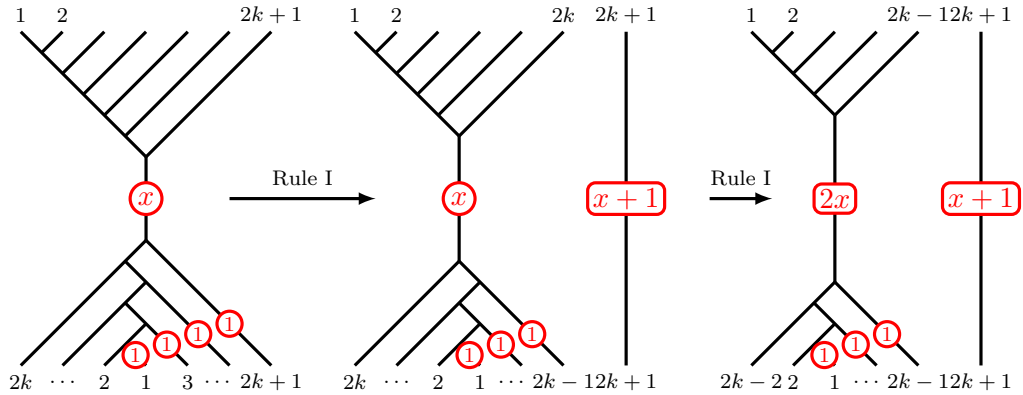
\captionof{figure}{Two consecutive steps of the algorithm when applied to $d_n$.}
    \label{fig:algorithm-V}
\end{center}
By iterating the algorithm $n$ more times on the left tree, we obtain $n+1$ cycles with labels
$2^{0}+1, 2^{1}+1,\ldots, 2^n+1$.
Therefore, the order of $d_n$ is
$$\lcm(2^0+1, 2^1+1, \ldots, 2^n+1) \asymp \exp(n^2), $$
see the estimate on page 410 of \cite{Unary_CF}.
\end{proof}

\medskip

\begin{rem} \label{rem:clone}
    Another closely related family of groups with co-context-free Word Problem are the groups $V_{(G,\theta)}$ with $G$ finite \cite{ThompsonV_clones}, see also \cite{cloning_explained} for a more detailed construction. It remains unclear whether these groups embed in Thompson group $V$, hence whether they could be counter-examples to Lehnert's conjecture. It is therefore natural to wonder about their period growth. Recall that these groups are instances of cloning systems as defined in \cite{cloning}, for $G_n=G^n \rtimes \Sym(n)$ with the obvious map $\rho_n\colon G_n\to\Sym(n)$. Using \cite[Lemma 3.12]{cloning}, we get a short exact sequence
    \begin{center}\begin{tikzcd}
        1 \arrow[r] & \bigoplus_{i=1}^\infty G \arrow[r] & V_{(G,\theta)} \arrow[r, "\mathcal T(\rho_*)"] &  V \arrow[r] & 1.
    \end{tikzcd}\end{center}
    Moreover, this sequences splits since the maps $\rho_n$ split, using Observation 2.28 and 2.29 from the same article. When $G$ is finite, by Lemma~\ref{lem:subgroup} and \ref{item-semi:1}, it easily follows that $p_{V_{(G,\theta)}}^D(n)\asymp p_{V}^D(n)$ for all $D\subset\N$.
    % Let $\pi=\mathcal T(\rho_*)$ be the projection $V_{(G,\theta)}\to V$. Pick a generating set $S$ of $V_{(G,\theta)}$, and consider the generating set $T=\pi(S)$ of $V$.
    % For $g\in V_{(G,\theta)}$, we have $\norm{g}_S\ge \norm{\pi(g)}_T$ and $\ord(g)= d\ord(\pi(g))$ for some divisor $d$ of $\abs G$. This shows that
    %\[ p_{V_{(G,\theta)}}^Q(n)\preceq p_V^Q(n). \]
    % For the reverse inequality, start with $h\in V$ of length $\norm{h}_T\le n$ such that $\ord(h)=p_V^Q(n)$ (in particular, the order of $h$ is $Q$-smooth). There exists $g\in V_{(G,\theta)}$ such that $\pi(g)=h$ and $\norm{g}_S=\norm{h}_T$. We have $\ord(g)=d\ord(h)$ for some divisor $d$ of $\abs G$. The element $g^d$ satisfies $\norm{g^d}_S\le d\norm g_S\le dn$ and $\ord(g^d)=\ord(h)=p_V^Q(n)$. This proves that
    %\[ p_{V_{(G,\theta)}}^Q(n)\succeq p_V^Q(n). \]
\end{rem}
\section{A \textbf{CF-TR} group with large period growth} \label{sec:CFTR}

In \cite{CFTR}, a new family of groups was introduced: transition groups of context-free graphs, also abbreviated \textbf{CF-TR} groups. We refer to \cite{CFTR} for the relevant definitions. Moreover, upper bounds for the period growth of \textbf{CF-TR} groups are proved. Looking into the details, one can extract the following results:
\begin{thm} \label{thm:CFTR_period}
	Consider $\Gamma$ a context-free graph and let $G\in {\normalfont\textbf{CF-TR}}$ be its transition group. Define $\bar\beta_\Gamma(n)=\sup_{v\in\Gamma}\abs{B_\Gamma(v,n)}$. Then
	\begin{enumerate}[leftmargin=7mm, label={\normalfont(\alph*)}]
		\item $p_G^{D_k}(n)\preceq \bigl(n\cdot \bar\beta_\Gamma(n)\bigr)^k\preceq \exp(n)$, and \hfill {\normalfont\cite[Lemma 5.3]{CFTR}}
		\item $p_G(n)\preceq \bar\beta_\Gamma(n)^n \preceq \exp(n^2)$. \hfill {\normalfont\cite[Proposition 5.4]{CFTR}}
	\end{enumerate}
\end{thm}
An important question is whether Thompson group $V$ is \textbf{CF-TR}. Combining the previous theorem with \cref{thm:period_V}, we get the following corollary:
\begin{cor} \label{cor:CF_graph_for_V_questionmark}
	If $\Gamma$ is a context-free graph whose transition group is $V$, then the volume growth $\bar\beta_\Gamma(n)$ is exponential.
\end{cor}
% One can get to the same conclusion using that $V$ is not virtually indicable. Indeed, context-free graphs have regular geodesic normal form, so they growth is either polynomial or exponential. Moreover, the growth is polynomial if and only if the graph is a bouquet of combs. If you take a sequence of points going down a \say{tooth} of the comb of the last layer, take the Chabauty limit of graphs you see around them, you get a quasi-transitive graph quasi-isometric to $\Z$ which is also a Schreier graph of the transition group, hence the group has a finite-index subgroup which surjects onto $\Z$.
Any improvement on \cref*{thm:CFTR_period} would directly imply that $V\notin\textbf{CF-TR}$. In the remainder of this section, we provide a new family of \textbf{CF-TR} groups for which \cref*{thm:CFTR_period}(a) is sharp. We consider instances of the groups
\[ B(G) = \bigl(\ldots\wr G)\wr G)\wr G)\wr \ldots \bigr)\rtimes\Z, \]
introduced by A.\ Erschler in \cite[\S 4]{Erschler}. The case $B(\Z)$ was also studied by Brin and Navas in \cite{Brin_Navas1, Brin_Navas2}, and is now called the Brin-Navas group. Other groups $B(G)$ are called Brin-Erschler-Navas (or BEN) groups.

\begin{rem}
Formally, we first define $H$ to be the group of finitary automorphisms of the $\abs G$-ary rooted tree whose local sections belong to $G$. Note that $H$ contains an isomorphic subgroup $H@1$ consisting of elements supported only below a given vertex $1$ on the first level. Then $B(G)$ is the ascending HNN-extension $B(G) = H *_{x_0}$
where $x_0\colon H \to H@1$.
\end{rem}

The group $B(C_2)$ can be defined via its action on the ternary tree (or \say{bi-infinite binary tree}). We fix bi-infinite path in the tree, which we call \say{the spine}. This defines a partition of the tree into \say{levels} (or horospheres, see \cref{fig:CFTR_story}). The group $B(C_2)$ is generated by two automorphisms:
\begin{itemize}[leftmargin=7mm]
    \item $x_0$, which translates by one along the spine, and moves rigidly all subtrees hanging below the spine.
    \item $y$, which exchanges rigidly the two subtrees below the intersection of the spine and the $0$th level, and act trivially elsewhere.
\end{itemize}

\begin{center}
   \begin{tikzpicture}[scale=.97]
        % Levels
        \draw[thin, dotted] (-.7,-1) -- (11.5,-1);
        \node at (12.15,-1) {\footnotesize level $-1$};
        \draw[thin, dotted] (0,0) -- (11.5,0);
        \node at (12,0) {\footnotesize level $0$};
        \draw[thin, dotted] (1,1) -- (11.5,1);
        \node at (12,1) {\footnotesize level $1$};

        % Tree
        \draw (0,0) -- (1,1) -- (2,0);
        \draw (3,0) -- (4,1) -- (5,0);
        \draw (6,0) -- (7,1) -- (8,0);
        \draw (9,0) -- (10,1) -- (11,0);
        \draw (1,1) -- (2.5,2) -- (4,1);
        \draw (7,1) -- (8.5,2) -- (10,1);
        \draw (2.5,2) -- (5.5,3) -- (8.5,2);
        \draw (5.5,3) -- (10,4);
        
        \draw (-.7,-1) -- (0,0) -- (.7,-1);
        \draw (1.3,-1) -- (2,0) -- (2.7,-1);

        \draw (-1.1,-2) -- (-.7,-1) -- (-.3,-2);
        \draw (.3,-2) -- (.7,-1) -- (1.1,-2);
        
        \node[circle, fill=black, inner sep=1.5pt, label=below:$0$] at (0,0) {};
        \node[circle, fill=black, inner sep=1.5pt, label=below:$1$] at (2,0) {};
        \node[circle, fill=black, inner sep=1.5pt, label=below:$10$] at (3,0) {};
        \node[circle, fill=black, inner sep=1.5pt, label=below:$11$] at (5,0) {};

        \node[circle, fill=black, inner sep=1.5pt, label=below:$100$] at (6,0) {};
        \node[circle, fill=black, inner sep=1.5pt, label=below:$101$] at (8,0) {};
        \node[circle, fill=black, inner sep=1.5pt, label=below:$110$] at (9,0) {};
        \node[circle, fill=black, inner sep=1.5pt, label=below:$111$] at (11,0) {};

        % action!
        \draw[thick, latex-latex, purple, bend right=20] (-.65,-1) to (.65,-1);
        \node[purple] at (0,-1.4) {$y$};

        \draw[thick, -latex, bend right, blue] (2.5,2) to (1,1);
        \draw[thick, -latex, bend right, blue] (1,1) to (0,0);
        \draw[thick, -latex, bend right, blue] (0,0) to (-.7,-1);
        \node[blue] at (0.15,.9) {$x_0$};
    \end{tikzpicture}
   \captionsetup{font=small}
   
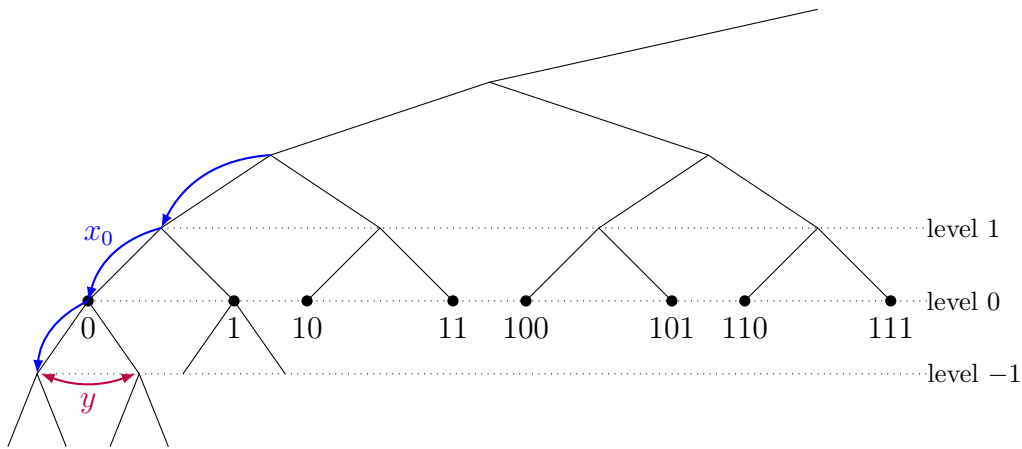
\captionof{figure}{The tree partitioned into levels, and the two automorphisms.}
   \label{fig:CFTR_story}
\end{center}

\bigskip

An interesting observation is that $B(\mathbb{Z})$ embeds into Thompson group $F$. Of course, the group $B(C_2)$ contains non-cyclic finite groups, hence it doesn't embed in $F$, or even in $T$. We prove the next best thing:
\medskip
\begin{prop}
    The group $B(C_2)$ embeds in Thompson group $V$.
\end{prop}
\begin{proof} We have
    \[ B(C_2)= \bigl(\ldots\wr C_2)\wr C_2)\wr C_2)\wr \ldots \bigr)\rtimes\Z \simeq \la x_0,y\ra \le V.  \]
    Alternatively, this is very clear once we know that $\mathrm{QAut}(\mathcal T_{2,c})\into V$ (see \cite{almost_automorphism} for a precise definition and statement).
\end{proof}

\begin{prop} \label{prop:BEN_has_large_period_growth}
    The period growth of $B(C_2)$ satisfies
    \[ \begin{cases} p_{B(C_2)}^D(n) \asymp \exp(n) & \text{if }D\supseteq D_{\{2\}}, \\ p_{B(C_2)}^D(n) \asymp 1 & \text{if }D\not\supseteq D_{\{2\}}. \end{cases} \]
\end{prop}
\begin{proof}
We first observe that $B(C_2)$ fits in a short exact sequence
\[ 1 \longto N \longto B(C_2) \longto \Z \longto 1,\]
where $N$ is the subgroup of level preserving automorphisms of $B(C_2)$. This group is a direct union of copies of the group of finitary automorphisms of the binary rooted tree. In particular, $N$ is a locally finite $2$-group and $\Z$ is torsionfree, hence $B(C_2)$ only admits $2$-torsion.

\begin{minipage}{.75\linewidth}
Adapting the computations from \cref{prop:flip_length} and \cref{exa:adding_machine}, we construct elements with large $2$-torsion. Specifically, the element
\[ g_n@1 = x_0^{-(n-2)}y(x_0y)^{n-2} \]
has order $2^{n-1}$. Combining with the upper bound from $V$ (\cref{lem:sub-over-group}), we conclude that $p_{B(C_2)}^D(n) \asymp\exp(n)$. \qedhere
\end{minipage}
\begin{minipage}{.24\linewidth}
    \centering
    \begin{tikzpicture}[scale = .57]
    \begin{scope}[very thick]
        \node at (.1,1.5) {$g_4@1$};
        
        \draw (0,0) -- (2,2) -- (4,0);
        \draw (1,0) -- (2.5,1.5);
        \draw (2,0) -- (3,1);
        \draw (3,0) -- (3.5,.5);

        \draw (0,-.8) -- (2,-2.8) -- (4,-.8);
        \draw (1,-.8) -- (2.5,-2.3);
        \draw (2,-.8) -- (1.5,-1.3);
        \draw (3,-.8) -- (2,-1.8);

        \draw[blue] (0,0) -- (0,-.8);
        \draw[blue] (1,0) -- (4,-.8);
        \draw[blue] (2,0) -- (3,-.8);
        \draw[blue] (3,0) -- (2,-.8);
        \draw[blue] (4,0) -- (1,-.8);
    \end{scope}
    \end{tikzpicture}
\end{minipage}
\end{proof}

Another interesting parallel is that the Brin-Navas group is \textbf{CF-TR}. Indeed, $F$ admits a faithful action with a context-free Schreier graph \cite{Savchuk}, that is, Thompson $F$ is \textbf{CF-TR}. Moreover the property \textbf{CF-TR} goes to finitely generated subgroups \cite[Proposition 4.17]{CFTR}. This result extends to $B(C_2)$.
\medskip
\begin{thm}
    The group $B(C_2)$ is {\normalfont\textbf{CF-TR}}.
\end{thm}
\begin{proof}
We consider the action of $B(C_2)$ on the vertices of the tree. This action is faithful by definition. Let us try to get a better picture of the associated Schreier graph. Vertices can be parametrized by pairs $(v,m)$ where $m$ is the level of the vertex, and $v\in\Z_{\ge 0}$ is its number on the level (starting from $0$, in binary, see \cref{fig:CFTR_story}). For each non-negative integer $v$, we define its length as its number of digits (in binary), that is,
\[ \ell(v) = \begin{cases} \hspace*{5mm}1 & \text{ if } n=0, \\  \lceil\log_2(n)\rceil & \text{ otherwise.} \end{cases} \]
The action is equivalently defined as $x_0^{\pm1}\cdot (v,m)=(v,m\pm1)$, and 
\[ y\cdot (v,m)= \begin{cases} (v,m) & \text{if }m< \ell(v), \\ (w,m) & \text{if }m\ge \ell(v), \end{cases} \]
where $w$ is obtained by changing the $m$-th digit of $v$ (see \cref{fig:CFTR_grid}).
\begin{center}
    \begin{tikzpicture}[xscale=1.2, yscale=-1.2]
        \clip (-1.7,-1.7) rectangle (8.5,3.55);
        \foreach \x in {0,...,7}{
            \foreach \y in {-2,...,3}{
                \node[circle, fill=black, inner sep=1.5pt] at (\x,\y) {};
                \draw[blue, very thick, -latex] (\x,{\y+.08}) -- (\x,{\y+.95});
            }
            
            \pgfmathsetmacro{\lx}{floor(log2(\x+1))}
            \foreach \y in {-2,...,\lx}{
            \draw[purple, very thick, loop right, -latex] ({\x+.05},\y) to ({\x+.05},\y);
            }
        }
        \draw[purple, very thick, bend left, latex-latex] (0.05,1) to (.95,1);
        \draw[purple, very thick, bend left=20, latex-latex] (0.05,2) to (1.95,2);
        \draw[purple, very thick, bend left=20, latex-latex] (1.05,2) to (2.95,2);
        \draw[purple, very thick, bend left=10, latex-latex] (0.05,3) to (3.95,3);
        \draw[purple, very thick, bend left=10, latex-latex] (1.05,3) to (4.95,3);
        \draw[purple, very thick, bend left=10, latex-latex] (2.05,3) to (5.95,3);
        \draw[purple, very thick, bend left=10, latex-latex] (3.05,3) to (6.95,3);

        \node at (-.9,1) {\footnotesize level $=-1$};
        \node at (-1,0) {\footnotesize level $=0$};
        \node at (-1,-1) {\footnotesize level $=1$};
        \node[circle, fill=Green, inner sep=1.5pt] at (0,0) {};
\end{tikzpicture}
    \captionsetup{font=small, margin=20mm}
    \captionof{figure}{The Schreier graph of $B(C_2)\acts\Z_{\ge 0}\times\Z$}
    \label{fig:CFTR_grid}
\end{center}
This graph can be re-drawn as follows:
\begin{center}
    \begin{tikzpicture}[scale=1.2]
    \begin{scope}[every node/.style={circle, very thick, inner sep=1.8pt}]
        \node[fill=Green] (0l0) at (1,0) {};
        \node[fill=black] (0l1) at (0,-1) {};
    
        \node[fill=black] (0l2) at (-1,-2) {};
        \node[draw=black] (1l1) at (1,-2) {};

        \node[fill=black] (0l3) at (-2,-3) {};
        \node[draw=black] (2l2) at (-.5,-3) {};
        \node[fill=black] (1l2) at (2,-3) {};

        \node[fill=black] (0l4) at (-3,-4) {};
        \node[draw=black] (4l3) at (-1.75,-4) {};
        \node[fill=black] (2l3) at (0,-4) {};
        \node[fill=black] (1l3) at (1.5,-4) {};
        \node[draw=black] (3l2) at (3,-4) {};

        \node[fill=black] (0l5) at (-4,-5) {};
        \node[draw=black] (8l4) at (-2.75,-5) {};
        \node[fill=black] (4l4) at (-1.5,-5) {};
        \node[fill=black] (2l4) at (-.5,-5) {};
        \node[draw=black] (6l3) at (.5,-5) {};
        \node[fill=black] (1l4) at (1,-5) {};
        \node[draw=black] (5l3) at (2,-5) {};
        \node[fill=black] (3l3) at (4,-5) {};

        \node[fill=gray] (1) at (2,0) {};
        \node[fill=gray] (2) at (3,0) {};
        \node (3) at (4,0) {};
    \end{scope}
    
    \begin{scope}[blue, very thick, -latex]
        \draw (0l0) -- (0l1);
        \draw (0l1) -- (0l2);
        \draw (0l2) -- (0l3);
        \draw (0l3) -- (0l4);
        \draw (0l4) -- (0l5);
        \draw (1l1) -- (1l2);
        \draw (2l2) -- (2l3);
        \draw (1l2) -- (1l3);
        \draw (4l3) -- (4l4);
        \draw (2l3) -- (2l4);
        \draw (1l3) -- (1l4);
        \draw (3l2) -- (3l3);

        \draw (3) -- (2);
        \draw (2) -- (1);
        \draw (1) -- (0l0);
    \end{scope}

    \begin{scope}[purple, very thick, latex-latex]
        \draw (0l1) -- (1l1);
        \draw (0l2) -- (2l2);
        \draw (0l3) -- (4l3);
        \draw (0l4) -- (8l4);
        \draw (1l2) -- (3l2);
        \draw (2l3) -- (6l3);
        \draw (1l3) -- (5l3);

        \draw[loop above, -latex] (0l0) to (0l0);
        \draw[loop above, -latex] (1) to (1);
        \draw[loop above, -latex] (2) to (2);
    \end{scope}

    \draw[thin, dotted, out=180, in=180, looseness=6] (3.9,.8) to (3.9,-.8);
    \node at (3.4, -.8) {$\Lambda$};
\end{tikzpicture}
    \captionsetup{font=small, margin=18mm}
    
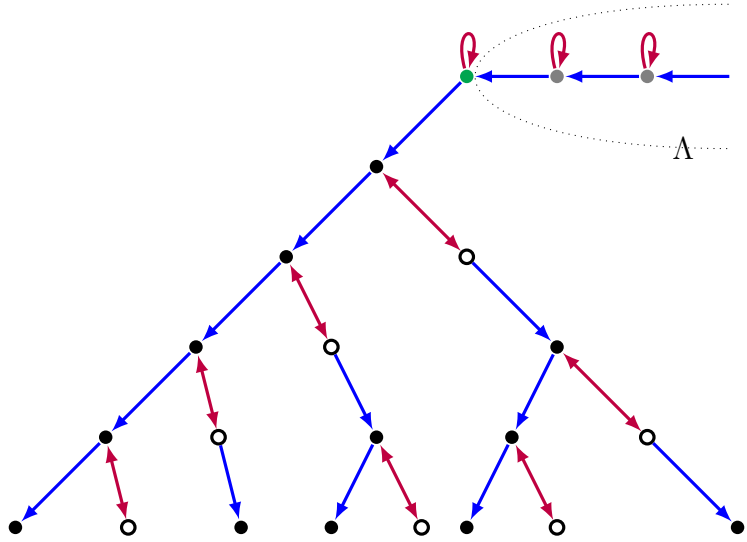
\captionof{figure}{The same Schreier graph, without copies of $\Lambda$ from each vertex of type $\mathbf\circ$. The root is marked in green. The other cone types are marked as $\bullet$, $\circ$ and $\color{gray}\bullet$.}
    \label{fig:CFTR_tree}
\end{center}
From Figure \ref*{fig:CFTR_tree}, it is easy to see that the Schreier graph is context-free in the sense of \cite{Muller_Schupp, CF_pairs}, hence $B(C_2)$ is \textbf{CF-TR}.
\end{proof}

%\medskip
%  Note that another example of \textbf{CF-TR} groups with unbounded torsion are the Houghton groups $H_m$ and similar constructions. In all these cases, we have $p_G^Q(n)\asymp n^{\abs Q}$ for $Q$ finite, and $p_G(n)\asymp \exp\!\big(\sqrt{n\log(n)}\big)$ for $Q=\P$.
% \begin{center}
%     \begin{tikzpicture}[scale=1.5]
%         \foreach \x in {-3,...,4}{
%         \node[circle, fill=black, inner sep=2pt] (P\x) at (\x,0) {};
%         \draw[-latex, thick, blue] ({\x+.05},0) -- ({\x+.95},0);}
        
%         \foreach \x in {-3,-2,-1,2,3,4}{
%             \draw[-latex, thick, magenta, loop, above] (P\x) to (P\x);}
            
%         \draw[-latex, thick, bend right, magenta] (.05,-.05) to (.95,-.05);
%         \draw[-latex, thick, bend right, magenta] (.95,.05) to (.05,.05);
%         \draw[-latex, thick, blue] (-3.7,0) -- (-3.05,0);
%     \end{tikzpicture}
%     \captionsetup{font=small}
%     \captionof{figure}{Context-free graph whose transition group is $H_2=\FSym(\Z)\rtimes\Z$.}
% \end{center}
%\begin{rem}
% This graph is remarkably similar to the Schreier graph of $F$ described in \cite{Savchuk}.
%\end{rem}

\begin{rem}
    The BEN group $B(D_\infty)$ contains  $B(C_2)$. Moreover, one can adapt the previous constructions to prove that it is \textbf{CF-TR} and embeds inside Thompson group $V$. The fact that $B(C_2)\le B(D_\infty)\le V$ and that all finite orders are powers of $2$ implies that $p_{B(D_\infty)}(n)\asymp\exp(n)$. More interestingly A.\ Erschler observed that, for appropriate generating sets, this group is isometric to the Brin-Navas group $B(\mathbb{Z})$ which is torsionfree, proving that period growth is not a quasi-isometry invariant \cite[Theorem 3]{Erschler}.
\end{rem}
\medskip
%\begin{rem}
% The language of ShortLex geodesics for $B$, $G$ and $H$ (with respect to \say{standard} generating sets) seems to be unambigously context-free, hence their (complete) growth series ought to be algebraic. Is it also rational?
%\end{rem}
\section{Further questions}

In all cases that we have encountered, elements of $V$ defined by diagrams with $n$ leaves have linear word length, this begs the following question:

\textbf{Open problem.} Find an \emph{explicit} sequence of element $g_n\in V$, where $g_n$ is defined by a reduced diagram with $n$ leaves, for which $\norm{g_n}\asymp n\log(n)$. 

This question might be surprising since elements satisfying $\norm{g}\asymp n\log(n)$ are exponentially generic (when $n\to\infty$) \cite[Proposition~3.8(2)]{ThompsonV_metric}. Despite this, we are not aware of any explicit family for which $\norm{g_n}$ grows even superlinearly.\footnote{To formalize this question and make a negative answer possible, does there exists a sequence $(g_n)$ and an algorithm taking as input $n$ and printing $g_n$ in $\mathrm{DTIME}\bigl(\exp(o(n))\bigr)$?} This situation is reminiscent of circuit complexity: it is known that most Boolean functions on $n$ variables can only be computed by circuits of size $\asymp\frac{2^n}n$, however it is a notoriously hard open problem to prove that any given function actually requires a circuit of superlinear size. One can also draw parallels with more elementary problems that remain open, for instance finding explicit elements in $\mathrm{GL}_n(\F_2)$ of superlinear length, see \cite[\S7.1]{wigderson2007p}.
\[ *** \]
Comparing \cref{thm:coCF_period} and the different examples studied throughout the article, a natural question is whether \cref*{thm:coCF_period} is sharp:

\textbf{Conjecture.} Let $G$ be a co-context-free group. Then
\begin{enumerate}[leftmargin=8mm, label={\normalfont(\alph*)}]
    \item The period growth satisfies $p_G(n)\preceq \exp(n^2)$. 
    \item If $D\subseteq D_k$ for some integer $k$, then $p_G^D(n)\preceq \exp(n)$. 
\end{enumerate}

Note that counter-examples to either part of the conjecture would disprove Lehnert's conjecture. Given \cref{rem:clone} and \cref{thm:CFTR_period}, the only examples in the literature for which these bounds are not known are FSS groups \cite{Farley_FFS}.

\bigskip

Another motivation was an approach to prove that Grigorchuk group $\mathfrak G$ is not co-context-free. Our results do not conclude since $p_{\mathfrak G}(n)\preceq n^\frac32$ \cite{Bartholdi_Sunic}. However, one realises that elements of small order should be thought as \say{very torsion} while elements of large finite order should be thought of as \say{barely torsion}. This leads to the following question:

\textbf{Question.} Does there exist a co-context-free group $G$ for which $p_G(n)$ is unbounded, but is still bounded above by a polynomial?

Currently, the slowest growing examples we have are Houghton $H_m$ (and similar constructions, eg.\ lampshuffler groups) with $p_{H_m}(n)\asymp \exp(\sqrt{n\log n})$.
\[ *** \]
A family of Thompson-like groups exhibiting mixed behaviors is given by the rearrangement groups of fractals introduced in \cite{belk2019rearrangement}. This family contains Thompson groups $F$ ($p_F(n)=1$) and $V$ ($p_V(n)\asymp\exp(n^2)$) at the two extremes. Most examples studied contain copies of Thompson $T$, with a notable exception for dendrite rearrangement groups \cite{tarocchi2025dendrites}. One could pursue results on period growth in this family by building on the techniques developed in this manuscript together with Tarocchi’s strand-diagrams \cite{tarocchi2023conjugacy}. \

\textbf{Open problem.} Find a criterion determining when a finitely generated rearrangement group has exponential period growth.

This would require a general statement about their metric properties, which might appear difficult; nevertheless, all explicit examples known so far seem to share the same nature as Thompson groups $F$ and $T$ (see \cref{thm:metric_all}).
\[ *** \]
% Remove this is keep old CF-TR
It remains open whether Thompson $T$ and $V$ are contained in the class of \textbf{CF-TR} groups. In particular, whether \textbf{CF-TR} groups coincides with the class of finitely generated subgroups of $V$. It would be interesting to know how fast $p_G^{\P}(n)$ can grow for \textbf{CF-TR} groups.

\textbf{Question.} Does there exist a \textbf{CF-TR} group $G$ such that $p_G^\P(n)\asymp\exp(n)$?

So far, the fastest growth comes from Houghton groups $H_m\in$ \textbf{CF-TR} (see \cref*{fig:Houghton_CFTR}) with $p_{H_2}^{\P}(n)\asymp n$, while $p_V^{\P}(n)\asymp\exp(n)$.
\begin{center}
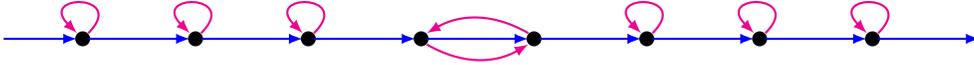

	\begin{tikzpicture}[scale=1.5]
		\foreach \x in {-3,...,4}{
			\node[circle, fill=black, inner sep=2pt] (P\x) at (\x,0) {};
			\draw[-latex, thick, blue] ({\x+.05},0) -- ({\x+.95},0);}
		
		\foreach \x in {-3,-2,-1,2,3,4}{
			\draw[-latex, thick, magenta, loop, above] (P\x) to (P\x);}
		
		\draw[-latex, thick, bend right, magenta] (.05,-.05) to (.95,-.05);
		\draw[-latex, thick, bend right, magenta] (.95,.05) to (.05,.05);
		\draw[-latex, thick, blue] (-3.7,0) -- (-3.05,0);
	\end{tikzpicture}
	\captionsetup{font=small}
	\captionof{figure}{Context-free graph whose transition group is $H_2=\FSym(\Z)\rtimes\Z$.}
	\label{fig:Houghton_CFTR}
\end{center}
\medskip
The class of BEN groups behaves similarly to wreath products. It would be interesting to investigate them further. We hightlight two questions:

\textbf{Question.} Is the growth series of $B(C_2)$ rational/algebraic?

If the growth series is rational, this could be a candidate of group with rational growth series and non $\Z G
$-rational complete growth series. If it is algebraic but non-rational, this would be the first amenable example with these combination of properties after groups of the form $L\wr F_2$.

\textbf{Question.} Is Subgroup Membership decidable in $B(C_2)$ and $B(\Z)$? 

Note that the decidability of Subgroup Membership is open in $F$, treating the Brin-Navas group $B(\Z)$ seems like a natural stepping stone.

\subsection*{Acknowledgment}

We thank Michal Ferov for suggesting \cref{prop:arbitrary_torsion_growth}. The second author would like to thank Giles Gardam and Cosmas Kravaris for conversations and references. The authors were supported by the Swiss NSF grant 200020-200400 during this project. The third and fourth authors were funded by the Swiss Government Excellence Scholarships No.\ 2022.0011 and 2024.0010. % respectively.

\emergencystretch=1em
\AtNextBibliography{\small}
\printbibliography
\Address

\end{document}